\documentclass[11pt]{extarticle}
\usepackage{amsmath,amsthm,amssymb,amsbsy,amsopn,mathtools,comment,bbm}
\usepackage{graphicx,color}
\usepackage[bookmarks=true,bookmarksnumbered=true,bookmarkstype=toc]{hyperref}
\usepackage{mathptmx,bm}
\usepackage{algorithm}
\usepackage{algorithmicx,tabularx,tabu,longtable}
\usepackage{epstopdf}
\usepackage[titletoc]{appendix}
\usepackage{lineno,datetime}

\def\essinf{\mathop{\rm essinf}}
\def\esssup{\mathop{\rm esssup}}

\newtheorem{Def}{Definition}[section]
\newtheorem{theorem}[Def]{Theorem}
\newtheorem{lemma}[Def]{Lemma}
\newtheorem{proposition}[Def]{Proposition}

\newtheorem{remark}[Def]{Remark}

\newtheorem{assumption}[Def]{Assumption}

\usepackage{caption}
\usepackage{subcaption}

\theoremstyle{definition}

\makeatletter
\@addtoreset{equation}{section}
\makeatother

\title{A recursive representation for decoupling time-state dependent jumps from jump-diffusion processes\footnotetext[0]{{\it Stochastics}, doi:10.1080/17442508.2023.2259534}}
\author{\sc Qinjing Qiu\footnote{Email address: qqiu9438@uni.sydney.edu.au. Postal Address: School of Mathematics and Statistics, University of Sydney, NSW 2006, Australia.}
\, and Reiichiro Kawai\footnote{Corresponding author.
Email Address: raykawai@g.ecc.u-tokyo.ac.jp. Postal Address: Graduate School of Arts and Sciences / Mathematics and Informatics Center, The University of Tokyo, Japan.
This work was initiated while RK was based in School of Mathematics and Statistics at the University of Sydney, Australia.}}
\date{}

\begin{document}
\maketitle

\begin{abstract}
\noindent 
We establish a recursive representation that fully decouples jumps from a large class of multivariate inhomogeneous stochastic differential equations with jumps of general time-state dependent unbounded intensity, not of L\'evy-driven type that essentially benefits a lot from independent and stationary increments.
The recursive representation, along with a few related ones, are derived by making use of a jump time of the underlying dynamics as an information relay point in passing the past on to a previous iteration step to fill in the missing information on the unobserved trajectory ahead.
We prove that the proposed recursive representations are convergent exponentially fast in the limit, and can be represented in a similar form to Picard iterates under the probability measure with its jump component suppressed.
On the basis of each iterate, we construct upper and lower bounding functions that are also convergent towards the true solution as the iterations proceed.
We provide numerical results to justify our theoretical findings.

\vspace{0.3em}
\noindent {\it Keywords:} jump-diffusion processes; time-state dependent jump rate; Picard iteration; partial integro-differential equations; first exit times.

\noindent {\it 2020 Mathematics Subject Classification:} 91B30, 60G51, 65M15, 65N15.
\end{abstract}

\section{Introduction}

Stochastic processes with jumps have long attracted a great deal of attention in a wide variety of applied areas since the distant past to the present day, ranging from social sciences \cite{10.1007/s00211-003-0511-8, https://doi.org/10.1111/1468-0262.00164} and natural sciences  \cite{https://doi.org/10.1016/0303-2647(96)01632-2, PhysRev.120.1093} to
information and computer sciences \cite{10.1109/34.809109}, to name just a few.
By introducing jumps, the underlying dynamics cannot only be made richer with possibly heavy-tailed marginals but also literally model characteristic jumps at a great sacrifice of analytical tractability.
For instance, it is essential in actuarial science that the underlying dynamics contains jumps to represent insurance claims \cite{GERBER1997129, EDPF}.
Most of the existing literature has focused on exploring analytical or, at least, semi-analytical solutions,
whereas, due to the inherent complexity of the dynamics with jumps, the analytical solution can only be found in a limited number of rather classical models. 
In practice, insurance companies are more concerned with the risks in a finite time horizon, where the availability of analytical solutions is even more significantly limited. 

Problems involving jump processes have thus been tackled by exploring other avenues, in most instances, through approximation.
For instance, 
in inversion-based approaches, computationally intensive numerical methods are often necessary even when the characteristic functions of the transition distribution are available in closed form \cite{https://doi.org/10.1111/1468-0262.00164}, where the inversion-based density evaluation suffers from a prohibitive load for each parameter set in numerical procedure.
Others include numerical methods for sample path generation \cite{casellaroberts2011, doi:10.1287/opre.2013.1191, GR2013, 10.3150/14-BEJ676}, for solving the relevant partial integro-differential equations \cite{10.1007/s00211-003-0511-8, 10.2307/2157318} and for constructing hard bounding functions \cite{doi:10.1137/110841497, doi:10.1137/151004070}. 
As is often the case, each has weaknesses in terms of problem dimension and constraints, regularity requirements, implementability and complexity.  


The aim of the present work is to establish recursive representations for a large class of stochastic differential equations with jumps.
The underlying class is large enough to accommodate the multivariate time-inhomogeneous dynamics, even in the absence of uniform ellipticity.
Equally important is that its jump component can be governed by a general time-state dependent unbounded intensity \cite{GR2013, 10.3150/14-BEJ676}, particularly not of rather restrictive L\'evy-driven type that essentially benefits a lot from independent and stationary increments.
In the proposed framework, a sequence (Theorem \ref{theorem wmrecurrsive}), in fact, two sequences (Theorem \ref{theorem vm wm}) of approximate solutions are set up with a jump time of the underlying dynamics as an information relay point in passing the past information on to a previous iteration step to fill in the missing information on the unobserved trajectory ahead.
We prove that the resulting iteration scheme is convergent, possibly slow in the pre-asymptotic regime yet exponentially fast in the limit (Theorem \ref{theorem convergence}) and, moreover, can even be monotonically convergent (Theorem \ref{theorem vm wm}) under suitable assumptions imposed on the input data.
In particular, if approximation solutions are smooth enough, then, by using the approximate solution at each iteration step, upper and lower bounding functions can be constructed (Theorem \ref{theorem when smooth}) in a similar spirit to \cite{doi:10.1137/151004070}, which form, so to speak, a $100\%$ confidence interval for the true solution.

It turns out that those approximate solutions can be represented under the probability measure with its jump component suppressed (Theorem \ref{theorem back to E0}), that is, the underlying dynamics with jumps is then measured as though the dynamics contained no jumps.
The iterative nature of those approximate solutions as a sequence in the original form does not look like a Picard iteration, whereas they, after this ``de-jump'' transform, may well do.
This result is not only possibly effective from a numerical point of view (for the obvious reason of the absence of cumbersome time-state dependent jumps), but also insightful from a broader perspective as the approximate solution can then be put in contrast with the Feynman-Kac representation of a suitable initial boundary value problem for partial differential equations (with the corresponding cumbersome integral term dropped), provided that an additional set of sufficient conditions is imposed on the problem parameter (Theorem \ref{theorem PDE equivalence}), for which a variety of numerical methods, such as finite element methods, are available in the literature.

As is always the case in the context of time-state dependent jump rates, the proposed framework can further benefit a lot from the uniform boundedness of the jump rate \cite{casellaroberts2011}.
For instance, the proposed iteration scheme can then be made more implementable via the Poisson thinning (Theorem \ref{theorem thinning}).
In addition, the modes of convergence of the sequences of approximate solutions and associated hard bounding functions can be strengthened (Theorem \ref{theorem uniform}).
Remarkably, all those convergence results are valid irrespective of problem dimension, whereas most existing numerical methods have mainly been focused on low, or even one, dimensions.
Despite the base approach is, in a broad sense, centered around the one in the (compound) Poisson disorder problem (for instance, \cite{37528059}), the significance of our recursive representation is that it succeeds to fully decouple general time-state dependent jumps from the underlying jump-diffusion processes by breaking the original dynamics into somewhat more tractable pieces, as in \cite{QIU2022109301, qiu2022iterative}, and provide the associated hard bounding functions at every step and the rate of convergence. 



The rest of this paper is set out as follows.
In Section \ref{section preliminaries}, we summarize and formulate the background materials.
We begin Section \ref{section main results} with the construction of approximate solutions and their convergence towards the true solution as a sequence (Theorem \ref{theorem convergence}).
We then present theoretical findings in turn (Theorems \ref{theorem wmrecurrsive}, \ref{theorem vm wm} and \ref{theorem back to E0}, respectively, in Sections \ref{section first few jumps}, \ref{section two routes} and \ref{section E0}), without additional technical conditions, and construct hard bounding functions (Theorem \ref{theorem when smooth}) in Section \ref{section hard bounding functions}.
We close the main section by developing the Poisson thinning (Theorem \ref{theorem thinning}) and strengthening the modes of convergence (Theorem \ref{theorem uniform}), provided that the jump rate is uniformly bounded.
In Section \ref{section numerical illustrations}, we examine two illustrative examples so as to justify out theoretical findings.
To maintain the flow of the paper, we collect all derivations and proofs in the Appendix.

\section{Problem formulation}\label{section preliminaries}

We first develop the notation that will be used throughout the paper and introduce the underlying stochastic process.
We denote the set of real numbers by $\mathbb{R}$, the set of natural numbers excluding $0$ by $\mathbb{N} := \{1,2,\cdots\}$, and $\mathbb{N}_0:=\mathbb{N}\cup \{0\}$.
We denote by $\mathbb{R}^d$ the $d$-dimensional Euclidean space with $\mathbb{R}^d_0:=\mathbb{R}^d\setminus \{0\}$.
We fix a non-empty open convex subset $D$ of $\mathbb{R}^d$ with Lipschitz boundary $\partial D$, and denote by $\overline{D}:=D\cup\partial D$ its closure.
We denote by $\partial_1$ the partial derivative with respect to the first argument and write $A^{\otimes 2}:=AA^{\top}$ for the matrix $A$.

Let $(\Omega,\mathcal{F},\mathbb{P})$ be a probability space and let $(\mathcal{F}_t)_{t\ge 0}$ be an increasing family of sub-$\sigma$-fields of $\mathcal{F}$.
All stochastic processes are assumed to be adapted to $(\mathcal{F}_t)_{t\ge 0}$.  
Consider the stochastic process $\{X_t:\,t\geq 0\}$ in $\mathbb{R}^d$ defined by the stochastic differential equation:
\begin{equation}\label{Xprocess}
dX_t = b(t,X_t)dt + \sigma(t,X_t)dW_t + \int_{\mathbb{R}_0^d} {\bf z}\mu(d{\bf z}, dt;X_{t-}), \quad X_0\in \overline{D},
\end{equation}
where $\{W_t: t\geq 0\}$ is the standard Brownian motion in $\mathbb{R}^{d_0}$, and the coefficients $b$ and $\sigma$ are locally bounded in time and Lipschitz continuous with a linear growth bound in state on $\mathbb{R}^d$.
For every ${\bf x}\in \overline{D}$, $\mu(d{\bf z},dt;{\bf x})$ denotes a Poisson counting measure on $\mathbb{R}_0^d \times[0,\infty)$, whose compensator is $\lambda(t,{\bf x})\nu(d{\bf z};{\bf x}) dt$, where $\lambda(t,{\bf x})$ is a non-negative function mapping from $[0,\infty)\times \overline{D}$ to $[0,\infty)$, and $\nu(d{\bf z};{\bf x})$ is a finite Levy measure on $\mathbb{R}_0^d$.
As usual, we assume that the Brownian motion $\{W_t:\,t\ge 0\}$ and the Poisson counting measure $\mu(d{\bf z},dt;{\bf x})$ are mutually independent for ${\bf x}\in \overline{D}$.

Since the L\'evy measure $\nu(d\cdot;{\bf x})$ is assumed to be finite \cite{Platen2010Numerical}, it loses no generality to impose the normalization condition $\int_{\mathbb{R}_0^d} \nu(d{\bf z};{\bf x})=1$ for all ${\bf x}\in \overline{D}$.
Moreover, since the singularity of jump sizes can be well captured by the L\'evy measure $\nu(d{\bf z};{\bf x})$ (for instance, a finite convolution of Poisson processes with different jump sizes), it is not restrictive to assume that the rate function $\lambda(t,\cdot)$ is locally Lipschitz for all $t\in [0,T]$.
In particular, as the rate function $\lambda$ and the L\'evy measure $\nu$ depend the state, the jump component in \eqref{Xprocess} affects itself in a non-trivial way, that is, each jump can have an impact, not only on the current and future jump intensity and magnitude, but also on the drift and diffusion components, and, in fact, vice versa.
Note that the rate function $\lambda$ does not need to be uniformly bounded \cite{doi:10.1287/opre.2013.1191, GR2013, 10.3150/14-BEJ676} until we impose so in Section \ref{section uniformly bounded jump rate}.

We denote by $\mathbb{P}_\lambda^{t,{\bf x}}(\cdot) := \mathbb{P}_\lambda(\cdot\vert X_t={\bf x})$ the probability measure under which $X_t={\bf x}$ almost surely and the rate of jumps is governed by the function $\lambda$, where $\lim_{\delta\to 0+}\delta^{-1}\mathbb{P}_{\lambda}^{t,{\bf x}}(X_{t+\delta}\neq X_t) = \lambda(t,{\bf x})$ for almost all $(t,{\bf x})\in [0,T]\times \overline{D}$, due to the standardized assumption $\int_{\mathbb{R}_0^d} \nu(d{\bf z};{\bf x})=1$ for all ${\bf x}\in \overline{D}$.
As usual, we denote by $\mathbb{E}_\lambda^{t,{\bf x}}$ the expectation taken under the probability measure $\mathbb{P}^{t,{\bf x}}_\lambda$. 
By formally setting $\lambda\equiv 0$ in the probability measure $\mathbb{P}_{\lambda}^{t,{\bf x}}$ and the expectation $\mathbb{E}_{\lambda}^{t,{\bf x}}$, one may switch on and off the jump component in the underlying process \eqref{Xprocess}.

\begin{remark}{\rm \label{first remark}
(i) The underlying process \eqref{Xprocess} can be thought of as an approximation of that with an infinite L\'evy measure $\int_{\mathbb{R}_0^d}\nu(d{\bf z};{\bf x})=+\infty$.
For instance, the jump component can be parameterized as  $\int_{\mathbb{R}_0^d} {\bf z}\mu_n(d{\bf z}, dt;X_{t-})$, where $\mu_n(d{\bf z},dt;{\bf x})$ denotes a Poisson counting measure on $\mathbb{R}_0^d \times[0,\infty)$ for ${\bf x}\in \overline{D}$, whose compensator is now $\lambda(t,{\bf x})\nu_n(d{\bf z};{\bf x}) dt$, with $\int_{\mathbb{R}_0^d}\nu_n(d{\bf z};{\bf x})=n$ for $n\in\mathbb{N}$.
Truncation by a sequence of L\'evy measures of increasing intensities (that is, $\{\nu_n(d{\bf z};{\bf x})\}_{n\in\mathbb{N}}$) has already been well studied for major L\'evy measures.
See \cite{2021arXiv210110533Y} for a comprehensive survey on this topic. 

\noindent (ii) The jump component can be generalized to, for instance, $\int_{\mathbb{R}_0^d} \mathfrak{m}(t,X_{t-},{\bf z})\mu(d{\bf z}, dt;X_{t-})$ as in \cite{doi:10.1137/151004070}, so as to modulate each jump at $(t,X_{t-})$ by a mapping $\mathfrak{m}(t,X_{t-},\cdot)$.
In the present paper, we do not adopt this generalization, as it does not seem to add essential value in the presence of the time-state dependent rate  $\lambda$ and the state dependent Poisson counting measure $\mu(d{\bf z},dt;{\bf x})$.
\qed}\end{remark}

Throughout, we reserve $T$ for the terminal time of the interval of interest $[0,T]$.
Define the first exit time of the underlying process $\{X_t:\,t\ge 0\}$ from the closure of the domain $\overline{D}$ later than a fixed time $t$, as well as its capped one by the terminal time $T$, by
\begin{equation}\label{etadef}
 \eta_t:=\inf\left\{s>t:X_s \not\in \overline{D}\right\},\quad \eta_t^T:=\eta_t\land T.
\end{equation}
The aim of the present work is to construct a theoretical framework for iterative weak approximation schemes for the function given in the form of expectation:
\begin{equation}\label{udef}
u(t,{\bf x}) := \mathbb{E}_{\lambda}^{t,{\bf x}} \left[\mathbbm{1} (\eta_t > T) \Theta_{t,T} g(X_T) + \mathbbm{1} (\eta_t \leq T) \Theta_{t,\eta_t} \Psi(\eta_t, X_{\eta_t},X_{\eta_t-}) -\int_t^{\eta_t^T} \Theta_{t,s} \phi(s, X_s) ds\right], \quad (t,{\bf x})\in[0,T]\times\overline{D},
\end{equation}
with $g:\overline{D}\to \mathbb{R}$, $\Psi:[0,T]\times(\mathbb{R}^d\setminus D)\times \overline{D} \to \mathbb{R}$, $\phi:[0,T]\times \overline{D}\to \mathbb{R}$ and 
\begin{equation}\label{defTheta}
 \Theta_{t_1,t_2} := \exp \left[ -\int_{t_1}^{t_2} r(s,X_s) ds\right],\quad 0\le t_1\le t_2\le T.
\end{equation}
For instance, the (finite-time version of) expected discounting penalty function \cite{GERBER1997129, EDPF} 
is a special case, given by
\[
 u(t,x)=\mathbb{E}_\lambda^{t,x}\left[e^{-r(\eta_t^T-t)}w(X_{\eta_t^T -}, |X_{\eta_t^T}|)\right],
\quad (t,x)\in [0,T]\times (0,+\infty),
\]
with a constant discount rate $r(s,{\bf x})\equiv r\ge 0$ (that is, $\Theta_{t_1,t_2}=e^{-r(t_2-t_1)}$) and a spectrally negative L\'evy measure, where the function $w$ here is called the penalty function. 
Here, on the event $\{\eta_t>T\}$, we have $X_{\eta_t^T}=X_{\eta_t^T-}=X_T$, that is, $w(X_{\eta_t^T -},|X_{\eta_t^T}|)=w(X_T,X_T)$, corresponding to the terminal condition $g$ in \eqref{udef}.
On the event $\{\eta_t\le T\}$, the term $w(X_{\eta_t^T-},|X_{\eta_t^T}|)$ can be represented by the boundary condition $\Psi$ in \eqref{udef} with its time argument suppressed.

Before moving to the main section, we summarize the basic assumptions that stand throughout.

\begin{assumption}{\rm \label{standing assumption}
{\bf (a)} The domain $D$ is a non-empty convex subset of $\mathbb{R}^d$ with Lipschitz boundary $\partial D$.

\noindent {\bf (b)} For every ${\bf x}\in \overline{D}$, the intensity measure $\nu(d{\bf z};{\bf x})$ is a L\'evy and probability measure on $\mathbb{R}_0^d$.

\noindent {\bf (c)} The input data $g$, $\Psi$ and $\phi$ are bounded and Lipschitz continuous on their respective domains and $g$ has polynomially bounded derivatives.

\noindent {\bf (d)} The coefficients $b$ and $\sigma$ are locally bounded in time, and are Lipschitz continuous with a linearly growth bound in state.\qed

}\end{assumption}

\section{Main results}\label{section main results}

We are ready to start our main section.
First, define the time of the first jump occurring after time $t$ by
\begin{equation}\label{tau1}
\tau_t^{(1)} := \inf\left\{s>t: X_{s-} \neq X_{s}\right\},
\end{equation}
and then, recursively, the time of the $m$-th jump occurring after time $t$ by
\begin{equation}\label{taum}
\tau_t^{(m)} := \inf\left\{s>\tau_t^{(m-1)}: X_{s-} \neq X_{s}\right\},
\end{equation}
for $m\in \mathbb{N}$, with the time of the zero-th jump after time $t$ defined formally by $\tau_t^{(0)}:=t$.
Next, we define a sequence $\{w_m\}_{m\in \mathbb{N}_0}$ of functions $[0,T]\times \overline{D}\to \mathbb{R}$ that play a central role as iterative approximate solutions for the solution $u$ defined in \eqref{udef}.
To begin the iteration, define the function $w_0:[0,T]\times \overline{D}\to \mathbb{R}$ by  
\begin{equation}\label{w0def1}
w_0(t,{\bf x}):=
\mathbb{E}_{0}^{t,{\bf x}} \left[\mathbbm{1} (\eta_t > T) \Theta_{t,T} g(X_T) + \mathbbm{1} (\eta_t \leq T) \Theta_{t,\eta_t} \Psi(\eta_t, X_{\eta_t},X_{\eta_t-}) -\int_t^{\eta_t^T} \Theta_{t,s} \phi(s, X_s) ds\right], \qquad (t,{\bf x})\in[0,T]\times \overline{D},
\end{equation}
which we let play a role as the initial approximate solution for the proposed iteration scheme.
It is worth noting that the expectation \eqref{w0def1} is taken under the probability measure $\mathbb{P}_0^{t,{\bf x}}$, namely, the jump component is suppressed with zero rate $\lambda\equiv 0$.
In this case, thus, we observe no jumps, to say nothing of overshooting, that is, $\Psi(\eta^T_t, X_{\eta^T_t},X_{\eta^T_t-})=\Psi(\eta_t^T, X_{\eta_t^T},X_{\eta_t^T})$, $\mathbb{P}_0^{t,{\bf x}}$-$a.s.$
Then, on the basis of the initial approximate solution \eqref{w0def1}, define the function $w_m$ by
\begin{multline}\label{wmdef2}
w_m(t,{\bf x}):=
\mathbb{E}_{\lambda}^{t,{\bf x}} \Bigg[\mathbbm{1} (X_{\eta_t^T\land \tau_t^{(m)}}\in \overline{D}) \Theta_{t,\eta_t^T\land \tau_t^{(m)}} w_0(\eta_t^T\land \tau_t^{(m)}, X_{\eta_t^T\land \tau_t^{(m)}})\\
+ \mathbbm{1} (X_{\eta_t^T\land \tau_t^{(m)}}\not\in \overline{D} ) \Theta_{t,\eta^T_t} \Psi(\eta^T_t, X_{\eta^T_t},X_{\eta^T_t-})
-\int_t^{\eta_t^T\land \tau_t^{(m)}} \Theta_{t,s} \phi(s, X_s) ds\Bigg],
\end{multline}
for $m\in \mathbb{N}$ and $(t,{\bf x})\in[0,T]\times \overline{D}$.
We first claim that the sequence $\{w_m\}_{m\in \mathbb{N}_0}$ is convergent to $u$, that is, the function of interest \eqref{udef}. 

\begin{theorem}\label{theorem convergence}
It holds, as $m\to +\infty$, that $w_m(t,{\bf x})\to u(t,{\bf x})$ for all $(t,{\bf x})\in [0,T]\times \overline{D}$.
\end{theorem}

\subsection{Recursive representations by paths until the first few jumps}
\label{section first few jumps}

On the basis of Theorem \ref{theorem convergence}, the function $w_m$, with a sufficiently large $m$,  can be employed as an approximation for $u$, while a large $m$ is nothing but observing the sample path until the $m$-th jump time $\tau_t^{(m)}$, according to the representation \eqref{wmdef2}.
In the presence of jumps, there would be no essential difference from the original form \eqref{udef} anymore.
In particular, when the underlying stochastic process \eqref{Xprocess} is employed as a finite-intensity approximation of an infinite-intensity counterpart (Remark \ref{first remark} (i)), an extremely large number of jumps would be required for reasonable accuracy of weak approximation.

This issue can be addressed in part with the aid of the following result, which leads to the construction of an iterative weak approximation scheme. 

\begin{theorem}\label{theorem wmrecurrsive}
Let $m\in \mathbb{N}$ and $n\in\{0,1,2,\cdots,m\}$.
It holds that
\begin{multline}\label{wmrecur2}
w_m(t,{\bf x}) =\mathbb{E}_\lambda^{t,{\bf x}} \Bigg[ \mathbbm{1}(X_{\eta_t^T\land \tau_t^{(n)}}\in\overline{D})\Theta_{t, \eta_t^T \land \tau_t^{(n)} } w_{m-n} (\eta_t^T\land \tau_t^{(n)} , X_{\eta_t^T\land \tau_t^{(n)} }) \\
+\mathbbm{1}(X_{\eta_t^T\land \tau_t^{(n)}}\not\in\overline{D})\Theta_{t,\eta_t} \Psi(\eta_t, X_{\eta_t},X_{\eta_t-}) 
- \int_t^{\eta_t^T \land \tau_t^{(n)}} \Theta_{t,s} \phi(s,X_s) ds \Bigg],
\end{multline}
for all $(t,{\bf x})\in[0,T]\times \overline{D}$.
\end{theorem}

In comparison with the representation \eqref{wmdef2} of the approximate solutions $\{w_m\}_{m\in\mathbb{N}_0}$, a few evident yet interesting features of the equivalent representation \eqref{wmrecur2} can be summarized as follows.
First of all, the representation \eqref{wmrecur2} can obviously recover the original one \eqref{wmdef2} by setting $n=m$.
The opposite extreme is the case $n=0$, which yields the trivial identity $w_m=w_m$, due to $\tau_t^{(0)}=t$ by definition.
When $n\in \{1,\cdots,m-1\}$ in between, the representation \eqref{wmrecur2} requires one to make an observation of the sample path until, at furthest, its $n$-th jump time $\tau_t^{(n)}$ (no longer all the way up to the $m$-th jump as for the previous representation \eqref{wmdef2}), on the basis of the previous approximate solution $w_{m-n}$ at $n$ steps ago.
Here, the $n$-th jump time $\tau_t^{(n)}$ can be thought of as an information relay point of passing the sample path behind to the previous step $w_{m-n}$, which fills in the missing information on the unobserved trajectory ahead. 

Before moving on to further theoretical developments, we provide an illustrative insight into the approximate solution $w_m$ with the aid of the representation \eqref{wmrecur2} in the context of the finite-time ruin probability, which is of significant interest in insurance mathematics.
We remark that the additional conditions on the drift and diffusion coefficients are imposed so as to ensure that an exit can only occur by a jump.

\begin{proposition}\label{proposition pure jump ruin probability}
Suppose that, for every $t\in [0,T]$, the drift $b(t,\cdot)$ is not outward pointing at the boundary $\partial D$, assume that there exists $\epsilon>0$ such that $\int_{[0,T]\times \partial D_{\epsilon}}\sigma(t,{\bf x})dtd{\bf x}=0$ where $\partial D_{\epsilon}:=\{{\bf x}\in \overline{D}:\,{\rm dist}({\bf x},\partial D)\le \epsilon\}$, and let $r\equiv 0$, $g\equiv 0$, $\Psi\equiv 1$ and $\phi\equiv 0$ in \eqref{udef}.
It then holds that for $m\in\mathbb{N}$ and $(t,{\bf x})\in [0,T]\times \overline{D}$,
\begin{multline*}
 w_m(t,{\bf x})=\mathbb{E}_\lambda^{t,{\bf x}}\bigg[
\mathbbm{1}(X_{T\land \tau_t^{(1)}}\not\in\overline{D})
+\mathbbm{1}(\{X_{T\land \tau_t^{(1)}}\in\overline{D}\}\cap \{X_{T\land \tau_t^{(2)}}\notin\overline{D}\})
\\
+\cdots +\mathbbm{1}(\{X_{T\land \tau_t^{(1)}}\in\overline{D}\}\cap
\cdots\cap\{X_{T\land \tau_t^{(m-1)}}\in\overline{D}\}\cap
\{X_{T\land \tau_t^{(m)}}\not\in\overline{D}\})
\bigg],
\end{multline*}
and $w_m(t,{\bf x})\ge w_{m-1}(t,{\bf x})$, with $w_0\equiv 0$.
\end{proposition}

In words, in the context of the finite-time ruin probability, the approximate solution $w_m$ represents the probability of ruin earlier than or right by the $m$-th jump.
Hence, the convergence towards the true ruin probability $u$ as well as its monotonicity are quite convincing. 
The monotonicity here is revisited in a general setting later in Theorem \ref{theorem vm wm}.

\subsection{Iterative approximation along two routes with monotonicity}
\label{section two routes}

Apart from the iteration scheme based on the representation \eqref{wmrecur2}, we may construct, so to speak, a side road, which allows the iteration to proceed along an alternative route, involving the sequence $\{v_m\}_{m\in \mathbb{N}}$, defined by
\begin{multline}\label{def of vm}
v_m(t,{\bf x}):=\mathbb{E}_\lambda^{t,{\bf x}}\Bigg[
\mathbbm{1}(\tau_t^{(m)} \geq \eta_t^T)
\mathbbm{1}(X_{\eta_t^T}\in\overline{D})\Theta_{t,\eta_t^T}w_0(\eta_t^T,X_{\eta_t^T})\\
+\mathbbm{1}(X_{\eta_t^T\land\tau_t^{(m)}}\not\in\overline{D})\Theta_{t,\eta_t^T}\Psi(\eta_t^T,X_{\eta_t^T},X_{\eta_t^T-})-\int_t^{\eta_t^T\land\tau_t^{(m)}}\Theta_{t,s}\phi(s,X_s)ds\Bigg],
\end{multline}
for $(t,{\bf x})\in [0,T]\times \overline{D}$, which initiates itself with $v_1$, unlike the sequence $\{w_m\}_{m\in\mathbb{N}_0}$ begins with $w_0$.
Other than that, the sequence $\{v_m\}_{m\in\mathbb{N}}$ here differs only slightly from the original one \eqref{wmdef2} by the event $\{\tau_t^{(m)}<\eta_t^T\}\cap \{X_{\tau_t^{(m)}}\in \overline{D}\}$ in the first term of \eqref{wmdef2},
on which the sample path makes the $m$-th jump before the terminal time $T$ or the first exit time $\eta_t$, and does not exit due to the $m$-th jump.
The sequence $\{v_m\}_{m\in\mathbb{N}}$ is convergent to the true solution $u$, along with $\{w_m\}_{m\in\mathbb{N}_0}$, as follows.

\begin{theorem}\label{theorem vm wm}
{\bf (a)} It holds that for $(t,{\bf x})\in [0,T]\times \overline{D}$, $m\in \mathbb{N}$ and $n\in \{0,1,\cdots,m-1\}$,
\begin{align}
 v_{m+1}(t,{\bf x})&=v_{m-n}(t,{\bf x})+\mathbb{E}_\lambda^{t,{\bf x}}\left[\mathbbm{1}(\tau_t^{(m-n)}< \eta_t^T)\Theta_{t,\tau_t^{(m-n)}}v_{n+1}(\tau_t^{(m-n)},X_{\tau_t^{(m-n)}})\right],\label{iteration vm}\\
  w_m(t,{\bf x})&=v_{m-n}(t,{\bf x})+\mathbb{E}_\lambda^{t,{\bf x}}\left[\mathbbm{1}(\tau_t^{(m-n)}< \eta_t^T)\Theta_{t,\tau_t^{(m-n)}}w_n(\tau_t^{(m-n)},X_{\tau_t^{(m-n)}})\right],\label{wm vm}
\end{align}
and that $v_m(t,{\bf x})\to u(t,{\bf x})$ as $m\to +\infty$ for all $(t,{\bf x})\in [0,T]\times \overline{D}$.

\noindent {\bf (b)}
If $g\ge 0$, $\Psi\ge 0$ and $\phi\le 0$ (or, if $g\le 0$, $\Psi\le 0$ and $\phi\ge 0$), then the convergence in {\bf (a)} is monotonic.  
\end{theorem}

Hence, as long as the iteration begins with $w_0$, one may stick to a single route (either $\{w_m\}_{m\in\mathbb{N}_0}$ by Theorem \ref{theorem wmrecurrsive}, or $\{v_m\}_{m\in\mathbb{N}}$ based on \eqref{iteration vm}), or can even go back and forth between these two routes thanks to the interaction \eqref{wm vm}, and even further, may skip a few steps ($n$ in \eqref{iteration vm} and \eqref{wm vm}) in exchange for making longer observation to get more jumps accordingly. 
As it turns out shortly in Theorem \ref{theorem uniform}, one might make convergence efficient by skipping a few jumps because the convergence can be slow in the pre-asymptotic regime yet exponentially fast in the limit.

As has already been seen in Proposition \ref{proposition pure jump ruin probability}, the input data $r\equiv 0$, $g\equiv 0$, $\Psi\equiv 1$ and $\phi\equiv 0$ correspond to the finite-time ruin probability and yield a monotonically increasing sequence of approximate solutions $w_m\ge w_{m-1}$ with $w_0\equiv 0$.
Indeed, Theorem \ref{theorem vm wm} {\bf (b)} is fully consistent with those, since the input data for  the finite-time ruin probability satisfy the condition in {\bf (b)}, and moreover, since $w_0\equiv 0$ yields $w_m=v_m$ due to \eqref{wmdef2} and \eqref{def of vm} in comparison.

\subsection{Recursive representations by decoupling jumps}
\label{section E0}

The recursive relation \eqref{wmrecur2} allows one to represent the function $w_m$ on the basis of a previous one $w_{m-n}$ (at $n$ steps ago) and sample paths only up to their $n$-th jump. 
In this section, we further reduce the representation \eqref{wmrecur2} to an alternative one on the probability measure $\mathbb{P}_0^{t,{\bf x}}$, that is, under which the jump component is suppressed.
To this end, we prepare the following notation:
\begin{equation*}
 G_{m-1}(t,{\bf x}) := \int_{\{{\bf z}\in\mathbb{R}_0^d:{\bf x}+{\bf z}\in\overline{D}\}} w_{m-1}(t, {\bf x}+{\bf z}) \lambda (t,{\bf x}) \nu(d{\bf z};{\bf x}),
\quad H(t,{\bf x}):=\int_{\{{\bf z}\in\mathbb{R}_0^d:{\bf x}+{\bf z}\not\in\overline{D}\}} \Psi(t,{\bf x}+{\bf z}, {\bf x}) \lambda(t,{\bf x}) \nu(d{\bf z};{\bf x}),
\end{equation*}
for $(t,{\bf x})\in [0,T]\times\overline{D}$ and $m\in\mathbb{N}$, as well as
\[
\Lambda_{t_1,t_2}:=\exp\left[-\int_{t_1}^{t_2}\lambda(s,X_s)ds\right],\quad 0\leq t_1 \leq t_2\leq T.
\]
We note that the jump component still enters via the two quantities $G_{m-1}$ and $H$ above, yet only implicitly, because those are based on the previous step $w_{m-1}$ alone.

\begin{theorem}\label{theorem back to E0}
It holds that 
\begin{multline}\label{wme0}
w_m(t,{\bf x}) = \mathbb{E}_0^{t,{\bf x}} \Bigg[\mathbbm{1}(\eta_t>T) \Theta_{t,T} \Lambda_{t,T}g(X_T) +\mathbbm{1}(\eta_t\leq T) \Theta_{t,\eta_t} \Lambda_{t,\eta_t}\Psi(\eta_t, X_{\eta_t},X_{\eta_t})\\
-\int_t^{\eta_t\land T} \Theta_{t,s}\Lambda_{t,s} \left(\phi(s,X_s) - G_{m-1}(s,X_s) - H(s,X_s)\right) ds\Bigg],
\end{multline}
for $m\in\mathbb{N}$ and $(t,{\bf x})\in[0,T]\times\overline{D}$.
\end{theorem}

Let us stress again that the approximate solution \eqref{wme0} is written on the probability measure $\mathbb{P}_0^{t,{\bf x}}$, that is, with the jump component suppressed.
In a sense, the representation \eqref{wme0} can be thought of as a further step back from the shortest observation in \eqref{wmrecur2} (shortest, as with $n=1$), whereas not identical to its trivial identity $w_m=w_m$ with $n=0$.
In short, the sample path strictly before the first jump is identical in law under the two probability measures $\mathbb{P}_{\lambda}^{t,{\bf x}}$ and $\mathbb{P}_0^{t,{\bf x}}$, while the missing information caused by ceasing observation strictly before the first jump can still be compensated by the two terms $G_{m-1}$ and $H$ inside the integral, where the first term $G_{m-1}$ carries over the information from the previous step $w_{m-1}$.  
We note that the two terms $G_{m-1}(t,{\bf x})$ and $H(t,{\bf x})$ are deterministic functions and thus do not vanish even under $\mathbb{P}_0^{t,{\bf x}}$.

We note that the recursive representation \eqref{wme0} is a distinctive result in the sense that it succeeds to fully decouple cumbersome time-state dependent jumps (as with respect to the expectation $\mathbb{E}_0^{t,{\bf x}}$) from jump-diffusion processes, whereas the existing recursive representations, for instance, in the context of the (compound) Poisson disorder problem (such as \cite{37528059}) are written under the original probability measure (that is, $\mathbb{P}_{\lambda}^{t,{\bf x}}$ in our context), corresponding to our previous step (Theorem \ref{theorem wmrecurrsive}), that is, there is no such things as decoupling of jumps.

It is further worth stressing that no ellipticity condition has been imposed on the underlying stochastic process \eqref{Xprocess} whatsoever.
As an extreme example, if the diffusion component is fully suppressed ($\sigma\equiv 0$), then the underlying stochastic process \eqref{Xprocess} is fully degenerate under the probability measure $\mathbb{P}_0^{t,{\bf x}}$, that is, the sample path is a deterministic function under $\mathbb{P}_0^{t,{\bf x}}$.
The representation \eqref{wme0} can then reduce to 
\begin{multline}\label{fully deterministic case}
    w_m(t,{\bf x}) =\mathbbm{1}(\eta_t>T) \Theta_{t,T} \Lambda_{t,T}g(x(T)) +\mathbbm{1}(\eta_t\leq T) \Theta_{t,\eta_t} \Lambda_{t,\eta_t}\Psi(\eta_t, x(\eta_t),x(\eta_t))\\
    -\int_t^{\eta_t\land T} \Theta_{t,s}\Lambda_{t,s} \left(\phi(s,x(s)) - G_{m-1}(s,x(s)) - H(s,x(s))\right)ds,\quad (t,{\bf x})\in [0,T]\times \overline{D},
\end{multline}
where the function $x$ here denotes the unique solution to the deterministic ordinary differential equation $x'(s)=b(s,x(s))$ for $s\in [t,T]$ with initial condition $x(t)={\bf x}$.
Here, we have $\eta_t=\inf\{s\geq t:x(s)\not\in\overline{D}\}$, $\Theta_{t_1,t_2}=\exp[ - \int_{t_1}^{t_2} r(s,x(s))ds]$, and $\Lambda_{t_1,t_2}=\exp[ - \int_{t_1}^{t_2} \lambda(s,x(s)) ds]$, where the first exit time $\eta_t$ is also deterministic in this context.
In brief, if the diffusion component is suppressed, then the iteration can be conducted based on the recursive relation \eqref{fully deterministic case} in a fully deterministic manner.
From a broader perspective, the deterministic identity \eqref{fully deterministic case} reveals the fact that the proposed iterative scheme based on the representation \eqref{wme0} can thus be thought of as a (probabilistic) Picard iteration process.
In Section \ref{simple poisson example}, we provide a concrete example for comprehensive demonstration of this fully deterministic iteration.

\subsection{Hard bounding functions}\label{section hard bounding functions}

It would make the theoretical framework (Theorems \ref{theorem wmrecurrsive}, \ref{theorem vm wm} and \ref{theorem back to E0}) more comprehensive and definitive to provide hard bounds for the true solution in a similar spirit to \cite{doi:10.1137/151004070}.
To this end, we prepare some notation.
For $(t,{\bf x})\in [0,T]\times\mathbb{R}^d$, define 
\begin{equation}\label{def of xi and Mr}
\xi(t,{\bf x};\lambda):=\mathbb{E}_{\lambda}^{t,{\bf x}}\left[\int_t^{\eta_t^T} \Theta_{t,s} ds\right],\quad M(t,{\bf x}):=\int_{\{{\bf z}\in\mathbb{R}_0^d:{\bf x}+{\bf z}\not\in\overline{D}\}}
\xi(t,{\bf x}+{\bf z};0)\lambda(t,{\bf x})\nu(d{\bf z};{\bf x}) - \lambda(t,{\bf x})\xi(t,{\bf x};0),
\end{equation}
with $M^U(t):=\esssup_{(s,{\bf x})\in[t,T]\times \overline{D}}(M(s,{\bf x}))_+$ and $M^L(t):=\essinf_{(s,{\bf x})\in[t,T]\times \overline{D}}(M(s,{\bf x}))_-$,
and moreover,
\begin{equation}\label{def of Nm}
N_m(t,{\bf x}):=
\begin{dcases}
G_0(t,{\bf x})-\lambda(t,{\bf x})w_0(t,{\bf x})+H(t,{\bf x}),&\text{if }m=0,\\
G_{m}(t,{\bf x})-G_{m-1}(t,{\bf x}),&\text{if }m\in \mathbb{N},
\end{dcases}
\end{equation}
with the essential supremum $\{N_m^U(t)\}_{m\in\mathbb{N}_0}$ and infimum $\{N_m^L(t)\}_{m\in\mathbb{N}_0}$, defined as follows:
\begin{equation}\label{def of NmU and NmL}
N_m^U(t):=\esssup_{(s,{\bf x})\in[t,T]\times \overline{D}}\left(N_m(s,{\bf x})\right)_+, \quad 
N_m^L(t):=\essinf_{(s,{\bf x})\in[t,T]\times \overline{D}}\left(N_m(s,{\bf x})\right)_-,\quad t\in [0,T].
\end{equation}
Clearly, $N_m^U(\cdot)$ and $N_m^L(\cdot)$ are monotone for all $m\in \mathbb{N}_0$.
Also, with a sufficiently smooth $f$, define the differential operator $\mathcal{L}_{\cdot}$ by
\[
\mathcal{L}_t f(t,{\bf x}):=\langle b(t,{\bf x}),\nabla f(t,{\bf x})\rangle +\frac{1}{2}{\rm tr}\left[(\sigma (t,{\bf x}))^{\otimes 2} {\rm Hess}_{\bf x} f(t,{\bf x})\right],
\]
for $(t,{\bf x})\in [0,\infty)\times \mathbb{R}^d$ wherever it is finite valued.
We are now ready to state the result.

\begin{theorem}\label{theorem when smooth}
Let $m\in \mathbb{N}$.
If $w_m$ is smooth enough so that $\partial_1 w_m$ and $\mathcal{L}_{\cdot} w_m$ exist almost everywhere on $[0,T)\times D$, then the functions $w_m$ and $w_{m-1}$ together satisfy 
\begin{equation}\label{wmpde}
\partial_1 w_m(t,{\bf x}) + \mathcal{L}_t w_m(t,{\bf x}) = (r(t,{\bf x}) +\lambda(t,{\bf x})) w_m(t,{\bf x}) + \phi(t,{\bf x}) -G_{m-1}(t,{\bf x}) -H(t,{\bf x}),\quad a.e.\mbox{-}(t,{\bf x})\in [0,T)\times D,
\end{equation}
with $w_m(T,\cdot) = g(\cdot)$.
If, moreover, $M^U(0)<1$ and $\xi(\cdot,\cdot;0)$ is smooth enough so that $\partial_1$ and   $\mathcal{L}_{\cdot}$ exist almost everywhere on $[0,T)\times D$, then it holds that for $(t,{\bf x})\in [0,T]\times \overline{D}$,
\begin{equation}\label{boundswm0}
\frac{N_m^L(t)}{1-M^L(t)}\xi(t,{\bf x};0)\le u(t,{\bf x}) -w_m(t,{\bf x})\le \frac{N_m^U(t)}{1-M^U(t)}\xi(t,{\bf x};0).
\end{equation}
\end{theorem}

The inequalities \eqref{boundswm0} can be employed as hard bounding functions for the true solution $u(t,{\bf x})$ at an arbitrary step $m$ on the basis of the two consecutive approximate solutions $w_m$ and $w_{m-1}$ required for obtaining the function $N_m$ via the formula \eqref{def of Nm}, especially when all involving elements are accessible for numerical purposes.
From a theoretical point of view, the inequalities \eqref{boundswm0} play an essential role for convergence results (Theorem \ref{theorem uniform}).
We remark that the smoothness conditions on the functions $w_{m-1}$, $w_m$ and $\xi(\cdot,\cdot;0)$ above do not require the uniform ellipticity or regularity of the input data.
Indeed, shortly in Section \ref{simple poisson example}, we derive hard bounding functions \eqref{boundswm0} in a problem setting without diffusion component ($\sigma\equiv 0$), yet within the scope of Theorem \ref{theorem when smooth}.

Let us stress that the partial differential equation \eqref{wmpde} here is not a consequence of the Feynman-Kac formula, but holds true, simply based on the definition \eqref{wmdef2} without imposing a variety of additional sufficient conditions.
It may still not be of use on its own in the absence of concrete information about the boundary $\partial D$, while we still present it here because it plays a crucial role not only in deriving the hard bounding functions \eqref{boundswm0}, but also in proving the upcoming Theorem \ref{theorem uniform}, as well as reappears in Theorem \ref{theorem PDE equivalence} with more conditions imposed.

\subsection{Uniformly bounded jump rate}\label{section uniformly bounded jump rate}

We have now been equipped with a comprehensive recursive mechanism (Sections \ref{section first few jumps}, \ref{section two routes} and \ref{section E0}), along with hard bounding functions (Section \ref{section hard bounding functions}).
In this section, we show that the mechanism may be eased in exchange for imposing additional restrictions on the problem setting.
Namely, hereafter, the jump rate $\lambda$ is assumed to be uniformly bounded over the domain, that is, there exists a constant $\widetilde{\lambda}$ such that $\lambda(t,{\bf x})\le \widetilde{\lambda}$ for all $(t,{\bf x})\in [0,T]\times \overline{D}$.


The first benefit from the jump rate function being uniformly bounded is the Poisson thinning, as is standard, for instance, in the context of exact simulation of jump-diffusion processes \cite{casellaroberts2011}.
In our context, by introducing the Poisson thinning, one may simplify the computation of each approximate solution, particularly when Monte Carlo methods are employed, at the sacrifice of the iteration progress.

For its formulation, we prepare some notation.
For every $(t,{\bf x})\in[0,T]\times\overline{D}$, define the jump intensity measure $\widetilde{\nu}(d{\bf z};t,{\bf x})$ on $\mathbb{R}^d$ (here, including the origin) by
\[
\widetilde{\nu}(d{\bf z};t,{\bf x}):=\frac{\lambda(t,{\bf x})}{\widetilde{\lambda}}\nu(d{\bf z};{\bf x}), \quad {\bf z}\in \mathbb{R}_0^d,
\]
with a mass at the origin $\widetilde{\nu}(\{0\};t,{\bf x}):=1-\lambda(t,{\bf x})/\widetilde{\lambda}$.
Note that the measure $\widetilde{\nu}(d{\bf z};t,{\bf x})$ remains finite and standardized, that is, $\int_{\mathbb{R}^d}\widetilde{\nu}(d{\bf z};t,{\bf x})=\int_{\mathbb{R}_0^d}\widetilde{\nu}(d{\bf z};t,{\bf x})+\widetilde{\nu}(\{0\};t,{\bf x})=1$, just as the base L\'evy measure $\nu(d{\bf z};{\bf x})$.
In view of the identity $\int_{\mathbb{R}^d}{\bf z} \widetilde{\lambda}\widetilde{\nu}(d{\bf z};t,{\bf x}) =\int_{\mathbb{R}_0^d}{\bf z} \lambda(t,{\bf x})\nu(d{\bf z};{\bf x})$, one may then formulate a thinning method.
That is, when the underlying process is at $(t,{\bf x})$, a jump is generated with the largest rate $\widetilde{\lambda}$ and then accepted with time-state dependent probability $\lambda(t,{\bf x})/\widetilde{\lambda}$, whereas its rejection with the remaining probability $1-\lambda(t,{\bf x})/\widetilde{\lambda}$ is formally represented as a zero-sized jump.

Next, for every $(t,{\bf x})\in [0,T]\times \overline{D}$, denote by $\widetilde{\mathbb{P}}_{\widetilde{\lambda}}^{t,{\bf x}}$ and $\widetilde{\mathbb{E}}_{\widetilde{\lambda}}^{t,{\bf x}}$, respectively, the probability measure and the associated expectation under which the underlying process $\{X_s:\,s\ge\ 0\}$ is governed, not by \eqref{Xprocess} but, by the following stochastic differential equation
\begin{equation}\label{underlying process with jumps of size zero}
dX_s = b(s,X_s)ds + \sigma(s,X_s)dW_s + \int_{\mathbb{R}^d} {\bf z}\widetilde{\mu}(d{\bf z} ,ds;s,X_{s-}),
\end{equation}
with initial state $X_t={\bf x}$, where $\widetilde{\mu}(d{\bf z},dt;t,{\bf x})$ denotes an extended Poisson counting measure on $\mathbb{R}^d\times [0,\infty)$, whose compensator is given by $\widetilde{\lambda}\widetilde{\nu}(d{\bf z};t,{\bf x})dt$.
The random counting measure here should be said to be extended in the sense that it permits (more precisely, pretends to permit) zero-sized jumps with time-state dependent probability $\widetilde{\nu}(\{0\};t,{\bf x})$ out of all jumps generated with constant rate $\widetilde{\lambda}$ wherever the underlying process is.
The redefined underlying processes \eqref{underlying process with jumps of size zero} under the probability measure $\widetilde{\mathbb{P}}_{\widetilde{\lambda}}^{t,{\bf x}}$ is identical in law to the original process \eqref{Xprocess} under $\mathbb{P}_{\lambda}^{t,{\bf x}}$.
As such, the first exit time $\eta_t$ and the discount factor $\Theta_{t_1,t_2}$ can remain unaltered, respectively, as \eqref{etadef} and \eqref{defTheta}, between the two probability measures $\mathbb{P}_{\lambda}^{t,{\bf x}}$ and $\widetilde{\mathbb{P}}_{\widetilde{\lambda}}^{t,{\bf x}}$.

The real relevance of preparing the notation $\widetilde{\mathbb{P}}_{\widetilde{\lambda}}^{t,{\bf x}}$ and $\widetilde{\mathbb{E}}_{\widetilde{\lambda}}^{t,{\bf x}}$ lies in the way of counting jumps.
The redefined underlying process \eqref{underlying process with jumps of size zero} has been designed to make more jumps (due to the higher jump rate $\widetilde{\lambda}$), some of which are however invisible because of zero size.
Namely, the definition of the jump times \eqref{tau1} and \eqref{taum} is not capable of detecting such zero-sized jumps.
Hence, out of necessity of counting in zero-sized jumps, we redefine a sequence of jump times by $\{\widetilde{\tau}_t^{(m)}\}_{m\in\mathbb{N}}$ on the redefined underlying process \eqref{underlying process with jumps of size zero}, where $\widetilde{\tau}_t^{(m)}$ denotes its $m$-th jump time among all jumps, some of which can be zero-sized.
Accordingly, we redefine the discount function as $\widetilde{\Lambda}_{t_1,t_2}:=\exp[-\int_{t_1}^{t_2}\widetilde{\lambda}ds]=e^{-\widetilde{\lambda} (t_2-t_1)}$ for $0\le t_1\le t_2$, which is evidently no larger than the original one $\Lambda_{t_1,t_2}$, due to the uniform dominance $\lambda(t,{\bf x})\le \widetilde{\lambda}.$

Now, in light of the representations \eqref{wmdef2} and \eqref{def of vm}, we define the sequences $\{\widetilde{w}_m\}_{m\in\mathbb{N}_0}$ and $\{\widetilde{v}_m\}_{m\in\mathbb{N}}$ of approximate solutions with the common initial approximate solution $\widetilde{w}_0:=w_0$, respectively, by
\begin{multline}\label{defwmthinning1}
\widetilde{w}_m(t,{\bf x}):=
\widetilde{\mathbb{E}}_{\widetilde{\lambda}}^{t,{\bf x}} \Bigg[\mathbbm{1} (X_{\eta_t^T\land \widetilde{\tau}_t^{(m)}}\in \overline{D}) \Theta_{t,\eta_t^T\land \widetilde{\tau}_t^{(m)}} w_0(\eta_t^T\land \widetilde{\tau}_t^{(m)}, X_{\eta_t^T\land \widetilde{\tau}_t^{(m)}})\\
+ \mathbbm{1} (X_{\eta_t^T\land \widetilde{\tau}_t^{(m)}}\not\in \overline{D} ) \Theta_{t,\eta^T_t} \Psi(\eta^T_t, X_{\eta^T_t},X_{\eta^T_t-})
-\int_t^{\eta_t^T\land \widetilde{\tau}_t^{(m)}} \Theta_{t,s} \phi(s, X_s) ds\Bigg],
\end{multline}
and
\begin{multline}\label{defvmthinning}
\widetilde{v}_m(t,{\bf x}):=\widetilde{\mathbb{E}}_{\widetilde{\lambda}}^{t,{\bf x}}\Bigg[
\mathbbm{1}(\widetilde{\tau}_t^{(m)} \geq \eta_t^T)
\mathbbm{1}(X_{\eta_t^T}\in\overline{D})\Theta_{t,\eta_t^T}w_0(\eta_t^T,X_{\eta_t^T})\\
+\mathbbm{1}(X_{\eta_t^T\land\widetilde{\tau}_t^{(m)}}\not\in\overline{D})\Theta_{t,\eta_t^T}\Psi(\eta_t^T,X_{\eta_t^T},X_{\eta_t^T-})-\int_t^{\eta_t^T\land\widetilde{\tau}_t^{(m)}}\Theta_{t,s}\phi(s,X_s)ds\Bigg],
\end{multline}
for $(t,{\bf x})\in [0,T]\times \overline{D}$.
We also prepare the notation
\begin{align*}
\widetilde{G}_{m-1}(t,{\bf x})
:=\int_{\{{\bf z}\in\mathbb{R}^d:{\bf x}+{\bf z}\in\overline{D}\}} \widetilde{w}_{m-1}(t, {\bf x}+{\bf z})\widetilde{\lambda}\widetilde{\nu}(d{\bf z};t,{\bf x})
=(\widetilde{\lambda}-\lambda(t,{\bf x}))\widetilde{w}_{m-1}(t,{\bf x}) 
+ \int_{\{{\bf z}\in\mathbb{R}_0^d:{\bf x}+{\bf z}\in\overline{D}\}} \widetilde{w}_{m-1}(t, {\bf x}+{\bf z})\lambda(t,{\bf x})\nu(d{\bf z};{\bf x}),
\end{align*}
while it turns out unnecessary to introduce the notation $\widetilde{H}$ in a similar manner, due to
\[
\widetilde{H}(t,{\bf x}):=\int_{\{{\bf z}\in\mathbb{R}^d:{\bf x}+{\bf z}\not\in\overline{D}\}}\Psi(t,{\bf x}+{\bf z},{\bf x})\widetilde{\lambda}\widetilde{\nu}(d{\bf z};t,{\bf x})
=\int_{\{{\bf z}\in\mathbb{R}_0^d:{\bf x}+{\bf z}\not\in\overline{D}\}}\Psi(t,{\bf x}+{\bf z},{\bf x})\lambda(t,{\bf x})\nu(d{\bf z};t,{\bf x})=H(t,{\bf x}),
\]
since the event $\{X_{t-}\in\overline{D}\}\cap\{X_{t-}+{\bf z} \not\in\overline{D}\}$ can occur only if the jump size ${\bf z}$ is non-zero.
We are now ready to give the main result, which is, so to speak, the Poisson thinning version of Theorems \ref{theorem convergence}, \ref{theorem wmrecurrsive}, \ref{theorem vm wm}, \ref{theorem back to E0} and \ref{theorem when smooth}.
Note that the choice of the probability measure in the representation \eqref{defwmthinning}, $\widetilde{\mathbb{E}}_0^{t,{\bf x}}$ or $\mathbb{E}_0^{t,{\bf x}}$, does not affect the result in any way, as the jump component is suppressed in there.

\begin{theorem}\label{theorem thinning}
Suppose there exists a constant $\widetilde{\lambda}$ such that $\lambda(t,{\bf x})\le \widetilde{\lambda}$ for all $(t,{\bf x})\in [0,T]\times \overline{D}$.

\noindent {\bf (a)}
It holds that for $(t,{\bf x})\in [0,T]\times \overline{D}$, $m\in \mathbb{N}$ and $n\in \{0,1,\cdots,m-1\}$,
\begin{align}
\widetilde{w}_m(t,{\bf x})&=\widetilde{\mathbb{E}}_{\widetilde{\lambda}}^{t,{\bf x}} \Bigg[ \mathbbm{1}(X_{\eta_t^T\land \widetilde{\tau}_t^{(n)}}\in\overline{D})\Theta_{t, \eta_t^T \land \widetilde{\tau}_t^{(n)} } \widetilde{w}_{m-n} (\eta_t^T\land \widetilde{\tau}_t^{(n)} , X_{\eta_t^T\land \widetilde{\tau}_t^{(n)} }) \nonumber\\
&\qquad \qquad \qquad \qquad \qquad \qquad+\mathbbm{1}(X_{\eta_t^T\land \widetilde{\tau}_t^{(n)}}\not\in\overline{D})\Theta_{t,\eta_t} \Psi(\eta_t, X_{\eta_t},X_{\eta_t-}) - \int_t^{\eta_t^T \land \widetilde{\tau}_t^{(n)}} \Theta_{t,s} \phi(s,X_s) ds \Bigg]\nonumber\\
&=\mathbb{E}_0^{t,{\bf x}} \Bigg[\mathbbm{1}(\eta_t>T) \Theta_{t,T} \widetilde{\Lambda}_{t,T}g(X_T) +\mathbbm{1}(\eta_t\leq T) \Theta_{t,\eta_t} \widetilde{\Lambda}_{t,\eta_t}\Psi(\eta_t, X_{\eta_t},X_{\eta_t})\nonumber\\
&\qquad \qquad \qquad \qquad \qquad \qquad \qquad \qquad -\int_t^{\eta_t\land T} \Theta_{t,s}\widetilde{\Lambda}_{t,s} \left(\phi(s,X_s) - \widetilde{G}_{m-1}(s,X_s) - H(s,X_s)\right) ds\Bigg],\label{defwmthinning}
\end{align}
and 
\begin{align}
 \widetilde{v}_{m+1}(t,{\bf x})&=\widetilde{v}_{m-n}(t,{\bf x})+\widetilde{\mathbb{E}}_{\widetilde{\lambda}}^{t,{\bf x}}\left[\mathbbm{1}(\widetilde{\tau}_t^{(m-n)}< \eta_t^T)\Theta_{t,\widetilde{\tau}_t^{(m-n)}}\widetilde{v}_{n+1}(\widetilde{\tau}_t^{(m-n)},X_{\widetilde{\tau}_t^{(m-n)}})\right],\label{iteration vm thinning}\\
  \widetilde{w}_m(t,{\bf x})&=\widetilde{v}_{m-n}(t,{\bf x})+\widetilde{\mathbb{E}}_{\widetilde{\lambda}}^{t,{\bf x}}\left[\mathbbm{1}(\widetilde{\tau}_t^{(m-n)}< \eta_t^T)\Theta_{t,\widetilde{\tau}_t^{(m-n)}}\widetilde{w}_n(\widetilde{\tau}_t^{(m-n)},X_{\widetilde{\tau}_t^{(m-n)}})\right].\label{wm vm thinning}
\end{align}

\noindent {\bf (b)} It holds, as $m\to +\infty$, that $\widetilde{w}_m(t,{\bf x})\to u(t,{\bf x})$ and $\widetilde{v}_m(t,{\bf x})\to u(t,{\bf x})$ for all $(t,{\bf x})\in[0,T]\times\overline{D}$.
In particular, if $g\ge 0$, $\Psi\ge 0$ and $\phi\le 0$ (or, if $g\le 0$, $\Psi\le 0$ and $\phi\ge 0$), then the latter convergence is monotonic.  

\noindent {\bf (c)}
Let $m\in \mathbb{N}$.
If $M^U(0)<1$, and if $\widetilde{w}_{m-1}$, $\widetilde{w}_m$ and $\xi(\cdot,\cdot;0)$ are smooth enough so that $\partial_1$ and $\mathcal{L}_{\cdot}$ exist almost everywhere on $[0,T)\times D$, then it holds that
\begin{equation}\label{boundswm02}
\frac{\widetilde{N}_m^L(t)}{1-M^L(t)}\xi(t,{\bf x};0)\le u(t,{\bf x}) -\widetilde{w}_m(t,{\bf x})\le \frac{\widetilde{N}_m^U(t)}{1-M^U(t)}\xi(t,{\bf x};0),\quad (t,{\bf x})\in [0,T]\times \overline{D},
\end{equation}
where $\widetilde{N}_m^L(t)$ and $\widetilde{N}_m^U(t)$ are, respectively, the essential infimum and supremum of $\widetilde{N}_m(s,{\bf x}):=\widetilde{G}_m(s,{\bf x})-\widetilde{G}_{m-1}(s,{\bf x})$ on $[t,T]\times\overline{D}$.
\end{theorem}


Another benefit from the jump rate function being uniformly bounded is the uniform convergence of the sequences $\{w_m\}_{m\in\mathbb{N}_0}$ and $\{v_m\}_{m\in\mathbb{N}}$ towards the true solution $u$, and as a consequence, the convergence of the upper and lower bounding functions \eqref{boundswm0} to each other.
The results indicate that the convergence is possibly slow in the pre-asymptotic regime if the intensities $\lambda_0$ and $\widetilde{\lambda}T$ are not small, yet exponentially fast in the limit.
We stress that the latter is a benefit from the former uniform convergence, because  the supremum and infimum \eqref{def of NmU and NmL} in the hard bounding functions are taken over the entire domain.
In addition, the results here provide an important insight into a trade-off whether or not the Poisson thinning should be applied, in terms of iteration progress per step relative to the computational complexity.

\begin{theorem}\label{theorem uniform}
Suppose there exists a constant $\widetilde{\lambda}$ such that $\lambda(t,{\bf x})\le \widetilde{\lambda}$ for all $(t,{\bf x})\in [0,T]\times \overline{D}$, and suppose that the function $w_0$ is smooth enough so that $\partial_1$ and $\mathcal{L}_{\cdot}$ exist almost everywhere on $[0,T)\times D$.
If $|N_0^U(0)|+|N_0^L(0)|<+\infty$, then it holds, as $m\to +\infty$, that 
\[
\sup_{(t,{\bf x})\in [0,T]\times \overline{D}}|w_m(t,{\bf x})-u(t,{\bf x})|=\mathcal{O}(\lambda_0^m/m!),\quad \sup_{(t,{\bf x})\in [0,T]\times \overline{D}}|v_m(t,{\bf x})-u(t,{\bf x})|=\mathcal{O}(\lambda_0^m/m!),
\]
with $\lambda_0:=\sup_{t\in [0,T]}\sup_{(s,{\bf x})\in [t,T]\times \overline{D}}\{\lambda(s,{\bf x})(T-t)\}$, and that
\[
\sup_{(t,{\bf x})\in [0,T]\times \overline{D}}|\widetilde{w}_m(t,{\bf x})-u(t,{\bf x})|=\mathcal{O}((\widetilde{\lambda} T)^m/m!),\quad \sup_{(t,{\bf x})\in [0,T]\times \overline{D}}|\widetilde{v}_m(t,{\bf x})-u(t,{\bf x})|=\mathcal{O}((\widetilde{\lambda} T)^m/m!).
\]
Moreover, all the bounding functions of \eqref{boundswm0} and \eqref{boundswm02} tend pointwise to zero as $m\to +\infty$.
\end{theorem}

Due to the dominance $\lambda_0 \le \widetilde{\lambda}T$, the difference in convergence rate above, between $\mathcal{O}(\lambda_0^m/m!)$ and $\mathcal{O}((\widetilde{\lambda} T)^m/m!)$, seems to imply that the original pair $\{w_m\}_{m\in\mathbb{N}_0}$ and $\{v_m\}_{m\in\mathbb{N}}$ outperforms the pair with thinning $\{\widetilde{w}_m\}_{m\in\mathbb{N}_0}$ and $\{\widetilde{v}_m\}_{m\in\mathbb{N}}$ as an iteration scheme, whereas this is a typical trade-off as the ones with thinning are generally easier to compute.
Also, the difference in convergence rate indicates, just as is standard in the context of Poisson thinning, that the constant rate $\widetilde{\lambda}$ should be chosen as small as possible with the most ideal value $\sup_{(t,{\bf x})\in [0,T]\times \overline{D}}\lambda(t,{\bf x})$ whenever it is available. 
For further illustration, one can show, under the regularity condition on $w_0$ of Theorem \ref{theorem uniform}, that, for all $(t,{\bf x})\in [0,T]\times \overline{D}$,
\[
w_m(t,{\bf x})=w_{m-1}(t,{\bf x}) + \mathbb{E}_\lambda^{t,{\bf x}}\left[\int_{\eta_t^T\land \tau_t^{(m-1)}}^{\eta_t^T \land \tau_t^{(m)}}\Theta_{t,s}N_0(s,X_s)ds\right],\quad
\widetilde{w}_m(t,{\bf x})=\widetilde{w}_{m-1}(t,{\bf x}) + \widetilde{\mathbb{E}}_{\widetilde{\lambda}}^{t,{\bf x}}\left[\int_{\eta_t^T\land \widetilde{\tau}_t^{(m-1)}}^{\eta_t^T\land \widetilde{\tau}_t^{(m)}}\Theta_{t,s}N_0(s,X_s)ds\right],
\]
as a direct consequence of the proof of Theorem \ref{theorem uniform} in the Appendix.
Looking at those two expressions together, a single step via the thinning version consists of the integration (of the common integrand $\Theta_{t,s}N_0(s,X_s)$) on an earlier and shorter interval, again due to its higher jump rate based on the dominance $\lambda(t,{\bf x})\le \widetilde{\lambda}$.   
For the most direct comparison, by setting $m=1$, we get
\[
w_1(t,{\bf x})=w_0(t,{\bf x}) + \mathbb{E}_\lambda^{t,{\bf x}}\left[\int_t^{\eta_t^T\land \tau_t^{(1)}}\Theta_{t,s}N_0(s,X_s)ds\right],\quad
\widetilde{w}_1(t,{\bf x})=\widetilde{w}_0(t,{\bf x}) + \widetilde{\mathbb{E}}_{\widetilde{\lambda}}^{t,{\bf x}}\left[\int_t^{\eta_t^T\land \widetilde{\tau}_t^{(1)}}\Theta_{t,s}N_0(s,X_s)ds\right],
\]
which implies, with $w_0=\widetilde{w}_0$ by definition, that the original one $\{w_m\}_{m\in\mathbb{N}_0}$ is naturally expected to make more progress by one step, in consistency with the difference in convergence rate in Theorem \ref{theorem uniform}.

If, additionally, the uniform ellipticity is imposed on the diffusion component (the condition (e) in Theorem \ref{theorem PDE equivalence} below), then the recursive relation \eqref{wmrecur2} under probability measure $\mathbb{P}^{t,{\bf x}}_0$ can further be put into the context of the Feyman-Kac formula.
To this end, we prepare some notation.
For $a=l+\alpha$ with $l\in\mathbb{N}_0$ and $\alpha\in (0,1]$ and an open set $B\in \mathbb{R}\times \mathbb{R}^d$, we denote by $\mathcal{H}_a(B)$ the Banach space of functions in $\mathcal{C}^{\lfloor l/2\rfloor,l}(B)$ having $l$-th space derivatives uniformly $\alpha$-H\"older continuous and $\lfloor l/2\rfloor$-time derivatives uniformly $(a/2-\lfloor l/2\rfloor)$-H\"older continuous.
Then, the following is a standard result on the existence of the solution to the partial differential equation \eqref{wmpde} and its equivalence to the probabilistic representation \eqref{wme0}.
We refer the reader to, for instance, \cite[Chapter 5]{doi:10.1142/3302}.

\begin{theorem}\label{theorem PDE equivalence}
Let $m\in \mathbb{N}$ and assume that there exists $\alpha\in (0,1)$ for the following conditions hold:
\begin{description}
\setlength{\parskip}{0cm}
\setlength{\itemsep}{0cm}
\item[(a)] The domain $D$ is non-empty convex subset of $\mathbb{R}^d$ with Lipschitz boundary $\partial D$.
\item[(b)]
The coefficients $b$ and $\sigma$ are in $\mathcal{H}_{1+\alpha}((0,T)\times D)$;
\item[(c)] The functions $g({\bf x})$ and $\Psi(t,{\bf x},{\bf x})$ are bounded and continuous in $(t,{\bf x})$ with $\Psi(T,{\bf x},{\bf x})=g({\bf x})$ for all ${\bf x}\in \partial D$;
\item[(d)] The functions $r(t,{\bf x})$, $\lambda(t,{\bf x})$, $\phi(t,{\bf x})$ and $G_{m-1}(t,{\bf x})+H(t,{\bf x})$ are bounded and in $\mathcal{H}_{\alpha}((0,T)\times D)$.
\item[(e)] 
There exists a strictly positive constant $c>0$ such that $\langle {\bm \xi},\sigma^{\otimes 2}(t,{\bf x}){\bm \xi}\rangle \ge c\|{\bm \xi}\|^2$ for all $(t,{\bf x},{\bm \xi})\in [0,T]\times \overline{D}\times \mathbb{R}^d$;
\end{description}
Then, there exists a unique solution in $\mathcal{C}^{1,2}((0,T)\times D)\cap \mathcal{C}^{0,0}([0,T]\times \overline{D})$ to the partial differential equation \eqref{wmpde} with initial and boundary conditions
$w_m(T,{\bf x})=g({\bf x})$ for ${\bf x}\in \overline{D}$, and $w_m(t,{\bf x})=\Psi(t,{\bf x},{\bf x})$ for $(t,{\bf x})\in[0,T)\times \partial D$.
\end{theorem}

We note that under the conditions of Theorem \ref{theorem PDE equivalence}, the function \eqref{udef} (under the probability measure $\mathbb{P}_{\lambda}^{t,{\bf x}}$ in the presence of jump component) solves the partial integro-differential equation
\[
 \partial_1 u(t,{\bf x}) + \mathcal{L}_t u(t,{\bf x}) 
 =r(t,{\bf x})u(t,{\bf x}) + \phi(t,{\bf x}) -\int_{\mathbb{R}_0^d} (u(t,{\bf x}+{\bf z})-u(t,{\bf x})) \lambda(t,{\bf x})\nu(d{\bf z};{\bf x}),
\]
or, equivalently and more closely to the present notation,
\begin{equation}
\partial_1 u(t,{\bf x}) + \mathcal{L}_t u(t,{\bf x}) = (r(t,{\bf x})+\lambda(t,{\bf x}))u(t,{\bf x}) + \phi(t,{\bf x}) -\int_{\{{\bf z}\in\mathbb{R}_0^d:\,{\bf x}+{\bf z}\in\overline{D}\}} u(t,{\bf x}+{\bf z})
    \lambda(t,{\bf x})\nu(d{\bf z};{\bf x})-H(t,{\bf x}), \label{pide u}
\end{equation}
with $u(T,{\bf x})=g({\bf x})$ for ${\bf x}\in \overline{D}$, due to the normalization condition $\int_{\mathbb{R}_0^d}\nu(d{\bf z};{\bf x})=1$ for ${\bf x}\in \overline{D}$ and $u(t,{\bf x})=\Psi(t,{\bf x},{\bf x})$ for $(t,{\bf x})\in [0,T)\times \partial D$.
This boundary value problem does not lend itself well to practical use, as partial integro-differential equations are generally not easy to deal with, whereas this fact provides an interesting insight into the proposed framework of iterative nature.
That is, the representation \eqref{wme0} can be thought of as a consequence of replacing the unknown true solution $u$ inside the integral by an approximate solution, say $w_{m-1}$, so as to get a better approximate solution $w_m$ by iteration.
The key point here is that the partial integro-differential equation loses its integral term, that is, then a partial differential equation, once this replacement has taken place, corresponding to the probability measure $\mathbb{P}_0^{t,{\bf x}}$ with jump component suppressed in the representation \eqref{wme0}. 

Before closing this section, let us stress again that the proposed framework has been constructed for quite a large class of multidimensional stochastic differential equations with jumps, based on the previous developments (Sections \ref{section first few jumps}, \ref{section two routes}, \ref{section E0} and \ref{section hard bounding functions}) alone, without coming into the present subsection (Section \ref{section uniformly bounded jump rate}), that is, without imposing uniform boundedness on jump rate or a variety of regularity conditions of Theorem \ref{theorem PDE equivalence}, to say nothing of differentiability on the drift and diffusion coefficients as in the exact simulation of one-dimensional sample paths \cite{casellaroberts2011, doi:10.1287/opre.2013.1191, GR2013, 10.3150/14-BEJ676}.

\section{Numerical illustrations}\label{section numerical illustrations}

In this section, we provide numerical examples to justify our theoretical findings, namely the recursive formula (Theorem \ref{theorem back to E0}), hard bounding functions (Theorem \ref{theorem when smooth}), thinning (Theorem \ref{theorem thinning}) and the convergence (Theorem \ref{theorem uniform}).
We devote the present section to two simple examples with a particular focus on effective illustration of the correctness of the proposed theoretical framework, rather than assessing competitiveness of the resulting numerical methodology in terms of complexity relative to the existing methods in the literature. 
As such, the two examples in what follows have been chosen in such a way that sophisticated numerical approximation is not required.

\subsection{Pure-jump processes}\label{simple poisson example}

The first example we consider is based on the pure-jump process:
\begin{equation}\label{path representation 0}
X_s=x+b(s-t)+ (N_s-N_t), \quad s\in [t,+\infty),\quad X_t=x(>0).
\end{equation}
where $b>0$ and $\{N_s: s\in [t,+\infty)\}$ is a Poisson process with constant rate $\lambda(>0)$ and a single jump size $c(<0)$, certainly in the framework \eqref{Xprocess} with $b(\cdot,\cdot)\equiv b$, $\sigma\equiv 0$ and $\nu(dz;\cdot)=\delta_{\{c\}}(dz)$.
With the diffusion component suppressed $(\sigma\equiv 0)$, we intend to demonstrate the fully deterministic iteration \eqref{fully deterministic case}.
We have adopted this example, where the jump rate is fully independent of time and state, so as to demonstrate most of our theoretical developments in the most illustrative manner.
Hence, the next example in Section \ref{section example impact} concerns jump-diffusion processes with a time-state dependent jump rate.

Now, with the domain $D=(0,\infty)$, we investigate the finite-time ruin probability before terminal time $T$:
\[
u(t,x)= \mathbb{P}_\lambda^{t,x}(\eta_t \leq T) = \mathbb{E}_\lambda^{t,x}\left[\mathbbm{1}(\eta_t\leq T)\right],\quad (t,x)\in [0,T]\times \overline{D},
\]
within the framework of Proposition \ref{proposition pure jump ruin probability} with $r(\cdot,\cdot)\equiv 0$, $g(\cdot)\equiv 0$, $\Psi(\cdot,\cdot,\cdot) \equiv 1$ and $\phi(\cdot,\cdot)\equiv 0$.
Hence, the sequence of approximate solutions $\{w_m\}_{m\in\mathbb{N}_0}$ is monotonically increasing.
For comparison purposes, we first derive the closed-form solutions of the first few approximate solutions with the aid of the expression for $w_m$ in Proposition \ref{proposition pure jump ruin probability}.
Since $b>0$ and $x>0$, it is immediate that $w_0(t,x)\equiv 0$ for $(t,x)\in [0,T]\times \overline{D}$.
For the case $m=1$, it holds that for $(t,x)\in [0,T]\times [0,+\infty)$,
\[
 w_1(t,x)=\mathbb{E}_{\lambda}^{t,x}\left[\mathbbm{1}\left(X_{\eta_t^T\land \tau_t^{(1)}}\notin [0,+\infty)\right)\right]
 =\mathbb{E}_{\lambda}^{t,x}\left[\mathbbm{1}(\tau_t^{(1)}\le T)\mathbbm{1}\left(X_{\tau_t^{(1)}}\notin [0,+\infty)\right)\right]
 =\mathbb{P}_{\lambda}^{t,x}\left(\tau_t^{(1)}-t\le \frac{-c-x}{b}\land (T-t)\right),
\]
and thus
\begin{equation}\label{w1 example}
 w_1(t,x)=\begin{cases}
 1-e^{-\lambda ((-c-x)/b\land (T-t))},& \text{if }(-c-x)/b\land (T-t)\ge 0,\\
 0,& \text{otherwise},
 \end{cases}
\end{equation}
due to $\tau_t^{(1)}-t\sim {\rm Exp}(\lambda)$.
For $m=2$, observe that
\[
 \mathbbm{1}\left(\{X_{T\land \tau_t^{(1)}}\in[0,\infty)\}\cap \{X_{T\land \tau_t^{(2)}}\notin[0,\infty)\}\right)
=\mathbbm{1}\left(\{t+E_1+E_2\le T\}\cap \{x+bE_1+c>0\}\cap \{x+b(E_1+E_2)+2c <0\}\right),
\]
where $E_1$ and $E_2$ are iid exponential random variables with rate $\lambda$.
The three events indicate, respectively, that the second jump occurs before the terminal time ($t+E_1+E_2\le T$), the location right after the first jump remains inside the domain ($x+bE_1+c>0$), and the location right after the second jump is strictly outside the domain ($x+b(E_1+E_2)+2c<0$).
Hence, we obtain that for $(t,x)\in [0,T]\times [0,+\infty),$
\begin{align}
w_2(t,x)
&=w_1(t,x)+\begin{dcases}
e^{-\lambda \theta_1(x)}-e^{-\lambda \theta_2(t,x)}-\lambda e^{-\lambda \theta_2(t,x)}\left(\theta_2(t,x)-\theta_1(x)\right),&\text{if } \theta_2(t,x)\ge \theta_1(x),\\
0,&\text{otherwise},
\end{dcases}\label{w2 example}
\end{align}
where $\theta_1(x):=\max\{0,(-c-x)/b\}$ and $\theta_2(t,x):=\max\{0,(-2c-x)/b\land (T-t)\}$.
Note that $w_2(t,x)\ge w_1(t,x)$ for all $(t,x)\in [0,T]\times [0,+\infty)$, since the function $e^{-x}-e^{-y}-e^{-y}(y-x)$ is non-negative for $(x,y)\in [0,+\infty)^2$ with $x\le y$.
It is evidently cumbersome to go beyond $m=2$, whereas it is natural to conjecture, based on those closed-form expressions $w_0\equiv 0$, \eqref{w1 example} and \eqref{w2 example}, that the function $w_m$ is continuous with only a countable number of non-differentiable points for all $m\in\mathbb{N}_0$.

Alternatively, we may derive a recursion formula based on the identity \eqref{fully deterministic case}.
In the absence of the jump component under the probability measure $\mathbb{P}_0^{t,{\bf x}}$, the underlying process \eqref{path representation 0} reduces to a deterministic linear function $x(s)=x+b(s-t)$.
The deterministic function $x(s)$ here never exits the domain $(0,+\infty)$ due to $b>0$, that is, the first exit time $\eta_t$ is understood to be infinite.
With constant rate $\lambda$, $\nu(dz;\cdot)=\delta_{\{c\}}(dz)$ and $\Psi\equiv 1$, we get 
\[
 G_{m-1}(s,x(s))=\lambda w_{m-1}(s,x+b(s-t)+c)\mathbbm{1}(x+b(s-t)+c\geq 0),\quad H(s,x(s))=\lambda \mathbbm{1}(x+b(s-t)+c<0),
\]
for $s\in [t,T]$.
Therefore, with $r\equiv 0$, $g\equiv 0$, $\phi\equiv 0$ and $\Psi\equiv 1$, we get $w_0\equiv 0$, and for each $m\in\mathbb{N}$,
\begin{equation}\label{iteration formula in example}
w_m(t,x)=\int_t^T \lambda e^{-\lambda(s-t)}\left(w_{m-1}(s,x+b(s-t)+c)\mathbbm{1}(x+b(s-t)+c\geq 0)+\mathbbm{1}(x+b(s-t)+c<0) \right)ds,
\end{equation}
according to the result \eqref{fully deterministic case}.
Note that the closed-form expressions \eqref{w1 example} and \eqref{w2 example} are consistent with the recursion formula \eqref{iteration formula in example}.
Recalling $w_m=v_m$ due to $w_0\equiv 0$, as we have pointed out in Section \ref{section two routes}, there is nothing really to illustrate Theorem \ref{theorem vm wm} here.

As conjectured earlier, one can now ensure, by induction on the recursion \eqref{iteration formula in example} starting with $w_0\equiv 0$, that for every $m\in \mathbb{N}_0$, the function $w_m$ is continuous with an increasing yet only countable number of non-differentiable points on $[0,T)\times (0,+\infty)$.
Moreover, we have $\xi(t,x;0)=T-t$, $M^U(t)=0$ and $M^L(t)=-\lambda (T-t)$, due to $\eta_t=+\infty$, $r\equiv 0$ and $M(t,x)=\lambda(T-t)(\mathbbm{1}(x+c\ge 0)-1)$. 
Next, observe that
\[
N_m(t,x)=G_m(t,x)-G_{m-1}(t,x) = \lambda  (w_m(t,x+c)-w_{m-1}(t,x+c))\mathbbm{1}(x+c>0),
\]
which is non-negative, that is, $N_m^L(t)\equiv 0$ for all $t\in[0,T]$ and $m\in\mathbb{N}$, since the sequence $\{w_m\}_{m\in\mathbb{N}_0}$ is monotonically increasing with $w_0\equiv 0$ by Proposition \ref{proposition pure jump ruin probability}. 
One can also deduce $N_m^U(t)=\lambda w_m(t,-(m-1)c)$ and $|N_m^U(0)|+|N_m^L(0)|<+\infty$ for all $m\in\mathbb{N}_0$ since $\{w_m\}_{m\in\mathbb{N}}$ is bounded.
That is, this problem setting lies in the scope of Theorems \ref{theorem when smooth} and \ref{theorem uniform}, with deterministic upper and lower bounds given, through \eqref{boundswm0}, by
\begin{equation}\label{UL example 1}
w_m(t,x)\leq u(t,x) \leq w_m(t,x)+\lambda (T-t) w_m(t,-(m-1)c).
\end{equation}

\begin{figure}[ht]
  \centering
  \begin{subfigure}[b]{0.48\linewidth}
    \includegraphics[width=\linewidth]{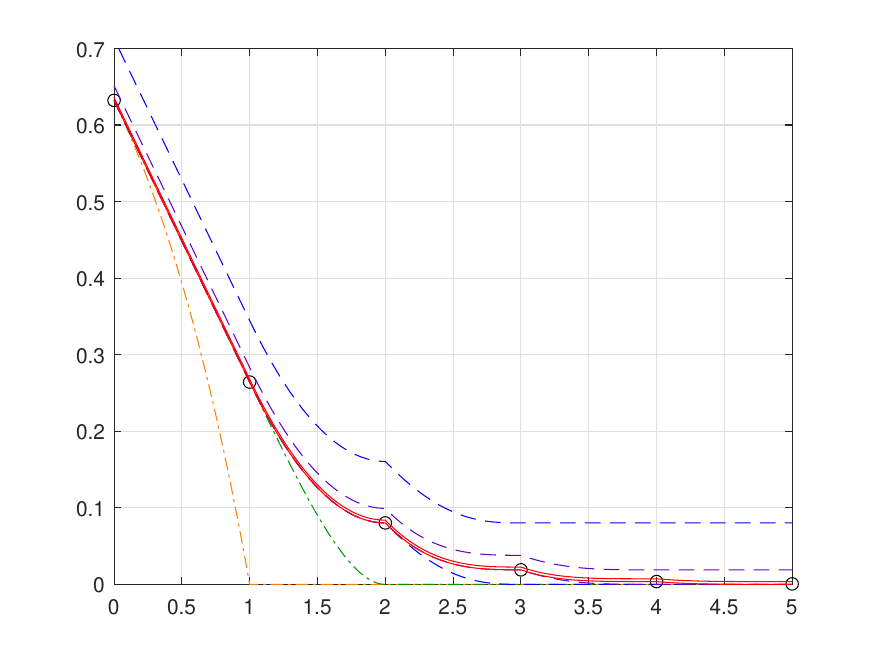}
    \caption{$\{w_m(0,x)\}_{m\in\{1,2,3,4,5\}}$ and upper  bounding functions}
    \end{subfigure}
      \begin{subfigure}[b]{0.48\linewidth}
    \includegraphics[width=\linewidth]{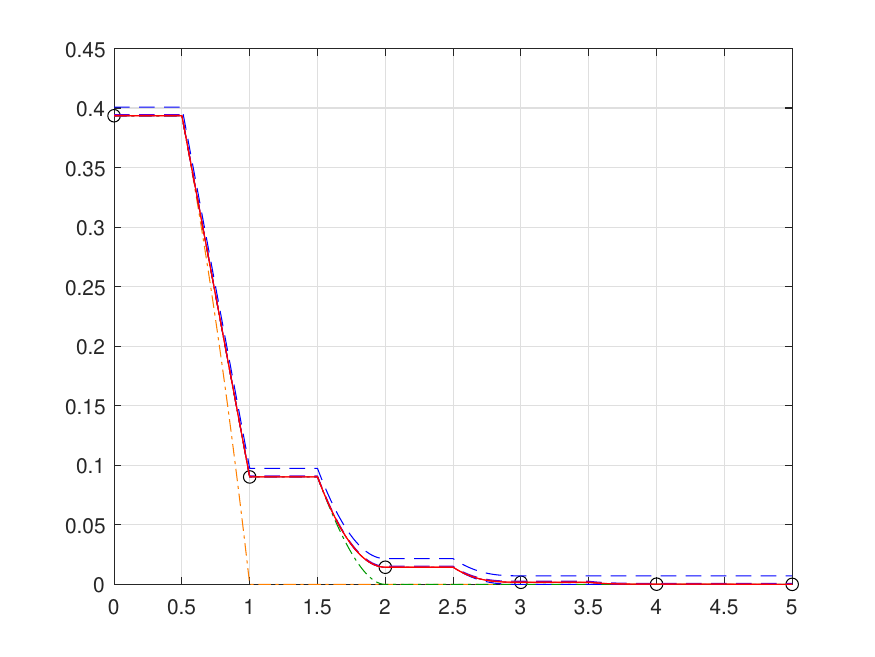}
    \caption{$\{w_m(0.5,x)\}_{m\in\{1,2,3,4,5\}}$ and upper bounding functions}
    \end{subfigure}
 \caption{Plots of the approximate solutions $\{w_m\}_{m\in\{1,2,3,4,5\}}$ at $t\in \{0,\,0.5\}$ and the associated upper bounding functions, according to \eqref{UL example 1}, at the third, fourth and fifth iterations, with $m=1$ (orange dash-dot), $m=2$ (green dash-dot), $m=3$ (blue dash), $m=4$ (purple dash) and $m=5$ (red solid).
  The six black unfilled circles in each figure indicate (very accurate Monte Carlo estimates of) the true values $u(\cdot,x)$ at $x=\{0,1,2,3,4,5\}$.}
\label{fig01}
\end{figure}

In what follows, with the unit drift $b=1$ fixed, we examine two distinctive situations; large yet rare jumps and small yet intense jumps.
For the first case of large yet rare jumps, with jump rate $\lambda=1$ and the drawdown jump $c=-1$ of unit length, we plot in Figure \ref{fig01} the first five approximate solutions and the associated upper bounding functions. 
To avoid overloading the figures, the upper bounding functions are only presented at the third, fourth and fifth iterations.
The values of the true solution $u$ indicated by unfilled circles at six points in each figure are estimates by Monte Carlo methods using $10^7$ iid replications, which are large enough to narrow down 99\% confidence intervals to almost a singleton.  

As the iteration proceeds, the approximate solutions are supporting the ruin probability gradually from the lower initial state.
This is not very surprising, as the present problem setting is built on a singular structure (positive linear drift with fatal drawdown jumps).
Even after an approximate solution has captured the true solution from below (around $m=4$), the upper bound is still approaching downwards to the true solution and specifies the true solution, almost exactly and uniformly when $m=5$, just as expected in Theorem \ref{theorem uniform}.

\begin{figure}[ht]
  \centering
    \begin{subfigure}[b]{0.45\linewidth}
    \includegraphics[width=\linewidth]{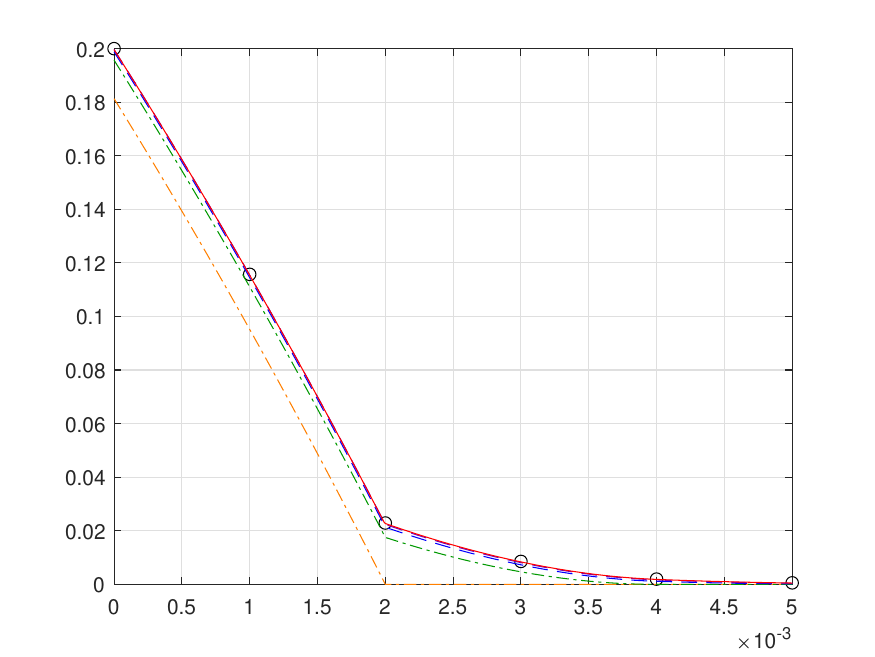}
    \caption{$c=-1/500$}
    \end{subfigure}
  \begin{subfigure}[b]{0.45\linewidth}
    \includegraphics[width=\linewidth]{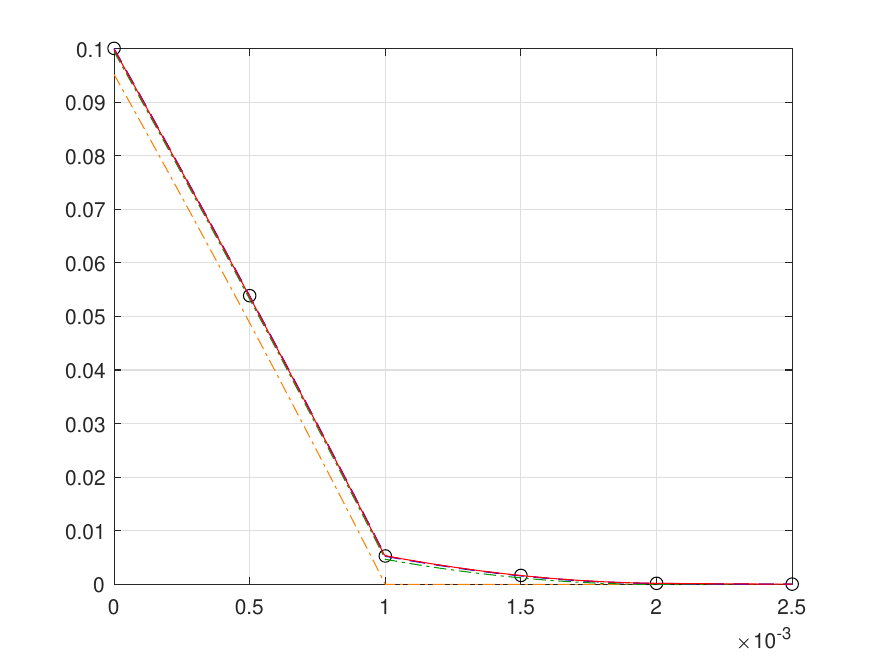}
    \caption{$c=-1/1000$}
    \end{subfigure}
 \caption{Plots of the approximate solutions $\{w_m(0,x)\}_{m\in\{1,2,3,4,5\}}$, with $m=1$ (orange dash-dot), $m=2$ (green dash-dot), $m=3$ (blue dash), $m=4$ (purple dash) and $m=5$ (red solid).
  The six black unfilled circles in each figure indicate (very accurate Monte Carlo estimates of) the true values $u(0,x)$ at six states.}
  \label{fig03}
\end{figure}

Next, for the case of small yet intense jumps, we plot in Figure \ref{fig03} the approximate solutions with intense jump rate $\lambda=100$ and small drawdown jumps of (a) $c=-1/500$ and (b) $c=-1/1000$, in such a way that each jump is nearly invisible relative to the ascending drift $b=1$.
In this situation, the associated upper bounding functions are not useful (and thus not plotted in the figures) because the jump rate makes direct impact as is clear in \eqref{UL example 1}. 
Still, it is encouraging that the proposed framework provides a fast convergent sequence of the approximate solutions in both distinctive situations of large yet rare jumps and small yet intense jumps. 

\subsection{Survival probability with time-state dependent jump rate}\label{section example impact}

The previous example (Section \ref{simple poisson example}) has concerned pure-jump processes on the basis of, in fact, already most of the presented theoretical results (Sections \ref{section first few jumps}, \ref{section two routes}, \ref{section E0} and \ref{section hard bounding functions}).
We next turn to the case of jump-diffusion processes with a time-state dependent jump rate, so as to justify the remaining developments, particularly, the Poisson thinning of Section \ref{section uniformly bounded jump rate}.
Note that this example is borrowed from \cite[Example 2.2.1]{doi:10.1137/151004070}, as semi-analytical expressions are available for the target solution, in order to avoid getting off the main track.

Consider a bounded space domain $D=(x_L,x_U)$ with $-\infty<x_L<x_U<\infty$.
The underlying stochastic process is a Brownian motion with drift and jumps:
\[
dX_s=b ds+\sigma dW_s + \int_{\mathbb{R}_0} z \mu(dz, ds; X_{s-}),\quad s\in [t,+\infty),\quad X_t=x,
\]
where $b\in\mathbb{R}$, $\sigma>0$ and $\{W_s:\, s\in[t,+\infty)\}$ is the standard Brownian motion in $\mathbb{R}$.
We let the jump rate time-state dependent with its jump size distributed Gaussian as follows: 
\[
\lambda(t,x)=5t(T-t)(x_U-x)(x-x_L), \quad
\nu(dz;x)=\frac{1}{\sqrt{2\pi \rho}}\exp\left(-\frac{z^2}{2\rho}\right) dz,\quad 
(t,x,z)\in [0,T]\times (x_L,x_U)\times \mathbb{R}_0,
\]
for some $\rho >0$.
We now examine the survival probability until the terminal time $T$, that is, $u(t,x)=\mathbb{P}_\lambda^{t,x}\left(\eta_t>T\right)$, corresponding to the input data $(r,g,\Psi,\phi)=(0,1,0,0)$ in the representation \eqref{udef}.
The initial approximate solution is available in semi-analytical form  \cite[Example 2.2.1]{doi:10.1137/151004070}, as
\begin{equation}\label{w0 infinite}
w_0(t,x) = \mathbb{P}_0^{t,x}\left(\eta_t>T\right)
=e^{-\frac{b}{\sigma^2}x}\sum_{k\in \mathbb{N}}\alpha_k\sin\left(k\pi\frac{x-x_L}{x_U-x_L}\right)e^{-\left(
\frac{b^2}{\sigma^2}+\left(\frac{k\pi\sigma}{x_U-x_L}\right)^2\right)\frac{(T-t)}{2}},
\end{equation}
where
\[
\alpha_k:=\frac{2k\pi\sigma^2}{b^2(x_U-x_L)^2+k^2\pi^2\sigma^4}\left(e^{\frac{b}{\sigma^2}x_L}-(-1)^k e^{\frac{b}{\sigma^2}x_U}\right),\quad k\in \mathbb{N}.
\]
In order to apply the Poisson thinning for further iterations, we may find and thus set the smallest upper bound in this case, that is,  $\widetilde{\lambda}=\sup_{(t,x)\in[0,T]\times[x_L,x_U]}\lambda(t,x)=(5/16)T^2 (x_U-x_L)^2$.
With the input data $(r,g,\Psi,\phi)=(0,1,0,0)$, the identity \eqref{defwmthinning} under the probability measure $\mathbb{P}_0^{t,x}$ reduces to 
\[
\widetilde{w}_m(t,x)=\mathbb{E}_0^{t,x} \left[\mathbbm{1}(\eta_t>T)  \widetilde{\Lambda}_{t,T} +\int_t^{\eta_t\land T} \widetilde{\Lambda}_{t,s}  \widetilde{G}_{m-1}(s,X_s)  ds\right],
\]
for all $(t,x)\in[0,T]\times[x_L,x_U]$ and $m\in\mathbb{N}$, where
\begin{equation}
\widetilde{G}_{m-1}(t,x)=(\widetilde{\lambda}-\lambda(t,x))\widetilde{w}_{m-1}(t,x)+\lambda(t,x) \int_{x_L-x}^{x_U-x} \widetilde{w}_{m-1}(t,x+z) \frac{1}{\sqrt{2\pi \rho}}\exp\left(-\frac{z^2}{2\rho}\right) dz.
\end{equation}
Its semi-analytical solution is available as
\begin{equation}\label{wm infinite}
\widetilde{w}_m(t,x)=e^{-\widetilde{\lambda}(T-t)}w_0(t,x)+e^{-\frac{b}{\sigma^2}x}\sum_{k\in \mathbb{N}}\sin\left(k\pi\frac{x-x_L}{x_U-x_L}\right)\int_t^{T} \exp\left[
-\left(
\frac{1}{2}\left(\frac{k\pi\sigma}{x_U-x_L}\right)^2+\frac{b^2}{2\sigma^2}+\widetilde{\lambda}\right)(s-t)\right]\beta_{m,k}(s) ds,
\end{equation}
for all $(t,x)\in[0,T]\times[x_L,x_U]$ and $m\in\mathbb{N}$, where
\[
\beta_{m,k}(s):=\frac{2}{x_U-x_L}\int_{x_L}^{x_U}e^{\frac{b}{\sigma^2}y} \widetilde{G}_{m-1}(s,y)\sin\left(k\pi\frac{y-x_L}{x_U-x_L}\right)dy,
\]
for $s\in [0,T]$ and $m\in \mathbb{N}$.
Finally, for constructing the hard bounding functions \eqref{boundswm02}, the required component $\xi(t,x;0)=\mathbb{E}_0^{t,x}[\eta_t^T]-t$  is also available in semi-analytical form \cite[Example 2.2.1]{doi:10.1137/151004070}, as
\begin{equation}\label{xi infinite}
\xi(t,x; 0)=h_0(x)-e^{-\frac{b}{\sigma^2}x}\sum_{k\in\mathbb{N}}c_k\sin\left(k\pi\frac{x-x_L}{x_U-x_L}\right)e^{-\left(\frac{b^2}{\sigma^2}+\left(\frac{k\pi\sigma}{x_U-x_L}\right)^2\right)\frac{(T-t)}{2}},
\end{equation}
for all $(t,x)\in[0,T]\times[x_L,x_U]$, where
\[
h_0(x):=\frac{(x_U-x)(e^{-2b(x_L-x)/\sigma^2}-1)+(x_L-x)(1-e^{-2b(x_U-x)/\sigma^2})}{b(e^{-2b(x_L-x)/\sigma^2}-e^{-2b(x_U-x)/\sigma^2})},
\qquad
c_k:=\frac{2}{x_U-x_L}\int_{x_L}^{x_U}h_0(y)e^{\frac{b}{\sigma^2}y}\sin\left(k\pi\frac{y-x_L}{x_U-x_L}\right)dy,
\]
for $x\in [x_L,x_U]$ and $k\in\mathbb{N}$.

In consistency with \cite[Example 2.2.1]{doi:10.1137/151004070}, we set $T=1$, $(x_L,x_U)=(0,2)$, $b=2$, $\sigma=1$ and $\rho=0.1$, and truncate the infinite series \eqref{w0 infinite}, \eqref{wm infinite} and \eqref{xi infinite} by 500 summands to ensure adequate convergences on the computer.
When constructing hard bounding functions according to the inequalities \eqref{boundswm02}, we find the supremum and infimum of $\widetilde{N}_m(\cdot,\cdot)$ and $M(\cdot,\cdot)$ approximately by discretization of the time-state domain $[0,T]\times [x_L,x_U]$ into equisized grids, each of $0.0005\times0.0005$.

\begin{figure}[ht]
  \centering
  \begin{subfigure}[b]{0.45\linewidth}
    \includegraphics[width=\linewidth]{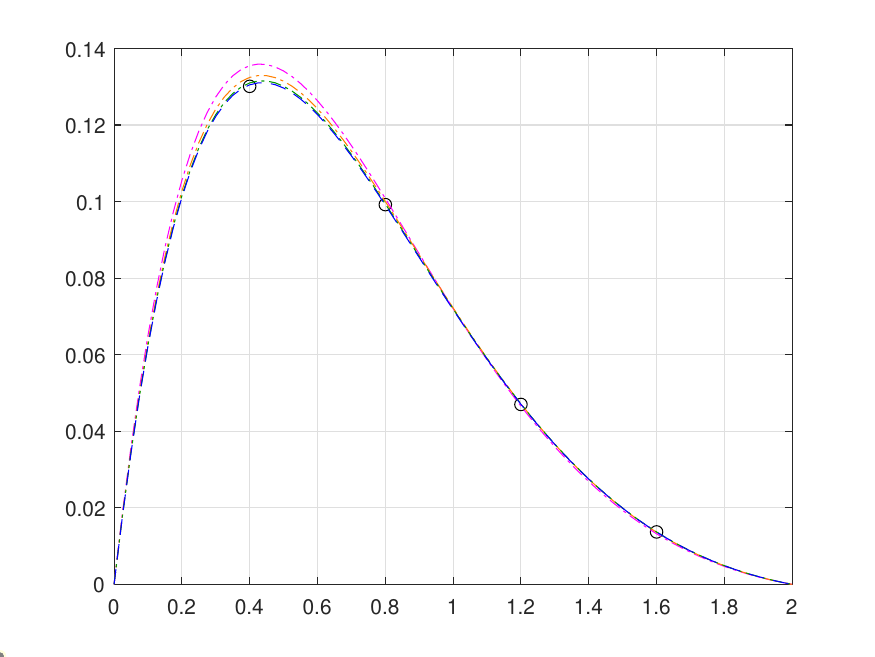}
    \caption{Approximate solutions $\{\widetilde{w}_m(0,x)\}_{m\in\{0,1,2,3\}}$}
    \end{subfigure}
     \begin{subfigure}[b]{0.45\linewidth}
    \includegraphics[width=\linewidth]{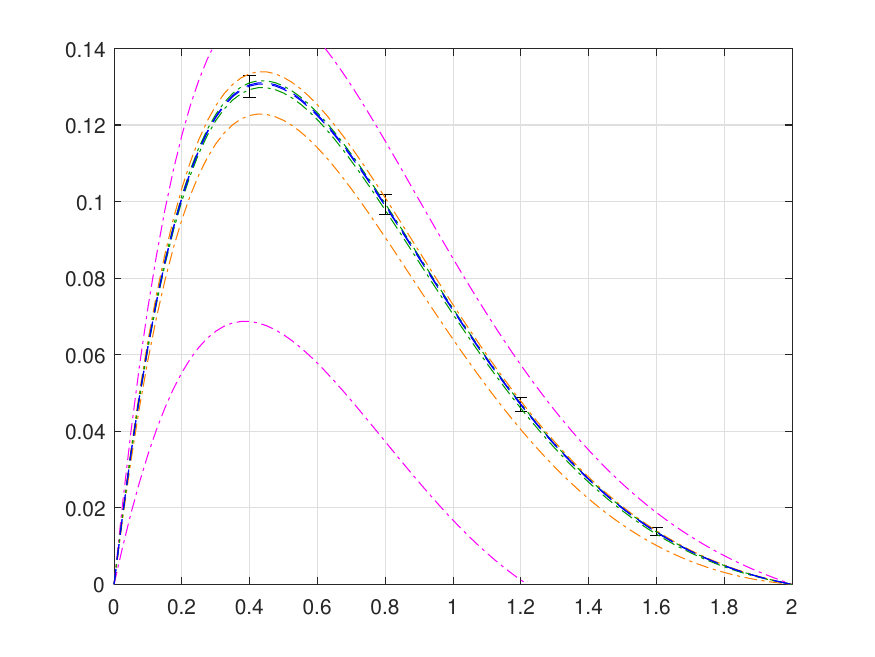}
    \caption{Hard bounding functions for $u(0,x)$ at $m\in \{0,1,2,3\}$}
    \end{subfigure}
      \begin{subfigure}[b]{0.45\linewidth}
    \includegraphics[width=\linewidth]{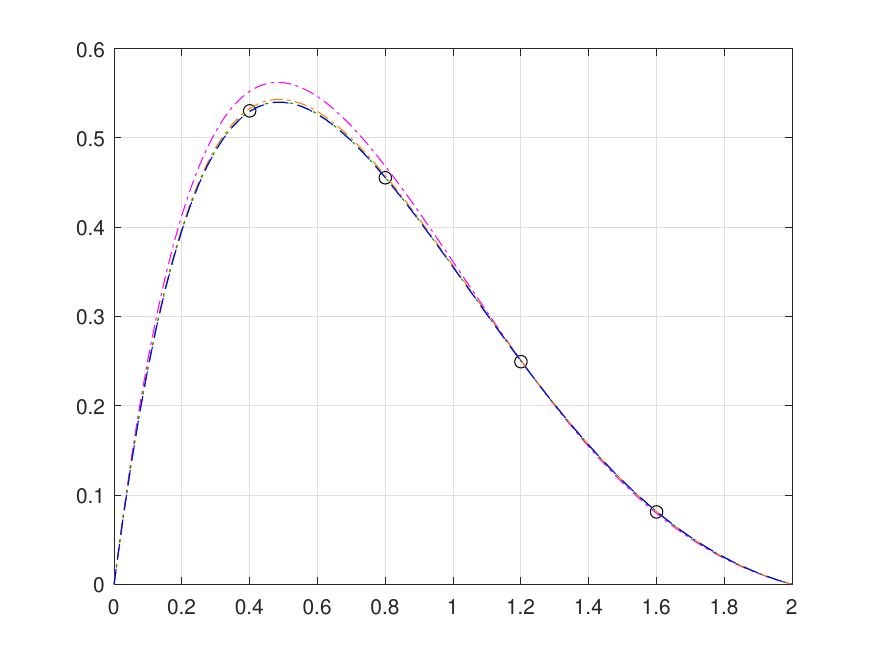}
    \caption{Approximate solutions $\{\widetilde{w}_m(0.5,x)\}_{m\in\{0,1,2,3\}}$}
    \end{subfigure}
      \begin{subfigure}[b]{0.45\linewidth}
    \includegraphics[width=\linewidth]{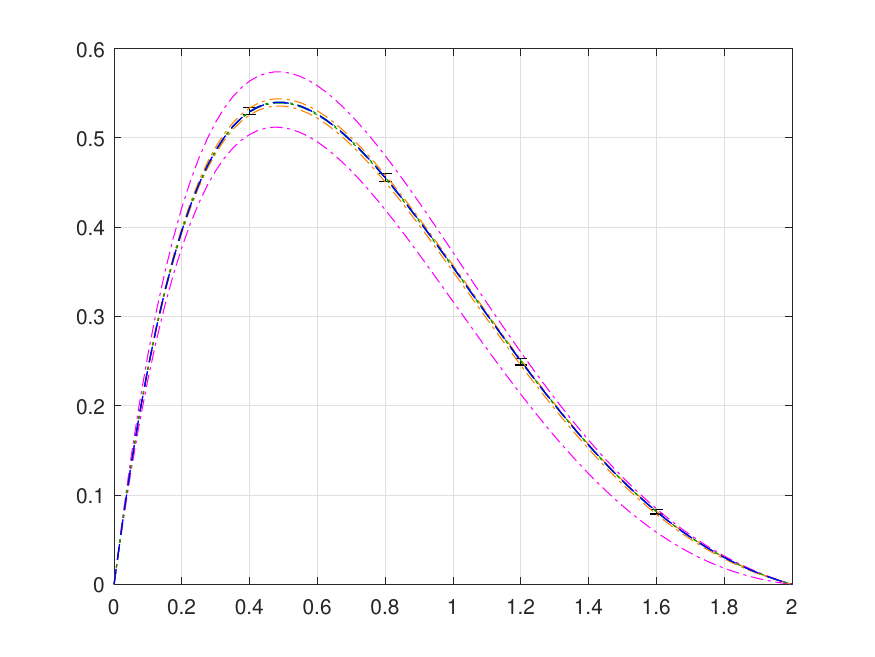}
    \caption{Hard bounding functions for $u(0.5,x)$ at $m\in \{0,1,2,3\}$}
    \end{subfigure}
 \caption{The approximate solutions $\{w_m(t,x)\}_{m\in\{0,1,2,3\}}$ and the associated upper and lower bounding functions at two timepoints $t\in \{0,\,0.5\}$ for $x\in [x_L,x_U]$ with $m=0$ (pink dash-dot), $m=1$ (orange dash-dot), $m=2$ (green dash) and $m=3$ (blue dash).
 The unfilled circles in (a) and (c) and the vertical whisker plots in (b) and (d) are, respectively, Monte Carlo estimates and 99\% confidence intervals, constructed based on $10^5$ iid sample paths by the exact simulation method for jump-diffusion processes \cite{casellaroberts2011}.}
\label{fig02}
\end{figure}

In Figure \ref{fig02}, we plot the approximate solutions $\{\widetilde{w}_m(t,x)\}_{m\in \{0,1,2,3\}}$ as well as the associated upper and lower hard bounding functions at $t\in \{0,\,0.5\}$ for $x\in [x_L,x_U]$.
For justification and comparison purposes, we add Monte Carlo estimates (as unfilled circles) in (a) and (c), while 99\% confidence intervals (as vertical whisker plots) in (b) and (d), at $x\in \{0.4, 0.8, 1.2, 1.6\}$, all constructed based on $10^5$ iid sample paths by the exact simulation method for jump-diffusion processes \cite{casellaroberts2011}.
Those estimates and confidence intervals can certainly be improved by, for instance, increasing the sample path and/or employing relevant variance reduction methods, whereas we do not go in such directions here as we are not particularly concerned with a direct competition with existing methodologies in computing time.
Instead, hard bounding functions (so to speak, 100\% confidence intervals) are sufficiently competitive with those 99\% confidence intervals already at the second ($m=2$) or third ($m=3$) iteration.
Let us, last but not least, remind that the proposed framework is general enough to accommodate quite a large class of multidimensional inhomogeneous stochastic differential equations with jumps and, particularly, does not require the drift and diffusion coefficients to be as smooth as in existing univariate exact simulation methods.

To close this work, we make some concluding remarks.
We have investigated rather simple examples with semi-analytical expressions available in order to avoid digressing into external numerical techniques and peripheral issues on the way.
Still, it has been demonstrated effectively that the proposed recursive representation and associated hard bounding function have the potential to act as a promising numerical method, possibly in high-dimensional problems as well, for instance, by appealing to the emerging deep learning type methods.
In principle, the main question would then be whether iterating the standard formation (for instance, PDEs without solution inside a integral term, like \eqref{wmpde}) is more advantageous than approximating an intricate formulation (for instance, PIDEs, like \eqref{pide u}) in high dimensions.
With convergence and its rate already available (Theorems \ref{theorem convergence}, \ref{theorem vm wm}, \ref{theorem thinning} and \ref{theorem uniform}), such numerical aspects would deserve separate attention.
Finally, the proposed framework may be extended to deal with backward stochastic differential equations with jumps, which generally require heavy computation \cite{BSDEsurvey}.
All those topics would be interesting future directions of research deserving of their own investigation.

\section{Proofs}

\footnotesize

To avoid overloading the paper, we skip nonessential details of somewhat routine nature in some instance, particularly when the derivation proceeds using the tower property with respect to the $\sigma$-field associated with a stopping time.
To define this, we denote by $(\Omega,\mathcal{F},(\mathcal{F}_t)_{t\ge 0},\mathbb{P})$ a filtered probability space satisfying the usual conditions to which the underlying process $\{X_t:\,t\ge 0\}$ is adapted.
For an $(\mathcal{F}_t)_{t\ge 0}$-stopping time $\xi$, we let $\mathcal{F}_{\xi}$ denote the stopped $\sigma$-field on the underlying process, that is the $\sigma$-field of events up to its stopping time $\xi$, defined by $\mathcal{F}_{\tau}:=\{B\in \mathcal{F}:\,B\cap \{\xi\le t\}\in \mathcal{F}_t \text{ for all }t\ge 0\}$.

\begin{proof}[Proof of Theorem \ref{theorem convergence}]
We use the identities
\begin{equation}
      \{\eta_t>T\}\cap \{X_T \not\in \overline{D}\}= \varnothing,\quad \{\eta_t>T\}\cap \{X_T\in \overline{D}\}= \{\eta_t>T\}.\label{p307}
\end{equation}
for $T>t$.
In words, on the event $\{\eta_t>T\}$, that is, no exit occurs before or at the terminal time $T$, the location $X_T$ cannot be strictly outside the closure $\overline{D}$.
The second identity is clearly the complement of the first one. 

Now, if the first exit is caused by a jump, that is, $X_{\eta_t}\neq X_{\eta_t-}$, then the location $X_{\eta_t}$ right after the jump is strictly outside the closure $\overline{D}$ almost surely.
If the first exit is caused by the diffusion component of infinite variation, then the exit must occur in the boundary $\partial D$ due to infinite variation, that is, that is, $X_{\eta_t}=X_{\eta_t-}\in \partial D$.
Since those two cases are disjoint, we obtain
\begin{equation}\label{p309}
\{X_{\eta_t}\neq X_{\eta_t-}\} = \{X_{\eta_t}\not\in\overline{D}\}, \quad \{X_{\eta_t}=X_{\eta_t-}\}=\{X_{\eta_t}\in\partial D\}.
\end{equation}
With the aid of \eqref{p307} and \eqref{p309}, we get
\begin{align*}
\mathbbm{1}(\{X_{\eta_t^T}\in \overline{D}\}\cap \{X_{\eta_t}\neq X_{\eta_t-}\})
        &=\mathbbm{1} (\{\eta_t>T\}\cap \{X_T\in \overline{D}\}\cap \{X_{\eta_t}\neq X_{\eta_t-}\})
        +\mathbbm{1} (\{\eta_t\le T\}\cap \{X_{\eta_t}\in \overline{D}\}\cap \{X_{\eta_t}\neq X_{\eta_t-}\})\\
        &=\mathbbm{1}(\{\eta_t>T\}\cap \{X_{\eta_t}\neq X_{\eta_t-}\}),\\
\mathbbm{1}(\{X_{\eta_t^T} \not\in \overline{D}\}\cap \{X_{\eta_t}\neq X_{\eta_t-}\})
        &=\mathbbm{1} (\{\eta_t>T\}\cap \{X_T \not\in \overline{D}\}\cap \{X_{\eta_t}\neq X_{\eta_t-}\})
        +\mathbbm{1} (\{\eta_t\le T\}\cap \{X_{\eta_t} \not\in \overline{D}\}\cap \{X_{\eta_t}\neq X_{\eta_t-}\})\\
        &=\mathbbm{1}(\{\eta_t\le T\}\cap \{X_{\eta_t} \not\in \overline{D}\}),
\end{align*}
and, in a similar manner,
\begin{gather*}
        \mathbbm{1} (\{X_{\eta_t^T}\in \overline{D}\}\cap \{X_{\eta_t}= X_{\eta_t-}\})
        =\mathbbm{1} (\{\eta_t>T\}\cap \{X_{\eta_t}= X_{\eta_t-}\})+\mathbbm{1} (\{\eta_t\le T\}\cap \{X_{\eta_t}\in \partial D \}),\\
        \mathbbm{1} (\{X_{\eta_t^T} \not\in \overline{D}\}\cap \{X_{\eta_t}= X_{\eta_t-}\})
        =\mathbbm{1} (\{\eta_t\le T\}\cap \{X_{\eta_t} \not\in\overline{D} \}).
\end{gather*}
Combining those altogether, we obtain the identities
\begin{equation}
\mathbbm{1} (X_{\eta_t^T}\in \overline{D})=\mathbbm{1} (\eta_t>T)+\mathbbm{1} (\{\eta_t\le T\}\cap \{X_{\eta_t}\in \partial D\}),\quad 
\mathbbm{1}(X_{\eta_t^T}\notin \overline{D})=\mathbbm{1} (\{\eta_t\le T\}\cap \{X_{\eta_t}\notin \overline{D}\}).\label{p311}
\end{equation}

Since the rate function $\lambda(t,\cdot)$ is locally Lipschitz for all $t\in [0,T]$ and the intensity measure $\nu(d\cdot;{\bf x})$ is finite for all ${\bf x}\in \overline{D}$, the number of jumps on every time interval of finite length is almost surely finite.
It thus holds that $\tau_t^{(m)}>T$ eventually as $m\to +\infty$ under the probability measure $\mathbb{P}_{\lambda}^{t,{\bf x}}$.
Since $g$, $\Psi$, $\phi$ and $\Theta$ are bounded and $\eta_t^T\in [t,T]$, the function $w_0$ is bounded due to the representation \eqref{w0def1}.
Hence, for every $m\in \mathbb{N}$, the random element in the expectation \eqref{wmdef2} is almost surely bounded.
Hence, with the aid of the bounded convergence theorem, it holds, as $m\to +\infty$, that for each $(t,{\bf x})\in [0,T]\times \overline{D}$,
\begin{align*}
w_m(t,{\bf x})  &\to \mathbb{E}_{\lambda}^{t,{\bf x}} \left[\mathbbm{1}(
X_{\eta_t^T}\in \overline{D}) \Theta_{t,\eta_t^T} w_0(\eta_t^T, X_{\eta_t^T}) + \mathbbm{1} (X_{\eta_t^T}\not\in \overline{D} ) \Theta_{t,\eta^T_t} \Psi(\eta^T_t, X_{\eta^T_t},X_{\eta^T_t-})
-\int_t^{\eta_t^T} \Theta_{t,s} \phi(s, X_s) ds\right]\\
&=\mathbb{E}_{\lambda}^{t,{\bf x}} \Bigg[\mathbbm{1} \left(\eta_t>T\right)\Theta_{t,T}w_0(T,X_T) + \mathbbm{1} \left(\{\eta_t\leq T\}\cap\{X_{\eta_t}\in\partial D\}\right)\Theta_{t,\eta^T_t}w_0(\eta_t,X_{\eta_t})
\nonumber\\
&\qquad\qquad\qquad +\mathbbm{1} (\{\eta_t\le T\}\cap \{X_{\eta_t}\notin\overline{D}\}) \Theta_{t,\eta^T_t} \Psi(\eta^T_t, X_{\eta^T_t},X_{\eta^T_t-})-\int_t^{\eta_t^T} \Theta_{t,s} \phi(s, X_s) ds\Bigg],\nonumber
\end{align*}
which yields \eqref{udef}, due to $w_0(t,{\bf x})=g({\bf x})$ for ${\bf x}\in\overline{D}$ and $w_0(\eta_t,X_{\eta_t})=\Psi(\eta_t, X_{\eta_t},X_{\eta_t-})=\Psi(\eta_t, X_{\eta_t},X_{\eta_t})$ on the event $\{X_{\eta_t}\in \partial D\}$.
\end{proof}

\begin{proof}[Proof of Theorem \ref{theorem wmrecurrsive}]
First, by the definition \eqref{etadef}, we have 
\begin{equation}\label{etaeta is eta}
\eta_{\eta_t} =\inf\left\{s>\eta_t: X_s\not\in\overline{D}\right\} = \eta_t,
\end{equation}
because the trajectory has already hit the boundary $\partial D$ by the diffusion component, if any, at time $\eta_t$, or has already exited from the closure $\overline{D}$ by a jump at time $\eta_t$. 
Next, it is straightforward that 
\begin{equation}\label{etaTetaT is etaT}
  \eta_{\eta_t^T}^T= \eta_{\eta_t \land T} \land T = \left(\eta_{\eta_t} \land T\right)\mathbbm{1}\left(\eta_t\le T\right)+\left(\eta_T \land T\right)\mathbbm{1}\left(\eta_t> T\right)
  =\eta_t \mathbbm{1}\left(\eta_t\le T\right)+T\mathbbm{1}\left(\eta_t> T\right)
  =\eta_t^T,
\end{equation}
where we have applied \eqref{etaeta is eta} for the second equality.
Noting $\tau^{(m)}_t \geq t$ for all $m\in\mathbb{N}$ and $t\ge 0$, due to \eqref{tau1} and \eqref{taum}, we obtain $\eta_{\eta_t^T}^T \land \tau_{\eta_t^T}^{(m)} = \eta_t^T\land \tau_{\eta_t^T}^{(m)} =\eta_t^T$, where we have applied \eqref{etaTetaT is etaT} for the first equality.
In light of the expression \eqref{wmdef2}, it holds that for every $m\in\mathbb{N}$ and $n\in \{0,1,2,\cdots,m\}$,
\begin{align}
&\mathbbm{1}(X_{\eta_t^T}\in\overline{D})w_{m-n}(\eta_t^T,X_{\eta_t^T})\nonumber\\
&\qquad = \mathbbm{1}(X_{\eta_t^T}\in\overline{D})\mathbb{E}_{\lambda}^{\eta_t^T,X_{\eta_t^T}}\left[ \mathbbm{1}(X_{\eta_t^T}\in\overline{D})\Theta_{\eta_t^T,\eta_t^T}w_0(\eta_t^T,X_{\eta_t^T})
+\mathbbm{1}(X_{\eta_t^T}\not\in\overline{D})\Theta_{\eta_t^T,\eta_{\eta_t^T}}\Psi(\eta_{\eta_t^T},X_{\eta_{\eta_t^T}},X_{\eta_{\eta_t^T}-})- \int_{\eta_t^T}^{\eta_t^T}\Theta_{\eta_t^T,s}\phi(s,X_s)ds
\right]\nonumber \\
&\qquad =\mathbbm{1}(X_{\eta_t^T}\in\overline{D})w_0(\eta_t^T,X_{\eta_t^T}).\label{1wm is 1w0}
\end{align}
Hereafter, we employ the identity $\tau_t^{(m)}=\tau_{\tau_t^{(n)}}^{(m-n)}$ to indicate that the $m$-th jump after $t$ is the $(m-n)$-th jump based on the $n$-th jump after $t$.
Also, we have $\eta_t^T = \eta_{\tau_t^{(n)}}^T$ and $X_{\tau_t^{(n)}}\in\overline{D}$ on the event $\{\eta_t^T>\tau_t^{(n)}\}$, whereas $\eta_t^T \land \tau_t^{(m)} = \eta_t^T$ on the event $\{\eta_t^T\le \tau_t^{(n)}\}$.
The $n$-th jump time $\tau_t^{(n)}$ can be introduced into the first term in the expectation \eqref{wmdef2} as follows:
\begin{align}
&\mathbbm{1} (X_{\eta_t^T\land \tau_t^{(m)}}\in \overline{D}) \Theta_{t,\eta_t^T\land \tau_t^{(m)}} w_0(\eta_t^T\land \tau_t^{(m)}, X_{\eta_t^T\land \tau_t^{(m)}})\nonumber\\
&\quad =\mathbbm{1}(\eta_t^T>\tau_t^{(n)})\mathbbm{1} (X_{\eta_t^T\land \tau_t^{(m)}}\in \overline{D})
\Theta_{t,\tau_t^{(n)}} \Theta_{\tau_t^{(n)},\eta_t^T\land \tau_t^{(m)}} w_0(\eta_t^T\land \tau_t^{(m)}, X_{\eta_t^T\land \tau_t^{(m)}})
+\mathbbm{1}(\eta_t^T\le \tau_t^{(n)})\mathbbm{1}(X_{\eta_t^T}\in \overline{D}) \Theta_{t,\eta_t^T} w_0(\eta_t^T, X_{\eta_t^T})\nonumber\\
&\quad =\mathbbm{1}(\eta_t^T>\tau_t^{(n)})\mathbbm{1}(X_{\tau_t^{(n)}}\in\overline{D})
\Theta_{t,\tau_t^{(n)}} \mathbbm{1} \left(
X_{\eta_{\tau_t^{(n)}}^T\land \tau_{\tau_t^{(n)}}^{(m-n)}}\in \overline{D}\right)
\Theta_{\tau_t^{(n)},\eta_{\tau_t^{(n)}}^T\land \tau_{\tau_t^{(n)}}^{(m-n)}} w_0\left(\eta_{\tau_t^{(n)}}^T\land \tau_{\tau_t^{(n)}}^{(m-n)}, X_{\eta_{\tau_t^{(n)}}^T\land \tau_{\tau_t^{(n)}}^{(m-n)}}\right)\nonumber\\
&\qquad\qquad +\mathbbm{1}(\eta_t^T\le \tau_t^{(n)})\mathbbm{1}(X_{\eta_t^T}\in \overline{D}) 
\Theta_{t,\eta_t^T} w_{m-n}(\eta_t^T, X_{\eta_t^T}),\label{pwmrecur1}
\end{align}
where the equality for the second additive term holds true due to \eqref{1wm is 1w0}.
In a similar manner, the second term in the expectation \eqref{wmdef2} can be equipped with the $n$-th jump time $\tau_t^{(n)}$, as 
\begin{align}
    &\mathbbm{1}(X_{\eta_t^T\land \tau_t^{(m)}}\not\in \overline{D}) \Theta_{t,\eta^T_t} \Psi(\eta^T_t, X_{\eta^T_t},X_{\eta^T_t-})\nonumber \\
    &\quad =\mathbbm{1}(\eta_t^T>\tau_t^{(n)} )\mathbbm{1}(X_{\eta_t^T\land \tau_t^{(m)}}\notin \overline{D} ) \Theta_{t,\tau_t^{(n)}}\Theta_{\tau_t^{(n)},\eta_t} \Psi(\eta^T_t, X_{\eta^T_t},X_{\eta^T_t-}) 
    + \mathbbm{1}(\eta_t^T\le \tau_t^{(n)})\mathbbm{1} (X_{\eta_t^T}\notin \overline{D}) \Theta_{t,\eta^T_t} \Psi(\eta^T_t, X_{\eta^T_t},X_{\eta^T_t-})\nonumber\\
    &\quad =\mathbbm{1}(\eta_t^T>\tau_t^{(n)})\mathbbm{1}(X_{\tau_t^{(n)}}\in\overline{D})
\Theta_{t,\tau_t^{(n)}}
\mathbbm{1} \left(
X_{\eta_{\tau_t^{(n)}}^T\land \tau_{\tau_t^{(n)}}^{(m-n)}}\not\in \overline{D}\right)
\Theta_{\tau_t^{(n)},\eta_{\tau_t^{(n)}}^T} \Psi \left(\eta_{\tau_t^{(n)}}^T, X_{\eta_{\tau_t^{(n)}}^T}, X_{\eta_{\tau_t^{(n)}}^T-}\right)
+\mathbbm{1}(X_{\eta_t^T\land\tau_t^{(n)}}\not\in \overline{D}) \Theta_{t,\eta^T_t} \Psi(\eta^T_t, X_{\eta^T_t},X_{\eta^T_t-}),\label{pwmrecur2}
\end{align}
where the last equality holds by $\{\eta_t^T> \tau_t^{(n)}\}\cap \{X_{\tau_t^{(n)}}\not\in \overline{D}\}=\varnothing$, and thus
\[
 \mathbbm{1}(\eta_t^T\le \tau_t^{(n)})\mathbbm{1} (X_{\eta_t^T}\not\in \overline{D} )
 =\mathbbm{1}(\eta_t^T\le \tau_t^{(n)})\mathbbm{1} (X_{\eta_t^T}\not\in \overline{D} ) +
    \mathbbm{1}(\eta_t^T> \tau_t^{(n)})\mathbbm{1} (X_{\tau_t^{(n)}}\not\in \overline{D} )
    =\mathbbm{1} (X_{\eta_t^T\land\tau_t^{(n)}}\not\in \overline{D} ).
\]
Also, the third term in the expectation \eqref{wmdef2} can carry the $n$-th jump time $\tau_t^{(n)}$, as
\begin{align}
    \int_t^{\eta_t^T\land \tau_t^{(m)}} \Theta_{t,s} \phi(s, X_s) ds
    &=\mathbbm{1}(\eta_t^T>\tau_t^{(n)})\int_{\tau_t^{(n)}}^{\eta_t^T\land \tau_t^{(m)}} \Theta_{t,\tau_t^{(n)}}\Theta_{\tau_t^{(n)},s} \phi(s, X_s) ds
    +\mathbbm{1}(\eta_t^T>\tau_t^{(n)})\int_t^{\tau_t^{(n)}} \Theta_{t,s} \phi(s, X_s) ds
    +\mathbbm{1}(\eta_t^T\le \tau_t^{(n)})\int_t^{\eta_t^T} \Theta_{t,s} \phi(s, X_s) ds\nonumber\\
    &=\mathbbm{1}(\eta_t^T>\tau_t^{(n)}) \mathbbm{1}(X_{\tau_t^{(n)}}\in\overline{D})\Theta_{t,\tau_t^{(n)}}
    \int_{\tau_t^{(n)}}^{\eta_{\tau_t^{(n)}}^T\land \tau_{\tau_t^{(n)}}^{(m-n)}} \Theta_{\tau_t^{(n)},s} \phi(s, X_s) ds+
    \int_t^{\eta_t^T\land \tau_t^{(n)}}\Theta_{t,s} \phi(s, X_s) ds. \label{pwmrecur3}
\end{align}
By combining the first terms of \eqref{pwmrecur1}, \eqref{pwmrecur2} and \eqref{pwmrecur3} and taking the conditional expectation on the stopped $\sigma$-field $\mathcal{F}_{\tau_t^{(n)}}$, it holds by the strong Markov property of the underlying process \cite[Chapter III]{gihmanskorohod} that
\begin{align}
&\mathbb{E}_{\lambda}^{t,{\bf x}} \Bigg[
\mathbbm{1}(\eta_t^T>\tau_t^{(n)})\mathbbm{1}(X_{\tau_t^{(n)}}\in\overline{D})
\Theta_{t,\tau_t^{(n)}} \mathbbm{1} \left(
X_{\eta_{\tau_t^{(n)}}^T\land \tau_{\tau_t^{(n)}}^{(m-n)}}\in \overline{D}\right)
\Theta_{\tau_t^{(n)},\eta_{\tau_t^{(n)}}^T\land \tau_{\tau_t^{(n)}}^{(m-n)}} w_0\left(\eta_{\tau_t^{(n)}}^T\land \tau_{\tau_t^{(n)}}^{(m-n)}, X_{\eta_{\tau_t^{(n)}}^T\land \tau_{\tau_t^{(n)}}^{(m-n)}}\right)\nonumber\\
&\qquad\qquad +\mathbbm{1}(\eta_t^T>\tau_t^{(n)})\mathbbm{1}(X_{\tau_t^{(n)}}\in\overline{D})
\Theta_{t,\tau_t^{(n)}}
\mathbbm{1} \left(X_{\eta_{\tau_t^{(n)}}^T\land \tau_{\tau_t^{(n)}}^{(m-n)}}\not\in \overline{D}\right)
\Theta_{\tau_t^{(n)},\eta_{\tau_t^{(n)}}^T} \Psi \left(\eta_{\tau_t^{(n)}}^T, X_{\eta_{\tau_t^{(n)}}^T}, X_{\eta_{\tau_t^{(n)}}^T-}\right)\nonumber\\
&\qquad\qquad -\mathbbm{1}(\eta_t^T>\tau_t^{(n)})\mathbbm{1}(X_{\tau_t^{(n)}}\in\overline{D}) \Theta_{t,\tau_t^{(n)}}\int_{\tau_t^{(n)}}^{\eta_{\tau_t^{(n)}}^T\land \tau_{\tau_t^{(n)}}^{(m-n)}} \Theta_{\tau_t^{(n)},s} \phi(s, X_s) ds
\Bigg|\,\mathcal{F}_{\tau_t^{(n)}}\Bigg]=\mathbbm{1}(\eta_t^T>\tau_t^{(n)}) \mathbbm{1}(X_{\tau_t^{(n)}}\in\overline{D})\Theta_{t,\tau_t^{(n)}} w_{m-n}(\tau_t^{(n)},X_{\tau_t^{(n)}}),\label{pwmrecur4}
\end{align}
since the random variable $\mathbbm{1}(\eta_t^T>\tau_t^{(n)}) \mathbbm{1}(X_{\tau_t^{(n)}}\in\overline{D})\Theta_{t,\tau_t^{(n)}}$ is $\mathcal{F}_{\tau_t^{(n)}}$-measurable.
On the whole, by substituting \eqref{pwmrecur1}, \eqref{pwmrecur2}, \eqref{pwmrecur3} and \eqref{pwmrecur4} into the representation \eqref{wmdef2}, we obtain
\begin{multline*}
    w_m(t,{\bf x}) = \mathbb{E}_{\lambda}^{t,{\bf x}}\Bigg[
    \mathbbm{1}(\eta_t^T>\tau_t^{(n)}) \mathbbm{1}(X_{\tau_t^{(n)}}\in\overline{D})\Theta_{t,\tau_t^{(n)}} w_{m-n}(\tau_t^{(n)},X_{\tau_t^{(n)}})
    +\mathbbm{1}(\eta_t^T\le \tau_t^{(n)})\mathbbm{1}(
X_{\eta_t^T}\in \overline{D}) \Theta_{t,\eta_t^T} w_{m-n}(\eta_t^T, X_{\eta_t^T})\\
\qquad\qquad +\mathbbm{1}(X_{\eta_t^T\land\tau_t^{(n)}}\notin \overline{D}) \Theta_{t,\eta^T_t} \Psi(\eta^T_t, X_{\eta^T_t},X_{\eta^T_t-})
 -\int_t^{\eta_t^T\land \tau_t^{(n)}}\Theta_{t,s} \phi(s, X_s) ds\Bigg],
\end{multline*}
which reduces to the desired expression \eqref{wmrecur2} by combining the first two terms in the expectation.
\end{proof}

\begin{proof}[Proof of Proposition \ref{proposition pure jump ruin probability}]
Observe that $\{X_{\eta_t^T\land \tau_t^{(1)}}\in D\}=\{X_{T\land \tau_t^{(1)}}\in D\}$ as well as $\{X_{\eta_t^T\land \tau_t^{(1)}}\notin D\}=\{X_{T\land \tau_t^{(1)}}\notin D\}$ under $\mathbb{P}_{\lambda}^{t,{\bf x}}$, since the drift is not outward pointing at the boundary as well as there is no diffusion component in the neighborhood $\partial D_{\epsilon}$ of the boundary.
Hence, one can rewrite the representation \eqref{wmrecur2} as
\[
w_m(t,{\bf x}) =\mathbb{E}_\lambda^{t,{\bf x}} \left[\mathbbm{1}(X_{T\land \tau_t^{(1)}}\not\in\overline{D}) +\mathbbm{1}(X_{T\land \tau_t^{(1)}}\in\overline{D}) w_{m-1} (T\land \tau_t^{(1)} , X_{T\land \tau_t^{(1)} })\right].
\]
Moreover, it holds true that $\tau_t^{(1)}<\tau_t^{(2)}=\tau_{\tau_t^{(1)}}^{(1)}$ and $T\land \tau_{T\land \tau_t^{(1)}}^{(1)}=T\land \tau_t^{(2)}$, due to
\begin{align*}
T\land \tau_{T\land \tau_t^{(1)}}^{(1)} &=\mathbbm{1}(\tau_t^{(1)}\ge T)(T\land \tau_T^{(1)})
+\mathbbm{1}(\tau_t^{(1)}<T) (T\land \tau_{\tau_t^{(1)}}^{(1)})
=\mathbbm{1}(\tau_t^{(1)}\ge T)T
+\mathbbm{1}(\tau_t^{(1)}<T) (T\land \tau_t^{(2)})\\
&=\mathbbm{1}(\{\tau_t^{(1)}\ge T\}\cap \{\tau_t^{(2)}\geq T\})T
+\mathbbm{1}(\{\tau_t^{(1)}<T\}\cap\{\tau_t^{(2)}\geq T\})T 
+\mathbbm{1}(\{\tau_t^{(1)}<T\}\cap \{\tau_t^{(2)}< T\})\tau_t^{(2)}\\
&=\mathbbm{1}(\tau_t^{(2)}\ge T)T
+\mathbbm{1}(\tau_t^{(2)}< T)\tau_t^{(2)}= T\land \tau_t^{(2)}.
\end{align*}
Hence, we obtain
\begin{align*}
w_m(t,{\bf x}) {}
&=\mathbb{E}_\lambda^{t,{\bf x}} \left[ \mathbbm{1}(X_{\eta_t^T\land \tau_t^{(1)}}\not\in\overline{D})
+\mathbbm{1}(X_{\eta_t^T\land \tau_t^{(1)}}\in\overline{D}) w_{m-1} (\eta_t^T\land \tau_t^{(1)} , X_{\eta_t^T\land \tau_t^{(1)}}) 
 \right]\\
&=\mathbb{E}_\lambda^{t,{\bf x}} \left[ \mathbbm{1}(X_{T\land \tau_t^{(1)}}\not\in\overline{D})
+\mathbbm{1}(X_{T\land \tau_t^{(1)}}\in\overline{D}) 
\left(
\mathbbm{1}(X_{T\land \tau_t^{(2)}}\not\in\overline{D})+
\mathbbm{1}(X_{T\land \tau_t^{(2)}}\in\overline{D}) 
w_{m-2} (T\land \tau_t^{(2)} , X_{T\land \tau_t^{(2)} }) \right)
 \right],
\end{align*}
where we have applied the tower property and the strong Markov property.
This yields the desired identity by induction with the aid of $w_0(t,{\bf x})=\mathbb{E}_0^{t,{\bf x}}[\mathbbm{1}(\eta_t\le T)]\equiv 0$ due to the degeneracy of the jump component. 
The dominance is evident from the increasing number of indicator functions in $m$.
\end{proof}

\begin{proof}[Proof of Theorem \ref{theorem vm wm}]
{\bf (a)} Since overshooting can only be caused by a jump, the event $\{X_{\eta_t^T\land\tau_t^{(m+1)}}\not\in\overline{D}\}$ means that an exit occurs before the terminal time and, moreover, at latest at the $(m+1)$-st jump timing.
Hence, it holds that
\begin{align*}
\mathbbm{1}(X_{\eta_t^T\land\tau_t^{(m+1)}}\not\in\overline{D})
=\mathbbm{1}(\eta_t \leq \tau_t^{(m+1)}\land T)\mathbbm{1}(X_{\eta_t}\not\in\overline{D})
&=\mathbbm{1}(\eta_t \leq \tau_t^{(m)}\land T)\mathbbm{1}(X_{\eta_t}\not\in\overline{D})
+\mathbbm{1}(\eta_t^T=\tau_t^{(m+1)})\mathbbm{1}(X_{\tau_t^{(m+1)}}\not\in\overline{D})\\
&=\mathbbm{1}(X_{\eta_t^T\land\tau_t^{(m)}}\not\in\overline{D})
+\mathbbm{1}(\tau_t^{(m)}<\eta_t^T)\mathbbm{1}(X_{\eta_t^T\land \tau_t^{(m+1)}}\not\in\overline{D}).
\end{align*}
With this identity, the expression \eqref{def of vm} (with $m+1$) can be rewritten as follows:
\begin{align}
v_{m+1}(t,{\bf x})
&=\mathbb{E}_\lambda^{t,{\bf x}}\left[
\mathbbm{1}(\tau_t^{(m+1)} \geq \eta_t^T, \tau_t^{(m)}\geq \eta_t^T)
\mathbbm{1}(X_{\eta_t^T}\in\overline{D})\Theta_{t,\eta_t^T}w_0(\eta_t^T,X_{\eta_t^T})
+\mathbbm{1}(X_{\eta_t^T\land\tau_t^{(m)}}\not\in\overline{D})\Theta_{t,\eta_t}\Psi(\eta_t,X_{\eta_t},X_{\eta_t-})
-\int_t^{\eta_t^T\land\tau_t^{(m)}}\Theta_{t,s}\phi(s,X_s)ds\right]\nonumber\\
& \qquad +\mathbb{E}_\lambda^{t,{\bf x}}\Bigg[
\mathbbm{1}(\tau_t^{(m+1)} \geq \eta_t^T, \tau_t^{(m)} < \eta_t^T)
\mathbbm{1}(X_{\eta_t^T}\in\overline{D})\Theta_{t,\eta_t^T}w_0(\eta_t^T,X_{\eta_t^T})\nonumber\\
& \qquad \qquad\qquad \qquad \qquad\qquad
 +\mathbbm{1}(\tau_t^{(m)} < \eta_t^T)
\mathbbm{1}(X_{\eta_t^T\land \tau_t^{(m+1)}}\not\in\overline{D})\Theta_{t,\eta_t}
\Psi(\eta_t,X_{\eta_t},X_{\eta_t-})-\int_{\eta_t^T\land\tau_t^{(m)}}^{\eta_t^T\land\tau_t^{(m+1)}}\Theta_{t,s}\phi(s,X_s)ds\Bigg]\nonumber\\
&=v_m(t,{\bf x})+\mathbb{E}_\lambda^{t,{\bf x}}\Bigg[\mathbbm{1}(\tau_t^{(m)}< \eta_t^T)
\Theta_{t,\tau_t^{(m)}}\Bigg(\mathbbm{1}(X_{\tau_t^{(m+1)}}\in\overline{D})
\mathbbm{1}(\tau_t^{(m+1)}\geq\eta_t^T)
\Theta_{\tau_t^{(m)},\eta_t^T}w_0(\eta_t^T,X_{\eta_t^T})\nonumber\\
&\qquad \qquad\qquad \qquad\qquad \qquad
+\mathbbm{1}(X_{\eta_t^T\land\tau_t^{(m+1)}}\not\in\overline{D})\Theta_{\tau_t^{(m)},\eta_t}
\Psi(\eta_t,X_{\eta_t},X_{\eta_t-})-\int_{\tau_t^{(m)}}^{\eta_t^T\land\tau_t^{(m+1)}}\Theta_{\tau_t^{(m)},s}\phi(s,X_s)ds\Bigg)\Bigg],\nonumber\\
&=v_m(t,{\bf x})+\mathbb{E}_\lambda^{t,{\bf x}}\left[\mathbbm{1}(\tau_t^{(m)}< \eta_t^T)\Theta_{t,\tau_t^{(m)}}v_1(\tau_t^{(m)},X_{\tau_t^{(m)}})\right].\label{vm to vm+1}
\end{align}
By applying the identity \eqref{vm to vm+1} recursively, we further obtain that
\begin{align*}
    v_{m+1}(t,{\bf x})&= v_m(t,{\bf x})+\mathbb{E}_\lambda^{t,{\bf x}}\left[\mathbbm{1}(\tau_t^{(m)}< \eta_t^T)\Theta_{t,\tau_t^{(m)}}v_1(\tau_t^{(m)},X_{\tau_t^{(m)}})\right]\\
    &=v_{m-1}(t,{\bf x})+\mathbb{E}_\lambda^{t,{\bf x}}\left[\mathbbm{1}(\tau_t^{(m-1)}< \eta_t^T)\Theta_{t,\tau_t^{(m-1)}}\left(v_1(\tau_t^{(m-1)},X_{\tau_t^{(m-1)}})
    +\mathbbm{1}(\tau_t^{(m)}< \eta_t^T)\Theta_{\tau_t^{(m-1)},\tau_t^{(m)}}v_1(\tau_t^{(m)},X_{\tau_t^{(m)}})\right)\right]\\
    &=v_{m-1}(t,{\bf x})+\mathbb{E}_\lambda^{t,{\bf x}}\left[\mathbbm{1}(\tau_t^{(m-1)}< \eta_t^T)\Theta_{t,\tau_t^{(m-1)}}v_2(\tau_t^{(m-1)},X_{\tau_t^{(m-1)}})\right],
\end{align*}
which yields the first identity \eqref{iteration vm} for general $n\in \{0,1,\cdots,m-1\}$ by induction.

For the second identity \eqref{wm vm}, it is straightforward in light of \eqref{wmdef2} and \eqref{def of vm} to obtain that
\begin{equation}\label{vm to wm+1}
 w_m(t,{\bf x})=v_m(t,{\bf x})+\mathbb{E}_\lambda^{t,{\bf x}}\left[\mathbbm{1}(\tau_t^{(m)}< \eta_t^T)\Theta_{t,\tau_t^{(m)}}w_0(\tau_t^{(m)},X_{\tau_t^{(m)}})\right].
\end{equation}
In a similar manner to the derivation of the iteration \eqref{iteration vm}, it holds by the identity \eqref{vm to wm+1} that
\begin{align*}
    w_{m}(t,{\bf x})
    &=v_{m-1}(t,{\bf x})+\mathbb{E}_\lambda^{t,{\bf x}}\left[
    \mathbbm{1}(\tau_t^{(m-1)}< \eta_t^T)\Theta_{t,\tau_t^{(m-1)}}\left(v_1(\tau_t^{(m-1)},X_{\tau_t^{(m-1)}})+
    \mathbbm{1}(\tau_t^{(m)}< \eta_t^T)\Theta_{\tau_t^{(m-1)},\tau_t^{(m)}}w_0(\tau_t^{(m)},X_{\tau_t^{(m)}})\right)\right]\\
    &=v_{m-1}(t,{\bf x})+\mathbb{E}_\lambda^{t,{\bf x}}\left[\mathbbm{1}(\tau_t^{(m-1)}< \eta_t^T)\Theta_{t,\tau_t^{(m-1)}}w_1(\tau_t^{(m-1)},X_{\tau_t^{(m-1)}})\right],
\end{align*}
which yields the first identity \eqref{wm vm} for general $n\in \{0,1,\cdots,m-1\}$ by induction.

To prove the pointwise convergence, observe first that $\tau_t^{(m)}>\eta_t^T$, $\mathbb{P}_\lambda^{t,{\bf x}}$-$a.s.$ eventually as $m\to +\infty$, due to a finite jump intensity.
Hence, it holds as $m\to +\infty$ that
\[
v_m(t,{\bf x})\to\mathbb{E}_\lambda^{t,{\bf x}}\left[
\mathbbm{1}(X_{\eta_t^T}\in\overline{D})\Theta_{t,\eta_t^T}w_0(\eta_t^T,X_{\eta_t^T})
+\mathbbm{1}(X_{\eta_t^T}\not\in\overline{D})\Theta_{t,\eta_t^T}\Psi(\eta_t^T,X_{\eta_t^T},X_{\eta_t^T-})
-\int_t^{\eta_t^T}\Theta_{t,s}\phi(s,X_s)ds
\right],
\]
which yields \eqref{udef}, due to $\mathbbm{1}(\eta_t>T)w_0(T,X_T)=\mathbbm{1}(\eta_t>T)g(X_T)$ and $\mathbbm{1}(X_{\eta_t}\in\overline{ D})w_0(\eta_t,X_{\eta_t})=\mathbbm{1}(X_{\eta_t}\in\partial D)\Psi(\eta_t,X_{\eta_t},X_{\eta_t-})$.
The passage to the limit can be justified by the bounded convergence theorem in a similar manner to Theorem \ref{theorem convergence}.

\noindent {\bf (b)}
If $g\ge 0$, $\Psi\ge 0$ and $\phi\le 0$ (respectively, if $g\le 0$, $\Psi\le 0$ and $\phi\ge 0$), then $w_0$ and then $v_1$ are non-negative (respectively, non-positive) due to \eqref{w0def1} and \eqref{def of vm} over the domain.
Therefore, the desired monotonicity holds true in light of the identity \eqref{iteration vm}.
\end{proof}

\begin{proof}[Proof of Theorem \ref{theorem back to E0}]
We first couple the underlying process \eqref{Xprocess} with $\{Y_t:\,t\ge 0\}$ defined as a unique solution of the stochastic differential equation
\begin{equation}\label{sde of Y}
dY_t=b(t,Y_t)dt+\sigma(t,Y_t)dW_t,
\end{equation}
and let the probability measure $\mathbb{P}_{\lambda}^{t,{\bf x}}$ represent the condition $X_t=Y_t={\bf x}$.
We define discount functions on $\{Y_t:\,t\ge 0\}$ by $\Theta^Y_{t_1,t_2}:= \exp [ - \int_{t_1}^{t_2} r(s,Y_s) ds]$ and $\Lambda^Y_{t_1,t_2}:=  \exp[ - \int_{t_1}^{t_2} \lambda(s,Y_s) ds]$ for $t\le t_1\le t_2$, as well as the first exit time of the coupled process $\{Y_t:\,t\ge 0\}$ and its truncation, respectively, $\xi_t:=\inf\{s\geq t: Y_s\not\in\overline{D}\}$ and $\xi_t^T:=\xi_t\land T$, in a similar manner to \eqref{etadef} for the underlying process $\{X_t:\,t\ge 0\}.$
It then holds $\mathbb{P}_\lambda^{t,{\bf x}}$-$a.s.$ that $X_s=Y_s$ for all $s\in[t,\tau_t^{(1)})$, as well as $\Theta^Y_{t_1,t_2} = \Theta_{t_1,t_2}$ and $\Lambda^Y_{t_1,t_2} =\Lambda_{t_1,t_2}$ for all $t\le t_1\le t_2$ with $t_2\le \tau_t^{(1)}$.
Moreover, we have $\{\eta_t^T<\tau_t^{(1)}\}=\{\xi_t^T<\tau_t^{(1)}\}$, $\mathbb{P}_\lambda^{t,{\bf x}}$-$a.s.$, since the left-hand side indicates that the first exit occurs by the drift or diffusion component strictly before the first jump of $\{X_t:\,t\ge 0\}$, and vice versa.
Hence, we get $\eta_t^T\land\tau_t^{(1)}=\xi_t^T\land\tau_t^{(1)}$, $\mathbb{P}_\lambda^{t,{\bf x}}$-$a.s.$, with which the representation \eqref{wmrecur2} with $n=1$ can be rewritten as
\begin{align}
w_m(t,{\bf x})
& =\mathbb{E}_\lambda^{t,{\bf x}} \Bigg[ \mathbbm{1}(X_{\xi_t^T\land \tau_t^{(1)}}\in\overline{D})\Theta_{t, \xi_t^T \land \tau_t^{(1)} } w_{m-1} (\xi_t^T\land \tau_t^{(1)} , X_{\xi_t^T\land \tau_t^{(1)} })\nonumber \\
&\qquad \qquad +\mathbbm{1}(X_{\xi_t^T\land \tau_t^{(1)}}\not\in\overline{D})\Theta_{t,\tau_t^{(1)}} \Psi(\tau_t^{(1)}, X_{\tau_t^{(1)}},X_{\tau_t^{(1)}-}) 
- \int_t^{\xi_t^T \land \tau_t^{(1)}} \Theta_{t,s} \phi(s,X_s) ds \Bigg]\nonumber\\
&=\mathbb{E}_\lambda^{t,{\bf x}} \Bigg[ \mathbbm{1}(\tau_t^{(1)}>\xi_t^ T)\Theta^Y_{t, \xi_t^ T  } w_{m-1} (\xi_t^ T , Y_{\xi_t^T })
+\mathbbm{1}(\tau_t^{(1)}\leq\xi_t^ T)\mathbbm{1}(Y_{ \tau_t^{(1)}}+X_{\tau_t^{(1)}}-X_{\tau_t^{(1)}-}\in\overline{D})\Theta^Y_{t, \tau_t^{(1)} } w_{m-1} ( \tau_t^{(1)} , Y_{ \tau_t^{(1)} }+X_{\tau_t^{(1)}}-X_{\tau_t^{(1)}-})\nonumber\\
&\qquad\qquad
+\mathbbm{1}(\tau_t^{(1)}\leq\xi_t^T)\mathbbm{1}(Y_{ \tau_t^{(1)}}+X_{\tau_t^{(1)}}-X_{\tau_t^{(1)}-}\not\in\overline{D})\Theta^Y_{t,\tau_t^{(1)}} \Psi(\tau_t^{(1)}, Y_{\tau_t^{(1)}}+X_{\tau_t^{(1)}}-X_{\tau_t^{(1)}-},Y_{\tau_t^{(1)}}) 
- \int_t^{\xi_t^ T} \mathbbm{1}(\tau_t^{(1)}>s)\Theta^Y_{t,s} \phi(s,Y_s) ds \Bigg],\label{pwme011}
\end{align}
where we have applied the identity $X_{\tau_t^{(1)}-}=Y_{\tau_t^{(1)}-}=Y_{\tau_t^{(1)}}$ for the second equality.

Now, fix $(t,{\bf x})\in [0,T]\times \overline{D}$ and let $\mathcal{G}$ denote the (completed) $\sigma$-field generated by the underlying Brownian motion $\{W_s:\,s\in [t,T]\}$.
Note that the first exit time $\xi_t^T$ is $\mathcal{G}$-measurable and, on the $\sigma$-field $\mathcal{G}$, the random variable $\tau_t^{(1)}$ can be treated as the first jump time of an inhomogeneous Poisson process with rate function $\lambda(s,Y_s)$ for $s\in [t,T]$.
For the first term of \eqref{pwme011}, it holds that
\begin{align*}
\mathbb{E}_\lambda^{t,{\bf x}} \left[ \mathbbm{1}(\tau_t^{(1)}>\xi_t^T)\Theta^Y_{t, \xi_t^T  } w_{m-1} (\xi_t^T , Y_{\xi_t^T }) \right]
&=\mathbb{E}_\lambda^{t,{\bf x}} \left[ \mathbb{E}_\lambda^{t,{\bf x}} \left[\mathbbm{1}(\tau_t^{(1)}>\xi_t^T)\Big|\, \mathcal{G} \right]\Theta^Y_{t, \xi_t^T } w_{m-1} (\xi_t^T , Y_{\xi_t^T }) \right]\\
&=\mathbb{E}_\lambda^{t,{\bf x}} \left[\Lambda_{t,\xi_t^T}^Y\Theta^Y_{t, \xi_t^T  } w_{m-1} (\xi_t^T , Y_{\xi_t^T}) \right]\\
&=\mathbb{E}_\lambda^{t,{\bf x}} \left[ \mathbbm{1}(\xi_t>T) \Lambda_{t,T}^Y\Theta^Y_{t, T  } g(Y_{T})+
\mathbbm{1}(\xi_t\leq T) \Lambda_{t,\xi_t}^Y\Theta^Y_{t, \xi_t  } \Psi(\xi_t,Y_{\xi_t},Y_{\xi})
\right],
\end{align*}
where the last equality holds because \eqref{sde of Y} contains no jump component, and moreover,  $\mathbbm{1}(\xi_t>T)w_{m-1}(T,Y_T)=\mathbbm{1}(\xi_t>T)w_0(T,Y_T)=\mathbbm{1}(\xi_t>T)g(Y_T)$, and 
$\mathbbm{1}(\xi_t\leq T)w_{m-1}(\xi_t,Y_{\xi_t})=\mathbbm{1}(\xi_t\leq T)w_0(\xi_t,Y_{\xi_t})=\mathbbm{1}(\xi_t\leq T)\Psi(\xi_t,Y_{\xi_t},Y_{\xi_t})$, due to \eqref{w0def1} and \eqref{wmdef2}.
For the second term of \eqref{pwme011}, observe that 
\begin{align*}
    &\mathbb{E}_\lambda^{t,{\bf x}} \left[\mathbbm{1}(\tau_t^{(1)}\leq\xi_t^T)\mathbbm{1}(Y_{ \tau_t^{(1)}}+X_{\tau_t^{(1)}}-X_{\tau_t^{(1)}-}\in\overline{D})\Theta^Y_{t, \tau_t^{(1)} } w_{m-1} ( \tau_t^{(1)} , Y_{\tau_t^{(1)}}+X_{\tau_t^{(1)}}-X_{\tau_t^{(1)}-})\right]\\
    &\qquad\qquad
    =\mathbb{E}_\lambda^{t,{\bf x}} \left[\mathbbm{1}(\tau_t^{(1)}\leq\xi_t^T)\Theta^Y_{t, \tau_t^{(1)} }\mathbb{E}_\lambda^{\tau_t^{(1)},Y_{\tau_t^{(1)}}} \left[\mathbbm{1}(Y_{\tau_t^{(1)}}+X_{\tau_t^{(1)}}-X_{\tau_t^{(1)}-}\in\overline{D}) w_{m-1} ( \tau_t^{(1)} , Y_{ \tau_t^{(1)} }+X_{\tau_t^{(1)}}-X_{\tau_t^{(1)}-}) \right]\right]\\
    &\qquad\qquad
    =\mathbb{E}_\lambda^{t,{\bf x}} \left[\mathbbm{1}(\tau_t^{(1)}\leq\xi_t^T)\Theta^Y_{t, \tau_t^{(1)} }
    \int_{\{{\bf z}\in\mathbb{R}_0^d:Y_{\tau_t^{(1)}}+{\bf z}\in\overline{D}\}}  w_{m-1} ( \tau_t^{(1)} , Y_{\tau_t^{(1)} }+{\bf z}) \nu(d{\bf z}; Y_{\tau_t^{(1)}})
    \right]\\
     &\qquad\qquad
     =\mathbb{E}_\lambda^{t,{\bf x}} \left[
     \int_t^{\xi_t^T}
     \left[\Theta^Y_{t, s }
    \int_{\{{\bf z}\in\mathbb{R}_0^d:Y_{s}+{\bf z}\in\overline{D}\}}  w_{m-1} ( s , Y_{s }+{\bf z})\nu(d{\bf z}; Y_{s})\right] \lambda(s,Y_s)\Lambda^Y_{t,s}ds
     \right]=\mathbb{E}_\lambda^{t,{\bf x}} \left[\int_t^{\xi_t^T} \Lambda^Y_{t,s}\Theta^Y_{t,s}G_{m-1}(s,Y_s)ds\right],
\end{align*}
where the first and second equalities hold by the tower property on the stopped $\sigma$-field $\mathcal{F}_{\tau_t^{(1)}}$ and the strong Markov property, the third by the tower property on the $\sigma$-field $\mathcal{G}$.
In a similar manner, for the third term of \eqref{pwme011}, we obtain 
\[
\mathbb{E}_\lambda^{t,{\bf x}} \left[
\mathbbm{1}(\tau_t^{(1)}\leq\xi_t^T)\mathbbm{1}(Y_{ \tau_t^{(1)}}+X_{ \tau_t^{(1)}}-X_{ \tau_t^{(1)}-}\not\in\overline{D})\Theta^Y_{t,\tau_t^{(1)}} \Psi(\tau_t^{(1)}, Y_{\tau_t^{(1)}}+X_{ \tau_t^{(1)}}-X_{ \tau_t^{(1)}-},Y_{\tau_t^{(1)}}) 
\right] =\mathbb{E}_\lambda^{t,{\bf x}} \left[\int_t^{\xi_t^T} \Lambda^Y_{t,s}\Theta^Y_{t,s}H(s,Y_s)ds\right].
\]
For the last term of \eqref{pwme011}, by the tower property on the $\sigma$-field $\mathcal{G}$ and interchanging the conditional expectation and the integral by the Fubini theorem, we obtain
\begin{align*}
    \mathbb{E}_\lambda^{t,{\bf x}} \left[ \int_t^{\xi_t^T} \mathbbm{1}(\tau_t^{(1)}>s)\Theta^Y_{t,s} \phi(s,Y_s) ds\right]
    &=\mathbb{E}_\lambda^{t,{\bf x}} \left[ \int_t^{\xi_t^T} \mathbb{E}_\lambda^{t,{\bf x}}\left[\mathbbm{1}(\tau_t^{(1)}>s)\Theta^Y_{t,s} \phi(s,Y_s)\Big|\, \mathcal{G}\right] ds\right]\\
    &=\mathbb{E}_\lambda^{t,{\bf x}} \left[ \int_t^{\xi_t^T} \mathbb{E}_\lambda^{t,{\bf x}}\left[\mathbbm{1}(\tau_t^{(1)}>s)\Big|\, \mathcal{G}\right]\Theta^Y_{t,s} \phi(s,Y_s) ds\right]
    =\mathbb{E}_\lambda^{t,{\bf x}} \left[ \int_t^{\xi_t^T} \Lambda_{t,s}^Y\Theta^Y_{t,s} \phi(s,Y_s) ds\right],
\end{align*}
since the random variable $\Theta^Y_{t,s} \phi(s,Y_s)$ is $\mathcal{G}$-measurable for all $s\in [t,T]$.
On the whole, the representation \eqref{pwme011} can be rewritten as
\[
w_m(t,{\bf x})=\mathbb{E}_0^{t,{\bf x}} \left[
\mathbbm{1}(\xi_t>T) \Lambda_{t,T}^Y\Theta^Y_{t, T  } g(Y_{T})+
\mathbbm{1}(\xi_t\leq T) \Lambda_{t,\xi_t}^Y\Theta^Y_{t, \xi_t  } \Psi(\xi_t,Y_{\xi_t},Y_{\xi_t})  - \int_t^{\xi_t^T} \Lambda^Y_{t,s}\Theta^Y_{t,s} \left(\phi(s,Y_s)-G_{m-1}(s,Y_s)-H(s,Y_s)\right) ds\right],
\]
where the expectation $\mathbb{E}_\lambda^{t,{\bf x}}$ has been replaced by $\mathbb{E}_0^{t,{\bf x}}$ since the integrand is written on the coupled process $\{Y_s:\,s\in[t,T]\}$ with no jumps involved.
Finally, we get the desired result \eqref{wme0} since, under the probability measure $\mathbb{P}_0^{t,{\bf x}}$, the two processes \eqref{Xprocess} and \eqref{sde of Y} on $[t,T]$ are identical in law.
\end{proof}

Before continuing the proofs, we build a lemma, which will be useful hereafter from place to place. 

\begin{lemma}\label{intermediate lemma}
Let $f:[0,\infty)\times\mathbb{R}^d\to\mathbb{R}$ be bounded and such that $\partial_1 f$ and $\mathcal{L}_{\cdot} f$ exist almost everywhere on $[0,+\infty)\times \overline{D}$.
It holds that for $(t,{\bf x})\in [0,T]\times \overline{D}$ and an $(\mathcal{F}_s)_{s\in [t,T]}$-stopping time $\kappa_t$ such that $\kappa_t\ge t$,
\begin{multline}\label{lemma51a}
\mathbb{E}_\lambda^{t,{\bf x}}\left[\mathbbm{1}(X_{\kappa_t\land \eta_t}\in\overline{D})\Theta_{t,\kappa_t\land\eta_t}f(\kappa_t\land \eta_t,X_{\kappa_t\land \eta_t})\right] - \mathbbm{1}({\bf x}\in\overline{D})f(t,{\bf x})\\
= \mathbb{E}_\lambda^{t,x}\left[\int_{t}^{\kappa_t\land \eta_t}\Theta_{t,s}\left(\partial_1 f(s,X_s)+\mathcal{L}_s f(s,X_s) - \lambda(s,X_s)f(s,X_s)+ \int_{\{{\bf z}\in\mathbb{R}_0^d:X_s+{\bf z}\in\overline{D}\}}  f(s, X_s+{\bf z}) \lambda(s,X_s)\nu(d{\bf z};X_s)\right)ds\right],
\end{multline}
and 
\begin{equation}\label{lemma51b}
\mathbb{E}_\lambda^{t,{\bf x}}\left[\mathbbm{1}(X_{\kappa_t\land \eta_t}\not\in\overline{D})\Theta_{t,\eta_t}f(\eta_t,X_{\eta_t})\right]
=\mathbb{E}_\lambda^{t,{\bf x}}\left[\int_{t}^{\kappa_t\land\eta_t}\Theta_{t,s} \int_{\{{\bf z}\in\mathbb{R}_0^d:X_s+{\bf z}\not\in\overline{D}\}}  f(s, X_s+{\bf z}) \lambda(s,X_s)\nu(d{\bf z};X_s)ds\right].
\end{equation}
\end{lemma}

\begin{proof}[Proof of Lemma \ref{intermediate lemma}]
Let $(t,{\bf x})\in [0,T]\times \overline{D}$.
By applying Ito formula to $\Theta_{t,\cdot}f(\cdot,X_{\cdot})$ over the interval $(t,\kappa_t\land\eta_t]$, it holds that 
\begin{multline*}
\mathbb{E}_\lambda^{t,{\bf x}}\left[\Theta_{t,\kappa_t\land\eta_t}f(\kappa_t\land \eta_t,X_{\kappa_t\land \eta_t})\right] - f(t,{\bf x})\\
 = \mathbb{E}_\lambda^{t,{\bf x}}\left[\int_{t}^{\kappa_t\land \eta_t}\Theta_{t,s-}\left(\partial_1 f(s,X_{s-})+\mathcal{L}_s f(s,X_{s-}) - \lambda(s,X_{s-})f(s,X_{s-})+ \int_{\mathbb{R}^d_0}  f(s, X_{s-}+{\bf z}) \lambda(s,X_{s-})\nu(d{\bf z};X_{s-})\right)ds\right],
\end{multline*}
and 
\begin{multline*}
\mathbb{E}_\lambda^{t,{\bf x}}\left[\mathbbm{1}(X_{\kappa_t\land \eta_t}\in\overline{D})\Theta_{t,\kappa_t\land\eta_t}f(\kappa_t\land \eta_t,X_{\kappa_t\land \eta_t})\right] - \mathbbm{1}({\bf x}\in\overline{D})f(t,{\bf x})\\
 = \mathbb{E}_\lambda^{t,{\bf x}}\left[\int_{t}^{\kappa_t\land\eta_t}\Theta_{t,s-}\left(\partial_1 f(s,X_{s-})+\mathcal{L}_s f(s,X_{s-}) - \lambda(s,X_{s-})f(s,X_{s-})+ \int_{\{{\bf z}\in\mathbb{R}_0^d:X_{s-}+{\bf z}\in\overline{D}\}}  f(s, X_{s-}+{\bf z})\lambda(s,X_{s-}) \nu(d{\bf z};X_{s-})\right)ds\right],
\end{multline*}
where we have kept the left limits in time to be precise in the location of the trajectory.
Note that the term $\lambda(s,X_{s-})f(s,X_{s-})$ appears above since the L\'evy measure $\nu(d\cdot;{\bf x})$ here is a probability measure for all ${\bf x}\in \overline{D}$.
By the above two identities, we get
\begin{align*}
\mathbb{E}_\lambda^{t,{\bf x}}\left[\mathbbm{1}(X_{\kappa_t\land \eta_t}\not\in\overline{D})\Theta_{t,\kappa_t\land\eta_t}f(\kappa_t\land \eta_t,X_{\kappa_t\land \eta_t})\right]
&=\left(\mathbb{E}_\lambda^{t,{\bf x}}\left[\Theta_{t,\kappa_t\land\eta_t}f(\kappa_t\land \eta_t,X_{\kappa_t\land \eta_t})\right] - f(t,{\bf x})\right)\\
&\qquad -\left(\mathbb{E}_\lambda^{t,{\bf x}}\left[\mathbbm{1}(X_{\kappa_t\land \eta_t}\in\overline{D})\Theta_{t,\kappa_t\land\eta_t}f(\kappa_t\land \eta_t,X_{\kappa_t\land \eta_t})\right] - \mathbbm{1}({\bf x}\in\overline{D})f(t,{\bf x})\right)\\
&=\mathbb{E}_\lambda^{t,{\bf x}}\left[\int_{t}^{\kappa_t\land\eta_t} \Theta_{t,s-}\int_{\{{\bf z}\in\mathbb{R}_0^d:X_{s-}+{\bf z}\not\in\overline{D}\}}  f(s, X_{s-}+{\bf z}) \lambda(s,X_{s-})\nu(d{\bf z};X_{s-})ds\right],
\end{align*}
which concludes due to $\eta_t\le \kappa_t$ on the event $\{X_{\kappa_t\land \eta_t}\not\in\overline{D}\}$ and the integrand is absolutely continuous with respect to the Lebesgue measure $ds$.
\end{proof}

\begin{proof}[Proof of Theorem \ref{theorem when smooth}]
\noindent {\bf (a)} Let $(t,{\bf x})\in[0,T) \times D$ and 
$\delta\in(0,T-t)$.
We rewrite the function $w_m$ in the representation \eqref{wmrecur2} in a local form on the time interval $(t,\eta_t^{t+\delta}]$ as $\delta\to 0+$ by considering two cases separately as to whether or not there is at least one jump on the interval.

First, consider the case where there is no jump on the interval $(t,(t+\delta)\land \eta_t^T]$, that is, the first jump time $\tau_t^{(1)}$ after $t$ is strictly later than the end time $\eta_t^{t+\delta}(=(t+\delta)\land \eta_t^T)$.
Hence, we get $\tau_t^{(1)}=\tau_{\eta_t^{t+\delta}}^{(1)}$.
Moreover, observe that $\eta^T_{\eta_t^{t+\delta}}=\mathbbm{1}(\eta^T_t \leq t+\delta)\eta^T_{\eta^T_t}+\mathbbm{1}(\eta^T_t >t+\delta)\eta^T_{t+\delta}=\eta^T_t\land\eta^T_{t+\delta}=\eta^T_t,$ due to \eqref{etaTetaT is etaT}.
By applying those identities, it holds that
\begin{align}
&\mathbb{E}_\lambda^{t,{\bf x}} \left[\mathbbm{1}(\tau_t^{(1)} >\eta_t^{t+\delta}) \left(\mathbbm{1}(X_{\eta_t^T\land \tau_t^{(1)}}\in\overline{D})\Theta_{t, \eta_t^T \land \tau_t^{(1)} } w_{m-1} (\eta_t^T\land \tau_t^{(1)} , X_{\eta_t^T\land \tau_t^{(1)} }) 
+\mathbbm{1}(X_{\eta_t^T\land \tau_t^{(1)}}\not\in\overline{D})\Theta_{t,\eta^T_t} \Psi(\eta^T_t, X_{\eta^T_t},X_{\eta^T_t-}) 
- \int_t^{\eta_t^T \land \tau_t^{(1)}} \Theta_{t,s} \phi(s,X_s) ds\right) \right]\nonumber\\
&\qquad =\mathbb{E}_\lambda^{t,{\bf x}} \Bigg[\mathbbm{1}(\tau_t^{(1)} > \eta_t^{t+\delta}) \Theta_{t,\eta_t^{t+\delta}} 
\mathbb{E}_\lambda^{\eta_t^{t+\delta},X_{\eta_t^{t+\delta}}}
\Bigg[\mathbbm{1}\left(X_{\eta_{\eta_t^{t+\delta}}^T\land \tau_{\eta_t^{t+\delta}}^{(1)}}\in\overline{D}\right)\Theta_{\eta_t^{t+\delta}, \eta_{\eta_t^{t+\delta}}^T \land \tau_{\eta_t^{t+\delta}}^{(1)} } w_{m-1} \left(\eta_{\eta_t^{t+\delta}}^T\land \tau_{\eta_t^{t+\delta}}^{(1)} , X_{\eta_{\eta_{t}^{t+\delta}}^T\land \tau_{\eta_t^{t+\delta}}^{(1)} }\right)\nonumber\\
&\qquad\quad +
\mathbbm{1}\left(X_{\eta_{\eta_t^{t+\delta}}^T\land \tau_{\eta_t^{t+\delta}}^{(1)}}\not\in\overline{D}\right)\Theta_{{\eta_t^{t+\delta}},\eta^T_{\eta_t^{t+\delta}}} \Psi\left(\eta^T_{\eta_t^{t+\delta}}, X_{\eta^T_{\eta_t^{t+\delta}}},X_{\eta^T_ {\eta_t^{t+\delta}}-}\right) 
- \int_{\eta_t^{t+\delta}}^{\eta_{\eta_t^{t+\delta}}^T \land \tau_{\eta_t^{t+\delta}}^{(1)}} \Theta_{{\eta_t^{t+\delta}},s} \phi(s,X_s) ds\Bigg]
- \mathbbm{1}(\tau_t^{(1)} > \eta_t^{t+\delta})\int_t^{\eta_t^{t+\delta}}\Theta_{t,s}\phi(s,X_s)ds\Bigg]\nonumber\\
&\qquad =\mathbb{E}_\lambda^{t,{\bf x}} \left[\mathbbm{1}(\tau_t^{(1)} >\eta_t^{t+\delta})\left(\Theta_{t,\eta_t^{t+\delta}}w_m(\eta_t^{t+\delta},X_{\eta_t^{t+\delta}}) - \int_t^{\eta_t^{t+\delta}}\Theta_{t,s}\phi(s,X_s)ds\right)\right],\label{case 1 theorem 3.3}
\end{align}
where we have applied the Markov property, the tower rule and the representation \eqref{wmrecur2}.

Next, consider the case where there is at one jump on the interval $(t,\eta^{t+\delta}_t]$, that is, the first jump time $\tau^{(1)}$ after $t$ is before or at latest at the end time $\eta_t^{t+\delta}$.
It then holds that 
\begin{multline}
\mathbb{E}_\lambda^{t,{\bf x}} \left[\mathbbm{1}(\tau_t^{(1)} \leq \eta_t^{t+\delta}) \left(\mathbbm{1}(X_{\eta_t^T\land \tau_t^{(1)}}\in\overline{D})\Theta_{t, \eta_t^T \land \tau_t^{(1)} } w_{m-1} (\eta_t^T\land \tau_t^{(1)}, X_{\eta_t^T\land \tau_t^{(1)} }) 
+\mathbbm{1}(X_{\eta_t^T\land \tau_t^{(1)}}\notin\overline{D})\Theta_{t,\eta^T_t} \Psi(\eta^T_t, X_{\eta^T_t},X_{\eta^T_t-}) 
- \int_t^{\eta_t^T \land \tau_t^{(1)}} \Theta_{t,s} \phi(s,X_s) ds\right) \right] \\
=\mathbb{E}_\lambda^{t,{\bf x}} \left[\mathbbm{1}(\tau_t^{(1)} \leq \eta_t^{t+\delta}) \left(\mathbbm{1}(X_{\tau_t^{(1)}}\in\overline{D})\Theta_{t, \tau_t^{(1)} } w_{m-1} (\tau_t^{(1)}, X_{\tau_t^{(1)}}) 
+\mathbbm{1}(X_{\tau_t^{(1)}}\notin\overline{D})\Theta_{t,\tau_t^{(1)}} \Psi(\tau_t^{(1)}, X_{\tau_t^{(1)}},X_{\tau_t^{(1)} -})-\int_t^{ \tau_t^{(1)}} \Theta_{t,s} \phi(s,X_s) ds\right) \right],\label{case 2 theorem 3.3}
\end{multline}
due to $\eta^T_t=\tau_t^{(1)}$ on the event $\{\tau_t^{(1)}\leq \eta^T_t\}\cap\{X_{\tau_t^{(1)}}\not\in\overline{D}\}$.
By combining the two identities \eqref{case 1 theorem 3.3} and \eqref{case 2 theorem 3.3}, we obtain
\begin{multline}\label{prepde}
w_m(t,{\bf x}) = \mathbb{E}_\lambda^{t,{\bf x}} \Bigg[\mathbbm{1}(\tau_t^{(1)} >\eta_t^{t+\delta})
\Theta_{t,\eta_t^{t+\delta}}w_m(\eta_t^{t+\delta},X_{\eta_t^{t+\delta}})\\
+\mathbbm{1}(\tau_t^{(1)} \leq \eta_t^{t+\delta})\Theta_{t,\tau_t^{(1)}}\left(
\mathbbm{1}(X_{\tau_t^{(1)}}\in\overline{D})
w_{m-1}(\tau_t^{(1)},X_{\tau_t^{(1)}})
+\mathbbm{1}(X_{\tau_t^{(1)}}\not\in\overline{D})
\Psi(\tau_t^{(1)},X_{\tau_t^{(1)}},X_{\tau_t^{(1)}-})\right)- \int_t^{\eta_t^{t+\delta}\land\tau_t^{(1)}}\Theta_{t,s}\phi(s,X_s)ds\Bigg].
\end{multline}

Here, we first focus on the first term in the expectation \eqref{prepde}, that is, on the event $\{\tau_t^{(1)} >\eta^{t+\delta}_t\}$, where the underlying process starting at time $t$ has not made a jump until the endpoint $\eta^{t+\delta}_t$ of the interval of interest.
Hence, even when it exits from the domain, that is, $\{\eta^T_t\le t+\delta\}$, the exit is necessarily caused by the diffusion term in the boundary $\partial D$.
With the aid of the given regularity of $w_m$, it holds $\mathbb{P}_{\lambda}^{t,{\bf x}}$-$a.s.$ on the event $\{\tau_t^{(1)} >\eta^{t+\delta}_t\}$ that 
\[
\Theta_{t,\eta^{t+\delta}_t} w_m(\eta^{t+\delta}_t,X_{\eta^{t+\delta}_t})
=w_m(t,{\bf x})+\int_t^{\eta^{t+\delta}_t}\Theta_{t,s}
\left(-r(s,X_s)w_m(s,X_s)+ \partial_1 w_m(s,X_s) +\mathcal{L}_sw_m(s,X_s)\right)ds+\int_t^{\eta^{t+\delta}_t}\Theta_{t,s}\left\langle \nabla w_m(s,X_s),\sigma(s,X_s)dW_s\right\rangle,
\]
that is, the jump component is degenerate on the interval $(t,t+\delta]$.
Hence, it holds by taking the expectation $\mathbb{E}_{\lambda}^{t,{\bf x}}$ that
\begin{multline*}
 \mathbb{E}_{\lambda}^{t,{\bf x}}\left[\mathbbm{1}(\tau_t^{(1)} >\eta^{t+\delta}_t)\Theta_{t,\eta^{t+\delta}_t} w_m(\eta^{t+\delta}_t,X_{\eta^{t+\delta}_t})\right]-w_m(t,{\bf x})\\
  =\mathbb{E}_{\lambda}^{t,{\bf x}}\left[\mathbbm{1}(\tau_t^{(1)} >\eta^{t+\delta}_t)\int_t^{\eta^{t+\delta}_t}\Theta_{t,s}
\left(-r(s,X_s)w_m(s,X_s)+ \partial_1 w_m(s,X_s) +\mathcal{L}_sw_m(s,X_s)\right)ds\right]
-\mathbb{P}_{\lambda}^{t,{\bf x}}(\tau_t^{(1)} \le \eta^{t+\delta}_t)w_m(t,{\bf x}).
\end{multline*}
By dividing this identity by $\delta$ and taking a limit $\delta\to 0+$, we get
\begin{equation}\label{result wm}
 \frac{1}{\delta}\left(\mathbb{E}_{\lambda}^{t,{\bf x}}\left[\mathbbm{1}(\tau_t^{(1)} >\eta^{t+\delta}_t)\Theta_{t,\eta^{t+\delta}_t} w_m(\eta^{t+\delta}_t,X_{\eta^{t+\delta}_t})\right]-w_m(t,{\bf x})\right)
\to -(r(t,{\bf x})+\lambda(t,{\bf x}))w_m(t,{\bf x})+ \partial_1 w_m(t,{\bf x}) +\mathcal{L}_tw_m(t,{\bf x}),
\end{equation}
where we have applied $\mathbb{P}_{\lambda}^{t,{\bf x}}(\tau_t^{(1)} \le \eta^{t+\delta}_t)\sim \delta \lambda(t,{\bf x})$ and $\mathbbm{1}(\tau_t^{(1)} >\eta^{t+\delta}_t)\to 1$, $\mathbb{P}_{\lambda}^{t,{\bf x}}$-$a.s.$, as $\delta\to 0+$, due to $\eta^{t+\delta}_t=t+\delta$ eventually $\mathbb{P}_{\lambda}^{t,{\bf x}}$-$a.s.$ 
By the same reasoning as for \eqref{result wm}, we obtain 
\begin{equation}\label{result phi}
\frac{1}{\delta}\mathbb{E}_{\lambda}^{t,{\bf x}}\left[\int_t^{\eta_t^{t+\delta}\land\tau_t^{(1)}}\Theta_{t,s}\phi(s,X_s)ds\right]\to \phi(t,{\bf x}).
\end{equation}

Next, for the second term in the expectation \eqref{prepde}, consider the event $\{\tau_t^{(1)}\leq \eta_t^{t+\delta}\}\cap  \{X_{\tau_t^{(1)}}\in\overline{D}\}$, that is, the first jump occurs before (or right at) $\eta_t^{t+\delta}$ and does not bring the trajectory strictly out of the closure $\overline{D}$.
In the other words, it holds true that the trajectory remains in the closure $\overline{D}$ over the interval $[t,\tau_t^{(1)}\land \eta_t^{t+\delta}]$, that is, $\mathbbm{1}(X_s\in\overline{D})=1$ for $s\in [t,\tau_t^{(1)}\land \eta_t^{t+\delta}]$, $\mathbb{P}_{\lambda}^{t,{\bf x}}$-$a.s.$
Thus, this indicator function does not break down smoothness for the Ito formula.
If $w_{m-1}$ is as smooth as imposed, we thus have
\begin{multline*}
    \mathbb{E}_{\lambda}^{t,{\bf x}}\left[
    \mathbbm{1}(X_{\tau_t^{(1)}\land \eta_t^{t+\delta}}\in\overline{D}) \Theta_{t,\tau_t^{(1)}\land \eta_t^{t+\delta}}w_{m-1}(\tau_t^{(1)}\land \eta_t^{t+\delta} ,X_{\tau_t^{(1)}\land \eta_t^{t+\delta}})\right]-\mathbbm{1}(x\in\overline{D})w_{m-1}(t,{\bf x}) \\
     =\mathbb{E}_{\lambda}^{t,{\bf x}}\Bigg[
    \int_t^{\tau_t^{(1)}\land \eta_t^{t+\delta}}\mathbbm{1}(X_s\in\overline{D})\Theta_{t,s}
    \Bigg(-(r(s,X_s)+\lambda(s,X_s))w_{m-1}(s,X_s)+\partial_1w_{m-1}(s,X_s)+\mathcal{L}_sw_{m-1}(s,X_s)   \\ 
    +\int_{\{{\bf z}\in\mathbb{R}_0^d:X_s+{\bf z}\in\overline{D}\}}w_{m-1}(t,X_s+{\bf z})\lambda(t,X_s)\nu(d{\bf z};X_s)    \Bigg)ds\Bigg],
\end{multline*}
where we have applied Lemma \ref{intermediate lemma} with $\kappa_t=(t+\delta)\land \tau_t^{(1)}$.
This, with the aid of the convergence \eqref{result wm} with $m$ replaced by $m-1$, yields the desired result 
\begin{equation}\label{result G}
\begin{aligned}
&\frac{1}{\delta}\mathbb{E}_{\lambda}^{t,{\bf x}}\left[
    \mathbbm{1}(\tau_t^{(1)}\leq \eta_t^{t+\delta})\mathbbm{1}(X_{\tau_t^{(1)}\land \eta_t^{t+\delta}}\in\overline{D}) \Theta_{t,\tau_t^{(1)}\land \eta_t^{t+\delta}}w_{m-1}(\tau_t^{(1)}\land \eta_t^{t+\delta} ,X_{\tau_t^{(1)}\land \eta_t^{t+\delta}})\right] \\
    &\qquad = \frac{1}{\delta}\left(\mathbb{E}_{\lambda}^{t,{\bf x}}\left[
    \mathbbm{1}(X_{\tau_t^{(1)}\land \eta_t^{t+\delta}}\in\overline{D}) \Theta_{t,\tau_t^{(1)}\land \eta_t^{t+\delta}}w_{m-1}(\tau_t^{(1)}\land \eta_t^{t+\delta} ,X_{\tau_t^{(1)}\land \eta_t^{t+\delta}})\right]-w_{m-1}(t,{\bf x})\right)\\
    &\qquad\qquad -\frac{1}{\delta}\left(\mathbb{E}_{\lambda}^{t,{\bf x}}\left[
    \mathbbm{1}(\tau_t^{(1)}> \eta_t^{t+\delta})\mathbbm{1}(X_{\tau_t^{(1)}\land \eta_t^{t+\delta}}\in\overline{D}) \Theta_{t,\tau_t^{(1)}\land \eta_t^{t+\delta}}w_{m-1}(\tau_t^{(1)}\land \eta_t^{t+\delta},X_{\tau_t^{(1)}\land \eta_t^{t+\delta}})\right]-w_{m-1}(t,{\bf x})\right)\\
    &\qquad \to \int_{\{{\bf z}\in\mathbb{R}_0^d:{\bf x}+{\bf z}\in\overline{D}\}}w_{m-1}(t,{\bf x}+{\bf z})\lambda(t,{\bf x})\nu(d{\bf z};{\bf x}) = G_{m-1}(t,{\bf x}), 
\end{aligned}\end{equation}
as $\delta\to 0+$.

Finally, for the third term in the expectation \eqref{prepde}, consider the event $\{\tau_t^{(1)}\leq \eta_t^{t+\delta}\}\cap \{X_{\tau_t^{(1)}}\not\in\overline{D}\}$, that is, the first jump occurs before or at $\eta_t^{t+\delta}$ and brings the trajectory strictly out of the closure $\overline{D}$.
Hence, the first jump time $\tau_t^{(1)}$ is nothing but the first exit time before, or at latest, the endpoint of the interval $[t,\eta_t^{t+\delta}]$ of current interest, that is, $\tau_t^{(1)}=\eta_t^T(\le t+\delta)$.
We thus have $\mathbbm{1}(X_{\cdot}\not\in \overline{D})=1$ on the singleton $\{\tau_t^{(1)}\}$ (which is smaller than or equal to $\eta_t^{t+\delta}$) and $\mathbbm{1}(X_{\cdot}\not\in \overline{D})=0$ on $[t,\tau_t^{(1)})(=[t,\eta_t^{t+\delta}]\setminus \{\tau_t^{(1)}\})$.
Conversely, if $\mathbbm{1}(X_s\not\in \overline{D})=0$ for all $s\in [t,\tau_t^{(1)}\land \eta_t^{t+\delta}]$, then no jump is allowed to occur before or at $\eta_t^{t+\delta}$ to bring the trajectory strictly out of the closure $\overline{D}$.
Therefore, it holds that 
\begin{align*}
    \mathbb{E}_{\lambda}^{t,{\bf x}}\left[
    \mathbbm{1}(\tau_t^{(1)}\leq \eta_t^{t+\delta})\mathbbm{1}(X_{\tau_t^{(1)}}\not\in\overline{D}) \Theta_{t,\tau_t^{(1)}}\Psi(\tau_t^{(1)},X_{\tau_t^{(1)}-} ,X_{\tau_t^{(1)}})\right]
    =&\mathbb{E}_{\lambda}^{t,{\bf x}}\left[
    \int_t^{\tau_t^{(1)}\land\eta_t^{t+\delta}}\int_{\mathbb{R}_0^d}\Theta_{t,s} \Psi(s,X_{s-}+{\bf z},X_{s-})\mathbbm{1}(X_{s-}+{\bf z}\notin \overline{D}) \mu(d{\bf z}, ds;X_{s-})\right]\\
    =&\mathbb{E}_{\lambda}^{t,{\bf x}}\left[
    \int_t^{\tau_t^{(1)}\land\eta_t^{t+\delta}}
    \int_{\{{\bf z}\in \mathbb{R}_0^d:\,X_s+{\bf z}\notin \overline{D}\}
}
    \Theta_{t,s}\Psi(s,X_s+{\bf z},X_s)\lambda(s,X_s)\nu(d{\bf z};X_s) ds\right],
\end{align*}
with the aid of Lemma \ref{intermediate lemma}, again with $\kappa_t=(t+\delta)\land \tau_t^{(1)}$.
By dividing the identity by $\delta$ and taking $\delta\to 0+$, we get
\begin{equation}\label{result H}
 \frac{1}{\delta}\mathbb{E}_{\lambda}^{t,{\bf x}}\left[
    \mathbbm{1}(\tau_t^{(1)}\leq \eta_t^{t+\delta})\mathbbm{1}(X_{\tau_t^{(1)}}\not\in\overline{D}) \Theta_{t,\tau_t^{(1)}}\Psi(\tau_t^{(1)},X_{\tau_t^{(1)}} ,X_{\tau_t^{(1)}-})\right]\to 
    \int_{\{{\bf z}\in \mathbb{R}_0^d: {\bf x}+{\bf z}\notin \overline{D}\}}
    \Psi(t,{\bf x}+{\bf z},{\bf x})\lambda(t,{\bf x})\nu(d{\bf z};{\bf x})=H(t,{\bf x}),
\end{equation}
since $\eta_t^{t+\delta}=t+\delta$ eventually $\mathbb{P}_{\lambda}^{t,{\bf x}}$-$a.s.$, due to a finite jump intensity, as for \eqref{result wm}.
By rearranging \eqref{prepde} and combining the results \eqref{result wm}, \eqref{result phi}, \eqref{result G} and \eqref{result H}, we obtain the desired partial differential equation \eqref{wmpde}.

To derive the initial condition, observe that for every ${\bf x}\in D$,
\[
w_m(t,{\bf x})=
\mathbb{E}_{\lambda}^{t,{\bf x}} \left[\mathbbm{1} (X_{T}\in \overline{D}) \Theta_{T,T} w_{m-1}\left(T, X_{T}\right) 
+ \mathbbm{1} (X_{T}\not\in \overline{D} ) \Theta_{T,T} \Psi(\eta_T^T, X_{\eta^T_T},X_{\eta^T_T-})
-\int_T^{T} \Theta_{T,s} \phi(s, X_s) ds\right]=w_{m-1}(t,{\bf x}),
\]
based on the representation \eqref{wmrecur2}, due to $\eta_T>T$, $\tau_T^{(1)}>T$ and $X_T={\bf x}$, $\mathbb{P}_{\lambda}^{t,{\bf x}}$-$a.s$.
Hence, by induction, we get $w_m(t,{\bf x})=w_0(t,{\bf x})=g({\bf x})$ for all ${\bf x}\in D$, again due to $\eta_T>T$, $\mathbb{P}_{\lambda}^{t,{\bf x}}$-$a.s.$ in \eqref{w0def1}.

\noindent{\bf (b)}
Due to \eqref{w0def1}, we have $w_0(T,{\bf x})=g({\bf x})$ and the identity
\begin{equation}\label{w0identity}
\mathbbm{1} (X_{\eta_t^T}\in\overline{D}) 
\Theta_{t,\eta_t^T} w_0(\eta_t^T,X_{\eta_t^T})=
\mathbbm{1}(\eta_t>T)\Theta_{t,T}g(X_T)+ \mathbbm{1}(\eta_t\leq T)\mathbbm{1}(X_{\eta_t}\in\partial D)\Theta_{t,\eta_t}\Psi(\eta_t,X_{\eta_t},X_{\eta_t}),\quad \mathbb{P}_0^{t,{\bf x}}\text{-}a.s.,
\end{equation}
where $X_{\eta_t^T}\in\overline{D}$ with probability one.
Thus, the representation \eqref{w0def1} can be rewritten as 
\[
\mathbb{E}_{0}^{t,{\bf x}} \left[\Theta_{t,\eta_t^T}w_0(\eta_t^T,X_{\eta_t^T})\right]-w_0(t,{\bf x})=
\mathbb{E}_{0}^{t,{\bf x}} \left[\int_t^{\eta_t^T} \Theta_{t,s} \phi(s, X_s) ds\right].
\]
Equating this with the Dynkin formula 
\[
 \mathbb{E}_{0}^{t,{\bf x}} \left[\Theta_{t,\eta_t^T}w_0(\eta_t^T,X_{\eta_t^T})\right]=w_0(t,{\bf x})+
\mathbb{E}_{0}^{t,{\bf x}} \left[\int_t^{\eta_t^T} \Theta_{t,s} \left(-r(s,X_s) w_0(s,X_s)+\partial_1 w_0(s,X_s) + \mathcal{L}_s w_0(s,X_s)\right) ds\right],
\]
under the smoothness conditions yields the partial differential equation
\begin{equation}\label{w0pde}
\partial_1 w_0(t,{\bf x}) + \mathcal{L}_t w_0(t,{\bf x}) = r(t,{\bf x}) w_0(t,{\bf x}) + \phi(t,{\bf x}),
\end{equation}
for almost every $(t,{\bf x})\in [0,T]\times \overline{D}$.
By rearranging \eqref{w0pde} and then applying it to the representation \eqref{udef}, we get the following identity 
\begin{align*} 
   u(t,{\bf x}) &=\mathbb{E}_{\lambda}^{t,{\bf x}} \left[\mathbbm{1} (\eta_t > T) \Theta_{t,T} g(X_T) + \mathbbm{1} (\eta_t \leq T) \Theta_{t,\eta_t} \Psi(\eta_t, X_{\eta_t},X_{\eta_t-}) -\int_t^{\eta_t^T} \Theta_{t,s} \phi(s, X_s) ds\right]\\
    &=\mathbb{E}_{\lambda}^{t,{\bf x}} \left[\mathbbm{1} (X_{\eta_t^T}\in\overline{D}) \Theta_{t,\eta_t^T} w_0(\eta_t^T,X_{\eta_t^T}) + \mathbbm{1} (X_{\eta_t^T}\not\in\overline{D}) \Theta_{t,\eta_t} \Psi(\eta_t, X_{\eta_t},X_{\eta_t-}) -\int_t^{\eta_t^T} \Theta_{t,s} \phi(s, X_s) ds\right]\\
    &=\mathbb{E}_{\lambda}^{t,{\bf x}} \Bigg[\mathbbm{1} (X_{\eta_t^T}\in\overline{D}) \Theta_{t,\eta_t^T} w_0(\eta_t^T,X_{\eta_t^T})-\int_t^{\eta_t^T}\Theta_{t,s}
    \left(\partial_1 w_0(s,X_s)+\mathcal{L}_sw_0(s,X_s) -r(s,X_s)w_0(s,X_s)+G_0(s,X_s)-\lambda(s,X_s)w_0(s,X_s)
    \right)ds\\
    &\qquad \qquad 
    + \mathbbm{1} (X_{\eta_t^T}\not\in\overline{D}) \Theta_{t,\eta_t} \Psi(\eta_t, X_{\eta_t},X_{\eta_t-})
    -\int_t^{\eta_t^T}\Theta_{t,s}H(s,X_s)ds + \int_t^{\eta_t^T}\Theta_{t,s}N_0(s,X_s)ds\Bigg]\\
    &=w_0(t,{\bf x})+\mathbb{E}_{\lambda}^{t,{\bf x}} \left[ \int_t^{\eta_t^T}\Theta_{t,s}N_0(s,X_s)ds\right],
\end{align*}
where the third equality holds by Lemma \ref{intermediate lemma} with $f=w_0$ and $\kappa_t=T$.
Thus, by taking essential supremum and infimum over the time-state domain, the functions $N_0^U(t)$ and $N_0^L(t)$ can be taken out of the expectation to yield the inequalities
\[
N_0^L(t)\xi(t,{\bf x};\lambda)\leq u(t,{\bf x}) -w_0(t,{\bf x})\leq N_0^U(t)\xi(t,{\bf x};\lambda),\quad (t,{\bf x})\in [0,T]\times \overline{D}.
\]
By setting $g\equiv 0$, $\Psi\equiv 0$ and $\phi\equiv -1$ in \eqref{udef} and \eqref{w0def1}, we get $H\equiv 0$, $M=N_0$, $w_0(t,{\bf x})=\xi(t,{\bf x};0)$, and $u(t,{\bf x})=\xi(t,{\bf x};\lambda)$.
Hence, the last inequalities reduce to $M^L(t)\xi(t,{\bf x};\lambda)\leq \xi(t,{\bf x};\lambda) -\xi(t,{\bf x};0)\leq M^U(t)\xi(t,{\bf x};\lambda)$.
By rearranging this inequality under the assumption $M^L(t)\le  M^U(t)<1$, we get
\begin{equation}\label{Mr inequality}
\frac{1}{1-M^L(t)}\xi(t,{\bf x};0)\leq \xi(t,{\bf x};\lambda)\leq 
\frac{1}{1-M^U(t)}\xi(t,{\bf x};0),\quad (t,{\bf x})\in [0,T]\times \overline{D}.
\end{equation}
In a similar manner, by rearranging \eqref{wmpde} and with \eqref{def of Nm}, we get the identity
\[
    \phi(t,{\bf x}) 
=\partial_1 w_m(t,{\bf x})+\mathcal{L}_tw_m(t,{\bf x}) -r(t,{\bf x})w_m(t,{\bf x})-\lambda(t,{\bf x})w_m(t,{\bf x}) + G_{m}(t,{\bf x})+H(t,{\bf x})  - N_m(t,{\bf x}),
\]
for almost every $(t,{\bf x})\in [0,T)\times D$, with which the representation \eqref{udef} yields
    $u(t,{\bf x}) 
    =w_m(t,{\bf x})
    +\mathbb{E}_\lambda^{t,{\bf x}}[\int_t^{\eta_t^T}\Theta_{t,s}N_m(s,X_s)ds]$,
where we have applied $w_m(T,{\bf x})=g({\bf x})$ for all ${\bf x}\in D$, and Lemma \ref{intermediate lemma}.
Thus, by taking the essential supremum $N_m^U(t)$ and infimum $N_m^L(t)$ out of the expectation and then further applying the inequalities \eqref{Mr inequality}, we obtain the desired result \eqref{boundswm0}.
\end{proof}

For the proofs of Theorems \ref{theorem thinning} and \ref{theorem uniform}, we recall the notation $\widetilde{\mathbb{P}}_{\widetilde{\lambda}}^{t,{\bf x}}$ and $\widetilde{\mathbb{E}}_{\widetilde{\lambda}}^{t,{\bf x}}$, which respectively denote the probability measure and the associated expectation under which the underlying process is governed by the stochastic differential equation process \eqref{underlying process with jumps of size zero} with $X_t={\bf x}$ almost surely. 
Moreover, we have denoted by $\{\widetilde{\tau}_t^{(m)}\}_{m\in\mathbb{N}}$ a sequence of jump times of the redefined underlying process \eqref{underlying process with jumps of size zero}, where each $\widetilde{\tau}_t^{(m)}$ indicates its $m$-th jump time, which may be of a zero-sized jump.
 
\begin{proof}[Sketch of proof of Theorem \ref{theorem thinning}]
To avoid overloading the paper with simple repetition, we do not go through the proofs of Theorems \ref{theorem convergence}, \ref{theorem wmrecurrsive}, \ref{theorem vm wm}, \ref{theorem back to E0} and \ref{theorem when smooth} by highlighting required amendments.
Instead, we summarize a few key points to be taken care of.
Under the probability measure $\widetilde{\mathbb{P}}_{\widetilde{\lambda}}^{t,{\bf x}}$, the underlying process \eqref{underlying process with jumps of size zero} makes jumps based on the constant rate $\widetilde{\lambda}$, irrespective of the time and state.
Each jump is then rejected (by setting its size to zero) by the time-state dependent probability $\widetilde{\nu}(\{0\};t,{\bf x})$.
The jump times, denoted by $\{\widetilde{\tau}_t^{(m)}\}_{m\in\mathbb{N}}$, including those of zero-sized jumps, can thus be represented as a cumulative sum $\widetilde{\tau}_t^{(m)}\stackrel{\mathcal{L}}{=}t+\sum_{k=1}^m E_k$, where $\{E_k\}_{k\in\mathbb{N}}$ is a sequence of iid exponential random variables with rate $\widetilde{\lambda}$, that is, $\widetilde{\mathbb{P}}_{\widetilde{\lambda}}^{t,{\bf x}}(\widetilde{\tau}_t^{(1)}\geq s)=\widetilde{\Lambda}_{t,s}$ for $s\ge t$.
With those in mind, the result can be derived in a similar manner to the previous results.
\end{proof}

\begin{proof}[Proof of Theorem \ref{theorem uniform}]
Let $(t,{\bf x})\in [0,T]\times \overline{D}$ and $m\in\mathbb{N}_0$.
By combining two equalities based on \eqref{lemma51a} with $w_0$ as the test function, one with $\kappa_t=T$ and the other with $\kappa_t=\tau_t^{(m)}$, we obtain 
\begin{align*}
    &\mathbb{E}_{\lambda}^{t,{\bf x}} \left[\mathbbm{1} (X_{\eta_t^T}\in\overline{D}) \Theta_{t,\eta_t^T} w_0(\eta_t^T,X_{\eta_t^T})
    -\mathbbm{1} (X_{\eta_t^T\land\tau_t^{(m)}}\in\overline{D}) \Theta_{t,\eta_t^T\land\tau_t^{(m)}} w_0(\eta_t^T\land\tau_t^{(m)},X_{\eta_t^T\land\tau_t^{(m)}})\right]\\
    &\qquad\qquad=
    \mathbb{E}_{\lambda}^{t,{\bf x}} \left[
    \int_{\eta_t^T\land\tau_t^{(m)}}^{\eta_t^T} \Theta_{t,s}\left(
    \partial_1 w_0(s,X_s) +\mathcal{L}_s w_0(s,X_s) -(r(s,X_s)+\lambda(s,X_s))w_0(s,X_s) + G_0(s,X_s)\right)ds\right]\\
    &\qquad\qquad=
    \mathbb{E}_{\lambda}^{t,{\bf x}} \left[
    \int_{\eta_t^T\land\tau_t^{(m)}}^{\eta_t^T} \Theta_{t,s}\left(
    \phi(s,X_s)-\lambda(s,X_s)w_0(s,X_s) + G_0(s,X_s)\right)ds\right],
\end{align*}
where the second equality is due to the identity \eqref{w0pde}.
In a similar manner, by combining two equalities based on \eqref{lemma51b}, one with $(f,\kappa_t)=(u,T)$ and the other with $(f,\kappa_t)=(w_m,\tau_t^{(m)})$, we obtain
\begin{align*}
    &\mathbb{E}_{\lambda}^{t,{\bf x}} \left[\mathbbm{1} (X_{\eta_t^T}\not\in\overline{D})\Theta_{t,\eta_t} \Psi(\eta_t, X_{\eta_t},X_{\eta_t-})- \mathbbm{1} (X_{\eta_t^T\land\tau_t^{(m)}}\not\in\overline{D}) \Theta_{t,\eta_t} \Psi(\eta_t, X_{\eta_t},X_{\eta_t-})\right]
    =\mathbb{E}_{\lambda}^{t,{\bf x}} \left[\int_{\eta_t^T\land\tau_t^{(m)}}^{\eta_t^T} \Theta_{t,s} H(s,X_s) ds\right],
\end{align*}
where we have applied $\int_{\{{\bf z}\in \mathbb{R}_0^d:\,X_{s-}+{\bf z}\notin \overline{D}\}}f(s,X_{s-}+{\bf z})\nu(d{\bf z};X_{s-})=\int_{\{{\bf z}\in \mathbb{R}_0^d:\,X_{s-}+{\bf z}\notin \overline{D}\}}\Psi(s,X_{s-}+{\bf z},X_{s-})\nu(d{\bf z};X_{s-})$ for $s\in (t,\eta_t^T]$ if $f=u$ and for $s\in (t,\eta_t^T\land \tau_t^{(m)}]$ if $f=w_m$, respectively, in accordance with \eqref{udef} and \eqref{wmdef2}, each of which is absolutely continuous with respect to the Lebesgue measure $ds$.
By combining those identities along with \eqref{udef} and \eqref{wmdef2}, we get
\begin{align}
  u(t,{\bf x}) - w_m(t,{\bf x})
  &=\mathbb{E}_{\lambda}^{t,{\bf x}} \Bigg[\mathbbm{1} (X_{\eta_t^T}\in\overline{D}) \Theta_{t,\eta_t^T} w_0(\eta_t^T,X_{\eta_t^T}) + \mathbbm{1} (X_{\eta_t^T}\not\in\overline{D})\Theta_{t,\eta_t} \Psi(\eta_t, X_{\eta_t},X_{\eta_t-}) -\int_{\eta_t^T \land \tau_t^{(m)}}^{\eta_t^T} \Theta_{t,s} \phi(s, X_s) ds\nonumber\\
  &\qquad\qquad -\mathbbm{1} (X_{\eta_t^T\land\tau_t^{(m)}}\in\overline{D}) \Theta_{t,\eta_t^T\land\tau_t^{(m)}} w_0(\eta_t^T\land\tau_t^{(m)},X_{\eta_t^T\land\tau_t^{(m)}}) - \mathbbm{1} (X_{\eta_t^T\land\tau_t^{(m)}}\not\in\overline{D}) \Theta_{t,\eta_t} \Psi(\eta_t, X_{\eta_t},X_{\eta_t-}) \Bigg]\nonumber\\
  &= \mathbb{E}_{\lambda}^{t,{\bf x}} \left[
    \int_{\eta_t^T \land \tau_t^{(m)}}^{\eta_t^T} \Theta_{t,s}\left(
    -\lambda(s,X_s)w_0(s,X_s) + G_0(s,X_s) + H(s,X_s)\right)ds
    \right]\nonumber\\
    &= \mathbb{E}_{\lambda}^{t,{\bf x}} \left[
    \mathbbm{1}(\tau_t^{(m)}\leq \eta_t^T)
    \int_{\tau_t^{(m)}}^{\eta_t^T} \Theta_{t,s}N_0(s,X_s)ds
    \right],\label{E1 N0}
\end{align}
where we have applied \eqref{w0identity} and \eqref{def of Nm}.
Hence, due to $\Theta\in (0,1]$ and $[\tau_t^{(m)},\eta_t^T]\subseteq [t,T]$, it holds that for $(t,{\bf x})\in [0,T]\times \overline{D}$, 
\begin{equation}\label{wminequal}
|u(t,{\bf x}) - w_m(t,{\bf x})|\leq  (T-t) 
\max\left\{|N_0^U(t)|,|N_0^L(t)|\right\}\mathbb{P}_{\lambda}^{t,{\bf x}} (\tau_t^{(m)}\leq T).
\end{equation}
As for the rightmost probability, letting $N_{t,T}$ denote the number of jumps on $[t,T]$ and $\lambda_0(t):=\sup_{(s,{\bf x})\in [t,T]\times \overline{D}}\lambda(s,{\bf x})$, it holds that for every $m\in \{\lceil \lambda_0(t)(T-t)\rceil,\cdots\}$,
\[
 \mathbb{P}_\lambda^{t,{\bf x}}(\tau_t^{(m)}\le T)=\sum_{k=m}^{+\infty}\mathbb{P}_\lambda^{t,{\bf x}}(N_{t,T}=k)\le \sum_{k=m}^{+\infty}e^{-\lambda_0(t)(T-t)}\frac{(\lambda_0(t)(T-t))^k}{k!},
\]
for all $(t,{\bf x})\in [0,T]\times \overline{D}$,
which is independent of the initial state ${\bf x}$.
Recalling $\lambda_0=\sup_{t\in[0,T]}\lambda_0 (t)(T-t)$, there exists a strictly positive  constant $c_m$ such that
\begin{equation}\label{probability vanishing}
 \sup_{(t,{\bf x})\in [0,T]\times \overline{D}}\mathbb{P}_\lambda^{t,{\bf x}}(\tau_t^{(m)}\le T)\lesssim \sup_{t\in [0,T]}\sum_{k=m}^{+\infty}e^{-\lambda_0(t)(T-t)}\frac{(\lambda_0(t)(T-t))^k}{k!}
 \lesssim \sum_{k=m}^{+\infty}e^{-\lambda_0}\frac{\lambda_0^k}{k!}=e^{-c_m}\frac{\lambda_0^m}{m!}\to 0,
\end{equation}
as $m\to +\infty$.
Here, the asymptotic inequality holds since $\lambda_0$ is finite and $e^{-x}x^k$ is increasing on $(0,k)$.
The last equality holds 
\[
 \frac{m!}{\lambda_0^m}\sum_{k=m}^{+\infty}e^{-\lambda_0}\frac{\lambda_0^k}{k!}=e^{-\lambda_0}\sum_{k=m}^{+\infty}\frac{\lambda_0^{k-m}}{k!/m!}\le e^{-\lambda_0} \sum_{k=m}^{+\infty}\frac{\lambda_0^{k-m}}{(k-m)!}=1,
\]
due to $\binom{k}{m}\ge 1$ with equality only when $m=0$ or $k=m$.
Combining the inequalities \eqref{wminequal} and \eqref{probability vanishing} yields the desired uniform convergence.
Therefore, by the identity \eqref{wm vm}, the triangle inequality yields 
\[
|u(t,{\bf x}) - v_m(t,{\bf x})|\le |u(t,{\bf x}) - w_m(t,{\bf x})| +\left|\mathbb{E}_\lambda^{t,{\bf x}}\left[\mathbbm{1}(\tau_t^{(m)}< \eta_t^T)\Theta_{t,\tau_t^{(m)}}w_0(\tau_t^{(m)},X_{\tau_t^{(m)}})\right] \right|,
\]
where the last term has the same decay rate as \eqref{E1 N0} at most, due to the boundedness of the function $w_0$.
Finally, for each fixed $t\in[0,T]$, it holds by \eqref{def of Nm} that for every $(s,{\bf x})\in [t,T]\times\overline{D}$,
\[
\left|N_m(s,{\bf x})\right|=\left|G_m(s,{\bf x})-G_{m-1}(s,{\bf x})\right|
\leq  
\int_{\{{\bf z}\in\mathbb{R}_0^d:{\bf x}+{\bf z}\in\overline{D}\}}
\left|w_m(s,{\bf x}+{\bf z}) -w_{m-1}(s,{\bf x}+{\bf z})\right|\lambda(s,{\bf x})\nu(d{\bf z}; {\bf x}),
\]
which yields $N_m^L(t) \to 0$ and $N_m^U(t) \to 0$, as $m\to +\infty$, after taking supremum of $(s,{\bf x})$ over $[t,T]\times\overline{D}$, due to Theorem \ref{theorem uniform} and $\int_{\mathbb{R}_0^d}\nu(d{\bf z}; {\bf x})=1$ for all ${\bf x}$.
The results for $\widetilde{w}_m$ can be derived in a similar manner to \eqref{wminequal} under the probability measure $\widetilde{\mathbb{P}}_{\widetilde{\lambda}}^{(t,{\bf x})}$, yet with the condition $|N_0^U(0)|+|N_0^L(0)|<+\infty$ in common.
\end{proof}

\section*{Disclosure statement}

No potential conflict of interest was reported by the author(s).

\section*{Funding}

This work was partially supported by JSPS Grants-in-Aid for Scientific Research 20K22301 and 21K03347.

\bibliographystyle{abbrv}
\bibliography{combib.bib}

\end{document}